\documentclass[a4paper,12pt]{article}
\usepackage{amsmath}
\usepackage{amsthm}
\usepackage{amscd}
\usepackage{eufrak}

\theoremstyle{definition}

\theoremstyle{remark}

\numberwithin{equation}{section}

\newcommand{\mib}[1]{\mbox{\boldmath $#1$}}
\def\0{{\bf 0}}
\def\a{{\bf a}}

\def\f{{\mib{f}}}

\def\x{{\mib{x}}}
\def\y{{\mib{y}}}
\def\z{{\mib{z}}}
\def\h{{\mib{h}}}

\def\B{{\mib{B}}}
\def\C{{\bf{C}}}

\def\H{{\mathcal{H}}(N)}
\def\J{{\mib{J}}}
\def\M{{\mib{M}}}
\def\N{{\bf{N}}}
\def\R{{\bf{R}}}
\def\S{{\mathcal{S}}(N)}
\def\W{{\bf W}_{N}}
\def\X{{\mib{X}}}
\def\Y{{\mib{Y}}}

\def\mbK{{\bf K}}
\def\mbA{{\bf A}}
\def\Det{{\rm Det}}
\def\vtheta{\mib{\theta}}
\def\vlambda{{\mib{\lambda}}}
\def\vkappa{{\mib{\kappa}}}
\def\vomega{{\mib{\omega}}}
\def\vchi{{\mib{\chi}}}

\def\dH{{\mathcal{U}}(dH)}
\def\dS{{\mathcal{V}}(dS)}

\def\cNt{\widetilde{\mathcal{N}}_N}
\def\cNnk{{\mathcal{N}_N^{(\nu,\kappa)}}}

\def\WA{{\bf W}_{N}^{\rm A}}
\def\WC{{\bf W}_{N}^{\rm C}}
\def\WD{{\bf W}_{N}^{\rm D}}

\def\cR{{\rho}_N}

\begin{document}

\title{Noncolliding processes, 
matrix-valued processes \\
and determinantal processes
\footnote{
to be published in {\it Sugaku Expositions} (AMS)
}}

\author{
Makoto Katori
\footnote{
Department of Physics,
Faculty of Science and Engineering,
Chuo University, 
Kasuga, Bunkyo-ku, Tokyo 112-8551, Japan;
e-mail: katori@phys.chuo-u.ac.jp
}
and 
Hideki Tanemura
\footnote{
Department of Mathematics and Informatics,
Faculty of Science, Chiba University, 
1-33 Yayoi-cho, Inage-ku, Chiba 263-8522, Japan;
e-mail: tanemura@math.s.chiba-u.ac.jp
}}
%%%%%%%%%%%%%%%%%%%%%%%%%%%%%%%%
\date{23 March 2011}
%%%%%%%%%%%%%%%%%%%%%%%%%%%%%%%%
\pagestyle{plain}
\maketitle
\begin{abstract}
A noncolliding diffusion process is a conditional process of 
$N$ independent one-dimensional diffusion processes
such that the particles never collide with 
each other.
This process realizes an interacting particle system
with long-ranged strong repulsive forces acting between
any pair of particles.
When the individual diffusion process
is a one-dimensional Brownian motion, the noncolliding process
is equivalent in distribution with 
the eigenvalue process of an $N \times N$
Hermitian-matrix-valued process, which we call
Dyson's model.
For any deterministic initial configuration of $N$ particles, 
distribution of particle positions 
of the noncolliding Brownian motion on the real line
at any fixed time $t >0$ is a determinantal point process.
We can prove that 
the process is determinantal in the sense that
the multi-time correlation function
for any chosen series of times,
which determines joint distributions at these times,
is also represented by a determinant.
We study the asymptotic behavior of the system, 
when the number of Brownian motions $N$ 
in the system tends to infinity.
This problem is concerned with 
the random matrix theory on 
the asymptotics of eigenvalue distributions, 
when the matrix size becomes infinity.
In the present paper, we introduce a variety of
noncolliding diffusion processes
by generalizing the noncolliding Brownian motion,
some of which are temporally inhomogeneous.
We report the results of our research project
to construct and study finite and infinite particle systems
with long-ranged strong interactions
realized by noncolliding processes. 
%%%%%%%%%%%%%%%%%%%%%%%%%%%%%
\\
\begin{small}
\end{small}
%%%%%%%%%%%%%%%%%%%%%%%%%%%%%%%%%%

\end{abstract}
\maketitle

%%%%%%%%%%%%%%%%%%%%%%%%%%%%%%%%%%%%%%%%%%%%%%%%%%%%%%%%%%%%
%%%%%%%%%%Section 1%%%%%%%%%%%%%%%%%%%%%%%%%%%%%%%%%%%%%%%%%
%%%%%%%%%%%%%%%%%%%%%%%%%%%%%%%%%%%%%%%%%%%%%%%%%%%%%%%%%%%%
\section{Introduction}
%%%%%%%%%%%%%%%%%%%%%%%%%%%%%%%%%%%%%%%%%%%%%%%%%%%%%%%%%%%%
%%%%%%%%%%%%%%%%%%%%%%%%%%%%%%%%%%%%%%%%%%%%%%%%%%%%%%%%%%%%
%%%%%%%%%%%%%%%%%%%%%%%%%%%%%%%%%%%%%%%%%%%%%%%%%%%%%%%%%%%%
In a system of $N$ independent one-dimensional diffusion processes,
if we impose a condition such that the particles never collide with 
each other, then we obtain an interacting particle system
with long-ranged strong repulsive forces acting between
any pair of particles.
We call such a system a {\bf noncolliding diffusion process}.
In 1962 Dyson \cite{Dys62b} 
showed that, when the individual diffusion process
is a one-dimensional Brownian motion, the obtained 
noncolliding process,
the {\bf noncolliding Brownian motion},
is related to a {\bf matrix-valued process}.
He introduced a Hermitian-matrix-valued process having 
Brownian motions as its diagonal elements, 
and complex Brownian motions as off-diagonal elements.
The size of the matrix is supposed to be $N\times N$.
By virtue of the Hermitian property, all eigenvalues of the matrix
are real, and Dyson derived a system of 
$N$-simultaneous stochastic differential equations
for the process of $N$ eigenvalues on the real line $\R$.
In the present paper we call this stochastic process
of eigenvalues {\bf Dyson's model}.
(Strictly speaking, it is a special case
of Dyson's Brownian motion models 
with the parameter $\beta=2$ as explained below.)
If we regard each eigenvalue as a particle position
in one dimension, Dyson's model is 
considered to be a one-dimensional
system of interacting Brownian motions.
Dyson showed that this system is nothing but 
the noncolliding Brownian motion \cite{Bi95,Gra99}.

A probability distribution on the space of particle configurations
is called a {\bf determinantal point process} 
or a {\bf Fermion point process}, 
if its correlation functions are generally represented 
by determinants \cite{ST03,Sos00,HKPV09}.
The noncolliding Brownian motion provides us
examples of determinantal point processes:
for any deterministic initial configuration of $N$ particles, 
distribution of particle positions on $\R$ 
at any fixed time $t >0$
is a determinantal point process \cite{KT10}.
Moreover, by using the method developed by Eynard and Mehta
for multi-layer random matrix models 
\cite{EM98, Mehta},
we can show that the multi-time correlation functions
for any chosen series of times,
which determine joint distributions at these times,
are also represented by determinants \cite{KNT04,KT07b,KT09,KT10}.
In the present paper we call such a stochastic process 
that any multi-time correlation function
is given by a determinant a {\bf determinantal process}
\cite{KT07b}.

We study the asymptotic behavior of the system, 
when the number of Brownian motions $N$ 
in the system tends to infinity.
Since, as explained above, the noncolliding Brownian motion 
can be realized by the eigenvalue process and the correlation
functions are expressed by determinants of matrices,
this problem is concerned with the asymptotics
of eigenvalue distributions, 
when the matrix size becomes infinity.
The latter problem is one of the main topics of the random
matrix theory \cite{Mehta}. In other words, 
our research project reported 
in this paper is to construct infinite particle systems
with long-ranged strong interactions 
by applying the results of recent development of 
the random matrix theory \cite{KNT04,KT04,KT07a,KT07b,KT09,KT10}.

In the present paper, we introduce a variety of
noncolliding diffusion processes
by generalizing the noncolliding Brownian motion.
In Section 2 first we explain basic properties of diffusion processes
treated in this paper,
such as Brownian motions, Brownian bridges, 
absorbing Brownian motions, Bessel processes, Bessel bridges, 
and generalized meanders.
The transition probability density of a noncolliding diffusion process
is expressed by a determinant of a matrix,
each element of which is
the transition probability density of the individual diffusion process
in one dimension (the Karlin-McGregor formula). 
This formula provides a useful tool for us
to analyze noncolliding diffusion processes.
In Section 3 we state the Karlin-McGregor formula 
and present basic properties of
noncolliding diffusion processes.
When such a Hermitian-matrix-valued process is given that 
its elements are one-dimensional diffusion processes,
it will be a fundamental and interesting problem to determine 
a system of stochastic differential equations for
eigenvalue process of the given matrix-valued process.
Bru's theorem \cite{Bru89, Bru91} 
and its generalization \cite{KT03b,KT04} give
answers to this problem.
In Section 4 we give a generalized version of Bru's theorem 
and show its applications. 
The determinantal structures of multi-time correlation functions
of noncolliding processes are explained in Section 5.
There asymptotics in $N \to \infty$ are also discussed
\cite{KNT04,KT07a,KT07b,KT09,KT10}.

When we impose the noncolliding condition on 
a finite time-interval $(0,T), T \in (0, \infty)$
instead of an infinite time-interval $(0,\infty)$, 
the noncolliding diffusion processes become temporally inhomogeneous,
even if individual one-dimensional diffusion processes 
are temporally homogeneous.
In Section 6 we discuss these temporally inhomogeneous 
noncolliding processes.
These processes are not determinantal any more,
and make a new family of processes, which we call
{\bf Pfaffian processes} \cite{KNT04,KT07a}.
In the last section, Section 7, we list up the topics,
which are related to noncolliding processes, but
can not be discussed here.

%%%%%%%%%%%%%%%%%%%%%%%%%%%%%%%%%%%%%%%%%%%%%%%%%%%%%%%%%%%%
%%%%%%%%%%Section 2%%%%%%%%%%%%%%%%%%%%%%%%%%%%%%%%%%%%%%%%%
%%%%%%%%%%%%%%%%%%%%%%%%%%%%%%%%%%%%%%%%%%%%%%%%%%%%%%%%%%%%
\section{Brownian motion and its conditional processes}
%%%%%%%%%%%%%%%%%%%%%%%%%%%%%%%%%%%%%%%%%%%%%%%%%%%%%%%%%%%%
%%%%%%%%%%%%%%%%%%%%%%%%%%%%%%%%%%%%%%%%%%%%%%%%%%%%%%%%%%%
%%%%%%%%%%%%%%%%%%%%%%%%%%%%%%%%%%%%%%%%%%%%%%%%%%%%%%%%%%%

Let $(\Omega, \mathcal{F}, P)$ be a probability space.
The stochastic process called (one-dimensional or linear)
{\bf Brownian motion}, 
$\{ B(t,\omega)\}_{t\in [0,\infty)}$, 
satisfies the following conditions :
\begin{enumerate}
\item $B(0,\omega)=0$ with probability one.

\item For any fixed $\omega\in \Omega$, 
$B(t,\omega)$ is a real continuous function of $t$.
(This property is expressed by saying that
$B(t)$ has a {\bf continuous path}.)

\item For any sequence of times, $t_0\equiv 0<t_1<\cdots <t_M,
M=1,2,\dots$,
the increments $\{B(t_i)-B(t_{i-1})\}_{i=1,2,\dots,M}$ are independent,
and distribution of each increment 
is normal with mean zero and variance $t_i-t_{i-1}$.
\end{enumerate}
Then, the probability that the Brownian motion is observed 
in the interval $[a_i, b_i] \subset \R$ at time $t_i$ 
for each $i=1,2, \dots, M$,
$P(B(t_i)\in [a_i, b_i], i=1,2,\dots, M)$,
is given by
$$
\int_{a_1}^{b_1} dx_1 \int_{a_2}^{b_2} dx_2 \cdots 
\int_{a_M}^{b_M}dx_M
\prod_{i=1}^M G(t_i-t_{i-1},x_i-x_{i-1}),
$$
where $x_0\equiv 0$ and 
$$
G(t,x) = \frac{1}{\sqrt{2 \pi t}} 
\exp \left( - \frac{x^2}{2t}  \right), 
\quad t>0, \quad x \in\R.
$$
The integral kernel $G(s,x; t,y) \equiv G(t-s, y-x)$ is called 
the {\bf transition probability density function} of the Brownian motion.
For any fixed $s\ge 0$, under the condition that $B(s)$ is given,
$B(u), u\le s$ and $B(t), t>s$ are independent.
This property is called a {\bf Markov property}.
A positive random variable $\tau$ is called a {\bf Markov time}, 
if the event $\{\tau \le u\}$ is determined by the behavior of the process
until time $u$ and independent of the behavior
of the process after time $u$.
The first time that a Brownian motion visits a given domain $D$,
which is called the hitting time of $D$, 
is an example of a Markov time.
In the definition of Markov property mentioned above,
if a deterministic time $s$ is replaced by
a Markov time $\tau$, then we obtain a stronger property,
called a {\bf strong Markov property}.
In general, a stochastic process, which has a
strong Markov property and has a continuous path
almost surely, is called a diffusion process. 
See \cite{Wat75, RY98} for instance.
In the case that the transition probability density function
$G(s,x;t,y)$ 
does not depend on times $t$ and $s$ themselves
but only depends on the difference $t-s$, 
a Markov process is said to be {\bf temporally homogeneous}.
In this case we write the
transition probability density function as $G(t-s,y|x)$ 
instead of $G(s,x; t,y)$ to clarify its homogeneity in time 
in this paper.
The Brownian motion is an example of a 
temporally homogeneous diffusion process.
(It is also spatially homogeneous.) 

For $d\in \N \equiv \{1,2, \dots\}$, 
using independent one-dimensional 
Brownian motions $B_1(t),B_2(t), \dots, B_d(t), t \geq 0$, 
a $d$-dimensional Brownian motion is defined
by a vector-valued diffusion process
$\B(t)=(B_1(t),B_2(t), \dots, B_d(t)), t \geq 0$.

We want to consider the Brownian motion 
under the condition that it visits
the origin at a given time $T>0$.
Since the probability that this condition is satisfied is zero,
we first consider the Brownian motion under another condition 
such that it visits some point in an interval 
$(-\varepsilon, \varepsilon)$ at time $T$,
$\varepsilon > 0$, 
and then define the original conditional process by taking
the limit $\varepsilon \downarrow 0$.
The transition probability density function 
obtained in this limit is
$$
G^T(s,x;t,y) = \frac{G(T-t,0|y)G(t-s,y|x)}{G(T-s,0|x)},
\quad 0 \leq s < t \leq T, \quad x, y \in \R. 
$$
It is a temporally inhomogeneous diffusion process. 
We call this process a {\bf Brownian bridge} of duration $T$
and denote it by $\beta^{T}(t)$, $t\in [0,T]$.

Although one-dimensional Brownian motion can visit any point of $\R$,
we consider the Brownian motion conditioned to stay positive forever.
This conditional process $Y(t)$, $t\in [0, \infty)$ 
is temporally homogeneous process
with transition probability density function 
$G^{(1/2)}(t,y|x)$;
\begin{eqnarray}
G^{(1/2)}(t,y|x)
&=&\frac{y}{x} \Big\{ G(t,y|x) - G(t,-y|x) \Big\}, 
\quad t >0, \quad x>0, \, y \geq 0,
\label{eqn:Bessel1}
\\
\nonumber
G^{(1/2)}(t,y|0)
&=&\frac{2}{t} y^2 G(t,y|0), \quad t > 0, \quad y \geq 0.
\end{eqnarray}
The distance of a three-dimensional Brownian motion from the origin, 
$( B_1(t)^2+B_2(t)^2+B_3(t)^2 )^{1/2}, t > 0$
has exactly the same transition probability 
density as (\ref{eqn:Bessel1}),
and is called the {\bf three-dimensional Bessel process}.
In other words, the three-dimensional Bessel process 
$Y(t), t \in [0, \infty)$ 
has two different representations, 
`the representation by a Brownian motion conditioned 
to stay positive' and 
`the representation by 
a radial part of the three-dimensional Brownian motion'.
We also note that
$Y(t)$ solves the following stochastic differential equation \cite{RY98},
$$
\label{bessel}
Y(t)=B(t)+ \int_{0}^{t} \frac{1}{Y(s)} ds, \quad t>0.
$$

Consider the Brownian motion $X(t), t\in [0,T]$ 
under the condition that it stays positive during 
a finite time-interval $(0,T]$,
with $T\in (0,\infty)$.
Then the conditional process is a temporally inhomogeneous diffusion process 
with the transition probability density function $G^{(1/2,1)}_{T}(s,x;t,y)$;
\begin{eqnarray}
&&G^{(1/2,1)}_{T}(s,x;t,y)
= \frac{h(T-t,y)}{h(T-s,x)}
\Big\{ G(t-s,y|x) - G(t-s,-y|x) \Big\},
\label{eqn:GT2}
\\\nonumber
&&\qquad\qquad\qquad\qquad\qquad\qquad\qquad\qquad
\, 0 \leq s < t \leq T, \, x > 0, \, y \geq 0,
\nonumber\\
&&G^{(1/2,1)}_{T}(0,0;t,y)
= \frac{\sqrt{2\pi T}}{t} h(T-t,y) y G(t,y|0),
\quad t \in (0, T], \quad y \geq 0,
\nonumber
\end{eqnarray}
where $h(s,x)$, $x > 0, \, s >0 $ is
the probability that the Brownian motion starting from $x>0$
stays positive during the time-interval $[0,s]$.
This conditional process is called a {\bf Brownian meander}.
Using two independent Brownian motions $B_1(t), B_2(t)$ 
and a Brownian bridge $\beta^{T}(t)$ of duration $T$, 
which is independent of $B_1(t)$ and $B_2(t)$, 
we define a one-dimensional diffusion process by
$(B_1(t)^2+B_2(t)^2+\beta^{T}(t)^2)^{1/2}$,
$t\in [0,T]$. We can prove that
this process is identified with the Brownian meander.
That is, 
the Brownian meander has also two different representations, 
`the representation by a Brownian motion conditioned
to stay positive during
a finite time-interval $(0, T]$' and 
`the representation by 
a radial part of the three-dimensional 
diffusion process $(B_1(t), B_2(t), \beta^{T}(t)),
t \in (0, T]$' \cite{Yor92}.

By comparing (\ref{eqn:Bessel1}) with (\ref{eqn:GT2}),
we see that the distributions of
the three-dimensional Bessel process $Y(t)$ 
and the Brownian meander $X(t)$, both starting from the origin, 
are absolutely continuous and satisfy
\begin{equation}
\label{Imhof:original}
P(X(\cdot)\in dw) 
= \sqrt{\frac{\pi T}{2}}\frac{1}{w(T)}P(Y(\cdot)\in dw).
\end{equation}
The equality (\ref{Imhof:original}) 
is called {\bf Imhof's relation} \cite{Imhof84}.

The Brownian motion, which is killed at the origin, 
is called an {\bf absorbing Brownian motion} 
in the domain $(0,\infty)$.
Let $\widehat{G}(t-s,y|x)$ be
the transition probability density of this process.
It is the density of the Brownian motion at time $t$,
which starts from $x>0$ at time $s \,(< t)$, 
restricted on the event that it
stays positive in the time-interval $[s,t]$.
The {\bf reflection principle} of Brownian motion gives
$$
\widehat{G}(t-s,y|x)=G(t-s,y|x) - G(t-s,-y|x).
$$
The first formula in (\ref{eqn:Bessel1}) means that
the transformation of the transition probability density
$\widehat{G}(t,y|x)$, 
given by
$(y/x)\widehat{G}(t,y|x)$, 
is identified with the transition probability density
$G^{(1/2)}(t, y|x)$ of 
the three-dimensional Bessel process.
It implies that 
the three-dimensional Bessel process $Y(t)$ 
is the Doob {\bf $\h$-transformation}
of the absorbing Brownian motion in the domain $(0, \infty)$.

When $d \in \N$, the distance 
of the $d$-dimensional Brownian motion 
from the origin,
$( B_1(t)^2+B_2(t)^2+\cdots +B_d(t)^2 )^{1/2}$,
defines a one-dimensional diffusion process, which we call
the {\bf $d$-dimensional Bessel process}.
The Bessel process can be extended to the cases with
all positive real values of $d$ as follows.
With a parameter $\nu \in (-1,\infty)$, 
the transition probability density function of the
$2(\nu +1)$-dimensional Bessel process $Y^{(\nu)}(t)$, 
is given by
\begin{eqnarray*}
\label{eqn:Bessel1_G}
G^{(\nu)}(t,y|x)
&=&\frac{y^{\nu+1}}{x^{\nu}}\frac{1}{t} 
\exp \left( - \frac{x^2+y^2}{2t} \right)
I_{\nu} \left( \frac{x y}{t} \right), 
\quad t>0, \quad x>0, \ y\ge 0,
\\
G^{(\nu)}(t,y|0)
&=&\frac{y^{2\nu +1}}{2^{\nu} \Gamma (\nu +1)t^{\nu+1}} 
\exp \left( - \frac{y^2}{2t} \right),
\qquad t>0, \quad y\ge 0,
\end{eqnarray*}
where $\Gamma (z)$ is the Gamma function and $I_{\nu}(z)$ 
is the modified Bessel function with parameter $\nu$ \cite{RY98}.
The behavior of the Bessel process depends on 
the dimension $d$ (the parameter $\nu=(d-2)/2$).
When $d$ is greater than or equal to $2$ ($\nu \ge 0$),
the process has the origin as a transient point,
and when $d$ is less than $2$ ($-1<\nu <0$),
it has the origin as a recurrent point.
Moreover, if and only if $d$ is greater than or 
equal to $1$ ($\nu \ge -1/2$),
it is a semi-martingale \cite{RY98}.

Yor \cite{Yor92} introduced a family of diffusion processes 
with two parameters 
$(\nu, \kappa)$, $\nu \in (-1,\infty)$, $\kappa \in (0, 2(\nu+1) )$,
which includes the Brownian meander as a special case
$(\nu, \kappa)=(1/2,1)$, 
and he called each member of the family a 
{\bf generalized meander}.
The generalized meander $X^{(\nu,\kappa)}(t)$, 
$\nu \in (-1,\infty)$, $\kappa \in (0, 2(\nu+1) )$,
is the diffusion process with the transition probability density 
\begin{eqnarray*}
\label{eqn:GT1_G}
&&G^{(\nu,\kappa)}_{T}(s,x;t,y)
= \frac{h^{(\nu,\kappa)}_{T}(t,y)}{h^{(\nu,\kappa)}_{T}(s,x)}
G^{(\nu)}(t-s,y|x),
\quad 0 \le s < t \le T, \quad x > 0, \ y \ge 0,
\\
\label{eqn:GT2_G}
&&G^{(\nu,\kappa)}_{T}(0,0;t,y)
= \frac{\Gamma (\nu+1)(2T)^{\kappa/2}}{\Gamma (\nu+1-\kappa/2)}
h^{(\nu,\kappa)}_{T}(t,y) G^{(\nu)}(t,y|0),
\quad t\in (0,T], \quad y \ge 0,
\end{eqnarray*}
where $h^{(\nu,\kappa)}_{T}(t,x)
=\int_0^\infty dy \ G^{(\nu)}(T-t,y|x)y^{-\kappa}$,
$x \ge 0, \, t\in (0,T]$ \cite{Yor92}. 
Then Imhof's relation (\ref{Imhof:original}) 
between the three-dimensional Bessel process
$Y(t)=Y^{(1/2)}(t)$ and the Brownian meander 
$X(t)=X^{(1/2,1)}(t)$ 
is generalized as
$$
P(X^{(\nu,\kappa)}(\cdot)\in dw) 
= \frac{\Gamma (\nu+1)}{\Gamma (\nu+1-\kappa/2)}
\left( \frac{\sqrt{2T}}{w(T)}\right)^\kappa P(Y^{(\nu)}(\cdot)\in dw)
$$
for the $2(\nu+1)$-dimensional Bessel process 
$Y^{(\nu)}(t)$ and
the generalized meander $X^{(\nu, \kappa)}(t)$.
We remark that, though the parameter $\kappa$
of generalized meanders is in $(0, 2(\nu+1))$,
we can discuss the cases $\kappa=0$ and
$\kappa=2(\nu+1)$.
The former corresponds to the Bessel processes,
and the latter the Bessel bridges,
which are the conditional Bessel processes
to arrive at the origin at a fixed time $T >0$.

%%%%%%%%%%%%%%%%%%%%%%%%%%%%%%%%%%%%%%%%%%%%%%%%%%%%%%%%%%%%
%%%%%%%%%%Section 3%%%%%%%%%%%%%%%%%%%%%%%%%%%%%%%%%%%%%%%%%
%%%%%%%%%%%%%%%%%%%%%%%%%%%%%%%%%%%%%%%%%%%%%%%%%%%%%%%%%%%%
\section{Noncolliding diffusion processes}
%%%%%%%%%%%%%%%%%%%%%%%%%%%%%%%%%%%%%%%%%%%%%%%%%%%%%%%%%%%%
%%%%%%%%%%%%%%%%%%%%%%%%%%%%%%%%%%%%%%%%%%%%%%%%%%%%%%%%%%%%
%%%%%%%%%%%%%%%%%%%%%%%%%%%%%%%%%%%%%%%%%%%%%%%%%%%%%%%%%%%%

%%%%%%%%%%%%%%%%%%%%%%%%%%%%%%%%%%%%%%%%%%%%%%%%%%%
\subsection{Karlin-McGregor formula}
%%%%%%%%%%%%%%%%%%%%%%%%%%%%%%%%%%%%%%%%%%%%%%%%%%%

In order to analyze noncolliding diffusion processes,
it is useful to represent the transition probability 
density functions by means of determinants.
The representation is called the {\bf Karlin-McGregor formula}
\cite{KM59b} in probability theory,
and the {\bf Lindstr\"om-Gessel-Viennot formula}
\cite{Li73,GV85,Ste90} in combinatorics. 
It is also regarded as a stochastic-process version of the
{\bf Slater determinant}, which originally expresses
a many-body wave function of free Fermion particles
in quantum mechanics \cite{Suzuki,KT07b}.

%%%%%%%%%%%%%%%%%%%%%%%%%%%%%%%%%%%%%%%%%%%%%%%%%%%
%%% Karlin-McGregor formula %%%%%%%%%%%%%%%%%%%%%%
%%%%%%%%%%%%%%%%%%%%%%%%%%%%%%%%%%%%%%%%%%%%%%%%%%%
\vskip 3mm
\noindent
{\large\bf [Karlin-McGregor formula]}
{\rm (\cite{KM59b,Li73,GV85})}
\label{th21}
\; 
Let $G(s,x; t,y)$ be the transition probability density function
of a one-dimensional diffusion process.
On the line $\R$ set $N$ starting points $x_i$, $i=1,2,\dots,N$ and
$N$ terminal points $y_i$, $i=1,2,\dots,N$ with
$x_1< x_2 < \cdots < x_N$ and $y_1< y_2 < \cdots <y_N$,
respectively.
The transition probability density function 
of the system of $N$ diffusion processes restricted
on the event that they never collide with each other 
during the time-interval $[s,t]$ is given by
$$
G_0(s,{\x}; t,{\y})
=\det_{1 \leq i,j \leq N}
\Big(G(s,x_i; t,y_j) \Big).
$$
%%%%%%%%%%%%%%%%%%%%%%%%%%%%%%%%%%%%%%%%%%%%%%%%%%%
\vskip 3mm
When $N=2$, this claim is essentially equivalent to 
the reflection principle for a Brownian motion.
That is, this formula can be regarded as a generalization of
the reflection principle \cite{KT06}.

%%%%%%%%%%%%%%%%%%%%%%%%%%%%%%%%%%%%%%%%%%%%%%%%%%%
\subsection{Noncolliding Brownian motions
in a finite and an infinite time-intervals}
%%%%%%%%%%%%%%%%%%%%%%%%%%%%%%%%%%%%%%%%%%%%%%%%%%%
%%%%%%%%%%%%%%%%%%%%%%%%%%%%%%%%%%%%%%%%%%%%%%%%%%%

Let $\WA$ be a subset of $\R^N$ defined by
$\WA=\{\x \in \R^N: x_1<x_2< \, \cdots \, <x_N\}$,
which is called the Weyl chamber of type $A_{N-1}$
in representation theory \cite{FH91}. 
The transition probability density 
of the absorbing Brownian motion in $\WA$,
that is, the density function of an $N$-dimensional Brownian motion
at time $t$, which starts from $\x \in \WA$ at time 0, 
restricted on the event that it stays in $\WA$
during the time-interval $[0,t]$, is represented by
\begin{equation*}
f_N(t,\y|\x)=\det_{1 \leq i,j \leq N}\Big(G(t, y_j|x_i)\Big)
\end{equation*}
by the Karlin-McGregor formula.
Then the probability that the Brownian motion stays in $\WA$ 
until time $t$ is
\begin{equation*}
\mathcal{N}_N(t,\x)=\int_{\WA} f_N(t,\y|\x) d\y.
\end{equation*}

Now we consider the {\bf noncolliding Brownian motion
in a finite time-interval} $t \in (0,T]$,
$\X(t)=(X_1(t), X_2(t), \dots X_N(t))$.
The transition probability density  of the process
denoted by
$g_{N,T}(s,\x;t,\y)$ is the conditional density 
of $N$ Brownian motions at time $t$, which
started from the points 
$\x=(x_1, x_2, \dots, x_N) \in \WA$ at time $s \, (< t)$, 
under the condition that they never collide 
with each other in the time-interval $[s,t]$.
It is given by
\begin{equation}
\label{eqn:g(x,y)}
g_{N,T}(s,\x;t,\y) = \frac{\mathcal{N}_N(T-t, \y)}{\mathcal{N}_N(T-s, \x)}
f_N(t-s,\y|\x), \quad
0 \leq s < t \leq T,
\quad \x, \y \in \WA.
\end{equation}
The transition probability density in the case that
all $N$ particles start from the origin 
will be defined by taking the limit 
$\x \to {\bf 0} \equiv (0,0,\dots,0)$
in (\ref{eqn:g(x,y)}).
Since both of the numerator and the denominator in (\ref{eqn:g(x,y)})
tend to $0$ as $\x\to {\bf 0}$,
we have to know the asymptotic behavior of
$f_{N}(t,\y |\x)$ and $\mathcal{N}_{N}(t,\x)$
in $|\x|/\sqrt{t} \to 0$.
Performing bilinear expansions \cite{KT03a,KT04,KT07b} 
with respect to the multivariate symmetric functions
called the {\bf Schur functions} \cite{FH91,Ful97},
we have obtained
\begin{eqnarray}
f_{N}(t,\y |\x) &\sim&
\frac{t^{-N(N+1)/4}}{C_{1}(N)}
h_{N} \left( \frac{\x}{\sqrt{t}} \right) h_{N}(\y)
\exp \left( - \frac{|\y|^2}{2t} \right),
\label{eqn:asym}\\
\mathcal{N}_{N}(t,\x) &\sim&
\frac{C_{2}(N)}{C_{1}(N)}
h_{N}\left( \frac{\x}{\sqrt{t}} \right),
\qquad \frac{|\x|}{\sqrt{t}} \to 0.
\nonumber
\end{eqnarray}
Here
$h_{N}(\x)$ is the $N \times N$
{\bf Vandermonde determinant}, which is equal to
the product of differences 
of variables $x_1, x_2, \dots, x_N$, 
\begin{equation*}
h_{N}(\x)
= \det_{1 \leq i, j \leq N} \Big(
x_{j}^{i-1} \Big)
= \prod_{1 \leq i<j \leq N}(x_{j}-x_{i}),
\label{eqn:Vandermonde}
\end{equation*}
and 
$C_1(N)=(2 \pi)^{N/2} \prod_{i=1}^{N} \Gamma(i)$,
$C_2(N)=2^{N/2} \prod_{i=1}^{N} \Gamma(i/2)$.
By using (\ref{eqn:asym}) 
we obtain the transition probability density function
of the noncolliding Brownian motion,
when all $N$ particles start from the origin 
({\it i.e.} $\X(0)={\bf 0}$) as
\begin{eqnarray}
\label{eqn:g(0,y)}
&&g_{N,T}(0,\0;t,\y)=\frac{T^{N(N-1)/4}t^{-N^2/2}}{C_2(N)}
\mathcal{N}_N(T-t,\y) h_N(\y) 
\exp\left(-\frac{|\y|^2}{2t}\right),
\\
&&\qquad\qquad\qquad\qquad\qquad\qquad\qquad\qquad
\quad t \in (0, T], 
\quad \y \in \WA.
\nonumber
\end{eqnarray}
As illustrated in the left picture of Fig.1,
in this case the $N$ Brownian motions starting
from the origin at time $t=0$ 
rapidly separate from each other
to avoid collision.

%%%%%%%%%%%%%%%%%%%%%%%%%%%%%%%%%%%%%%%%%%
%%%%%%%%%%FIG 1%%%%%%%%%%%%%%%%%%%5
%%%%%%%%%%%%%%%%%%%%%%%%%%%%%%%%%%%%%%%%%%%
%\vskip 3mm
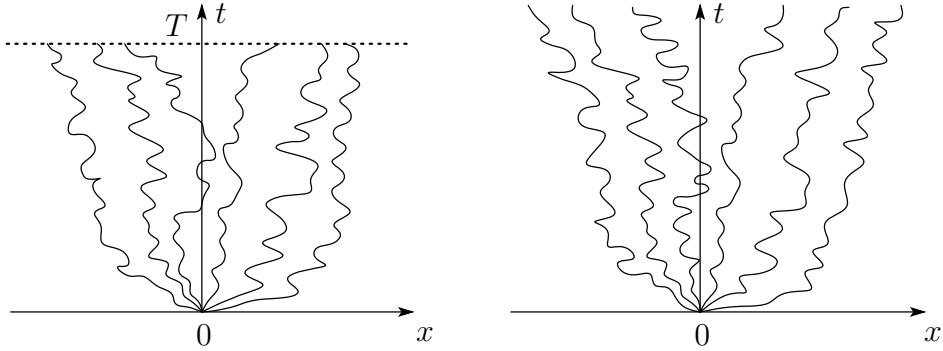
\begin{figure}[ch]
\begin{center}
%\input fig_ronsetu_1.tex
%WinTpicVersion3.08
\unitlength 0.1in
\begin{picture}( 47.5000, 16.8500)(  3.9000,-18.8500)
% VECTOR 2 0 3 0
% 2 3000 1880 5140 1880
% 
\special{pn 8}%
\special{pa 3000 1880}%
\special{pa 5140 1880}%
\special{fp}%
\special{sh 1}%
\special{pa 5140 1880}%
\special{pa 5074 1860}%
\special{pa 5088 1880}%
\special{pa 5074 1900}%
\special{pa 5140 1880}%
\special{fp}%
% STR 2 0 3 0
% 3 3950 1950 3950 2050 2 0
% $0$
\put(39.5000,-20.5000){\makebox(0,0)[lb]{$0$}}%
% STR 2 0 3 0
% 3 5140 1930 5140 2030 2 0
% $x$
\put(51.4000,-20.3000){\makebox(0,0)[lb]{$x$}}%
% STR 2 0 3 0
% 3 4050 270 4050 370 2 0
% $t$
\put(40.5000,-3.7000){\makebox(0,0)[lb]{$t$}}%
% VECTOR 2 0 3 0
% 2 3980 1880 3980 280
% 
\special{pn 8}%
\special{pa 3980 1880}%
\special{pa 3980 280}%
\special{fp}%
\special{sh 1}%
\special{pa 3980 280}%
\special{pa 3960 348}%
\special{pa 3980 334}%
\special{pa 4000 348}%
\special{pa 3980 280}%
\special{fp}%
% SPLINE 2 0 3 0
% 34 3980 1880 4030 1800 4070 1740 4030 1660 4070 1640 4100 1570 4070 1500 4100 1470 4100 1420 4060 1360 4070 1340 4120 1330 4150 1290 4170 1250 4190 1210 4140 1140 4100 1080 4100 1000 4150 1000 4170 950 4140 910 4180 880 4180 800 4150 760 4200 710 4220 670 4180 620 4200 570 4260 540 4410 480 4310 400 4380 380 4410 330 4410 280
% 
\special{pn 8}%
\special{pa 3980 1880}%
\special{pa 3994 1852}%
\special{pa 4010 1826}%
\special{pa 4032 1800}%
\special{pa 4056 1774}%
\special{pa 4070 1746}%
\special{pa 4060 1716}%
\special{pa 4036 1684}%
\special{pa 4030 1660}%
\special{pa 4058 1646}%
\special{pa 4088 1628}%
\special{pa 4102 1596}%
\special{pa 4098 1562}%
\special{pa 4078 1532}%
\special{pa 4070 1504}%
\special{pa 4088 1482}%
\special{pa 4106 1454}%
\special{pa 4102 1422}%
\special{pa 4078 1394}%
\special{pa 4062 1370}%
\special{pa 4070 1340}%
\special{pa 4102 1334}%
\special{pa 4130 1324}%
\special{pa 4148 1296}%
\special{pa 4162 1266}%
\special{pa 4178 1238}%
\special{pa 4190 1208}%
\special{pa 4180 1182}%
\special{pa 4156 1156}%
\special{pa 4132 1132}%
\special{pa 4114 1108}%
\special{pa 4100 1076}%
\special{pa 4090 1038}%
\special{pa 4094 1006}%
\special{pa 4118 1000}%
\special{pa 4152 1000}%
\special{pa 4172 974}%
\special{pa 4164 938}%
\special{pa 4142 914}%
\special{pa 4156 896}%
\special{pa 4186 874}%
\special{pa 4196 842}%
\special{pa 4186 808}%
\special{pa 4162 780}%
\special{pa 4150 756}%
\special{pa 4170 736}%
\special{pa 4200 712}%
\special{pa 4220 682}%
\special{pa 4212 656}%
\special{pa 4188 632}%
\special{pa 4180 604}%
\special{pa 4198 574}%
\special{pa 4222 554}%
\special{pa 4254 542}%
\special{pa 4292 534}%
\special{pa 4336 524}%
\special{pa 4376 514}%
\special{pa 4404 502}%
\special{pa 4412 486}%
\special{pa 4394 464}%
\special{pa 4360 442}%
\special{pa 4326 420}%
\special{pa 4310 404}%
\special{pa 4326 394}%
\special{pa 4362 388}%
\special{pa 4392 370}%
\special{pa 4408 342}%
\special{pa 4412 310}%
\special{pa 4410 280}%
\special{sp}%
% SPLINE 2 0 3 0
% 29 3980 1880 3920 1820 3880 1790 3880 1740 3790 1680 3810 1630 3700 1590 3760 1510 3700 1460 3710 1430 3670 1370 3740 1330 3660 1230 3800 1170 3690 1090 3780 1040 3600 910 3710 850 3600 770 3640 760 3590 710 3660 670 3420 560 3470 520 3440 480 3440 480 3430 400 3320 340 3320 280
% 
\special{pn 8}%
\special{pa 3980 1880}%
\special{pa 3962 1856}%
\special{pa 3938 1834}%
\special{pa 3910 1814}%
\special{pa 3884 1796}%
\special{pa 3880 1766}%
\special{pa 3878 1734}%
\special{pa 3848 1716}%
\special{pa 3812 1702}%
\special{pa 3790 1684}%
\special{pa 3802 1656}%
\special{pa 3810 1626}%
\special{pa 3784 1612}%
\special{pa 3742 1604}%
\special{pa 3706 1598}%
\special{pa 3704 1580}%
\special{pa 3730 1554}%
\special{pa 3758 1526}%
\special{pa 3754 1502}%
\special{pa 3722 1484}%
\special{pa 3700 1462}%
\special{pa 3710 1432}%
\special{pa 3696 1402}%
\special{pa 3672 1376}%
\special{pa 3682 1358}%
\special{pa 3718 1344}%
\special{pa 3742 1328}%
\special{pa 3732 1304}%
\special{pa 3702 1276}%
\special{pa 3670 1250}%
\special{pa 3662 1228}%
\special{pa 3682 1210}%
\special{pa 3722 1198}%
\special{pa 3764 1188}%
\special{pa 3796 1178}%
\special{pa 3798 1164}%
\special{pa 3770 1146}%
\special{pa 3732 1126}%
\special{pa 3700 1106}%
\special{pa 3692 1088}%
\special{pa 3716 1072}%
\special{pa 3752 1056}%
\special{pa 3780 1042}%
\special{pa 3780 1024}%
\special{pa 3754 1006}%
\special{pa 3714 986}%
\special{pa 3670 966}%
\special{pa 3630 946}%
\special{pa 3604 928}%
\special{pa 3600 910}%
\special{pa 3626 896}%
\special{pa 3666 882}%
\special{pa 3700 866}%
\special{pa 3710 848}%
\special{pa 3684 828}%
\special{pa 3642 806}%
\special{pa 3606 786}%
\special{pa 3598 772}%
\special{pa 3630 764}%
\special{pa 3636 752}%
\special{pa 3604 726}%
\special{pa 3592 704}%
\special{pa 3620 686}%
\special{pa 3656 672}%
\special{pa 3668 660}%
\special{pa 3648 650}%
\special{pa 3610 638}%
\special{pa 3560 626}%
\special{pa 3506 612}%
\special{pa 3458 600}%
\special{pa 3426 584}%
\special{pa 3416 568}%
\special{pa 3436 550}%
\special{pa 3466 530}%
\special{pa 3464 508}%
\special{pa 3440 480}%
\special{pa 3434 448}%
\special{pa 3434 414}%
\special{pa 3422 390}%
\special{pa 3394 376}%
\special{pa 3358 364}%
\special{pa 3328 350}%
\special{pa 3316 324}%
\special{pa 3318 292}%
\special{pa 3320 280}%
\special{sp}%
% SPLINE 2 0 3 0
% 26 3980 1880 4190 1850 4340 1830 4420 1780 4560 1750 4520 1670 4620 1640 4560 1560 4640 1530 4620 1430 4740 1330 4730 1190 4830 1090 4720 1000 4850 960 4790 890 4880 800 4850 710 4880 670 4810 610 4950 570 4870 480 4950 420 4880 370 5020 350 5020 280
% 
\special{pn 8}%
\special{pa 3980 1880}%
\special{pa 4012 1874}%
\special{pa 4042 1868}%
\special{pa 4074 1862}%
\special{pa 4104 1858}%
\special{pa 4136 1854}%
\special{pa 4170 1852}%
\special{pa 4204 1850}%
\special{pa 4238 1850}%
\special{pa 4270 1848}%
\special{pa 4302 1844}%
\special{pa 4332 1836}%
\special{pa 4356 1820}%
\special{pa 4380 1800}%
\special{pa 4410 1784}%
\special{pa 4448 1776}%
\special{pa 4492 1772}%
\special{pa 4530 1770}%
\special{pa 4556 1764}%
\special{pa 4558 1746}%
\special{pa 4540 1716}%
\special{pa 4520 1684}%
\special{pa 4528 1664}%
\special{pa 4564 1656}%
\special{pa 4604 1650}%
\special{pa 4620 1638}%
\special{pa 4602 1612}%
\special{pa 4572 1584}%
\special{pa 4560 1560}%
\special{pa 4584 1548}%
\special{pa 4622 1540}%
\special{pa 4646 1522}%
\special{pa 4640 1492}%
\special{pa 4626 1458}%
\special{pa 4622 1426}%
\special{pa 4636 1400}%
\special{pa 4666 1382}%
\special{pa 4698 1364}%
\special{pa 4728 1346}%
\special{pa 4744 1322}%
\special{pa 4742 1292}%
\special{pa 4734 1260}%
\special{pa 4726 1224}%
\special{pa 4730 1194}%
\special{pa 4750 1166}%
\special{pa 4782 1142}%
\special{pa 4812 1122}%
\special{pa 4830 1102}%
\special{pa 4824 1080}%
\special{pa 4794 1056}%
\special{pa 4758 1034}%
\special{pa 4728 1014}%
\special{pa 4720 1000}%
\special{pa 4744 990}%
\special{pa 4786 984}%
\special{pa 4828 976}%
\special{pa 4850 964}%
\special{pa 4840 944}%
\special{pa 4812 918}%
\special{pa 4790 892}%
\special{pa 4800 870}%
\special{pa 4830 848}%
\special{pa 4862 826}%
\special{pa 4880 802}%
\special{pa 4872 772}%
\special{pa 4854 742}%
\special{pa 4850 712}%
\special{pa 4872 684}%
\special{pa 4878 660}%
\special{pa 4850 638}%
\special{pa 4818 618}%
\special{pa 4814 604}%
\special{pa 4842 596}%
\special{pa 4888 590}%
\special{pa 4930 584}%
\special{pa 4950 570}%
\special{pa 4938 550}%
\special{pa 4906 526}%
\special{pa 4876 500}%
\special{pa 4872 476}%
\special{pa 4902 456}%
\special{pa 4936 438}%
\special{pa 4950 418}%
\special{pa 4926 398}%
\special{pa 4892 378}%
\special{pa 4880 368}%
\special{pa 4902 366}%
\special{pa 4942 366}%
\special{pa 4988 364}%
\special{pa 5020 350}%
\special{pa 5030 324}%
\special{pa 5022 288}%
\special{pa 5020 280}%
\special{sp}%
% SPLINE 2 0 3 0
% 33 3980 1880 3890 1850 3850 1800 3810 1750 3740 1740 3730 1690 3590 1680 3550 1620 3640 1600 3590 1540 3500 1540 3500 1450 3430 1410 3540 1370 3460 1250 3510 1200 3510 1150 3390 1110 3450 1020 3490 950 3350 890 3370 820 3460 800 3400 720 3280 670 3190 630 3330 610 3300 520 3220 490 3150 450 3200 420 3150 360 3090 280
% 
\special{pn 8}%
\special{pa 3980 1880}%
\special{pa 3948 1874}%
\special{pa 3916 1864}%
\special{pa 3890 1850}%
\special{pa 3868 1828}%
\special{pa 3850 1800}%
\special{pa 3834 1770}%
\special{pa 3810 1750}%
\special{pa 3774 1748}%
\special{pa 3744 1744}%
\special{pa 3734 1716}%
\special{pa 3728 1686}%
\special{pa 3708 1676}%
\special{pa 3674 1682}%
\special{pa 3636 1688}%
\special{pa 3596 1682}%
\special{pa 3564 1660}%
\special{pa 3548 1630}%
\special{pa 3564 1614}%
\special{pa 3602 1610}%
\special{pa 3636 1604}%
\special{pa 3640 1584}%
\special{pa 3618 1556}%
\special{pa 3580 1538}%
\special{pa 3540 1540}%
\special{pa 3510 1544}%
\special{pa 3498 1532}%
\special{pa 3504 1500}%
\special{pa 3506 1462}%
\special{pa 3482 1434}%
\special{pa 3448 1418}%
\special{pa 3430 1410}%
\special{pa 3450 1402}%
\special{pa 3492 1392}%
\special{pa 3530 1380}%
\special{pa 3542 1360}%
\special{pa 3528 1334}%
\special{pa 3498 1306}%
\special{pa 3470 1278}%
\special{pa 3460 1252}%
\special{pa 3478 1230}%
\special{pa 3506 1208}%
\special{pa 3518 1176}%
\special{pa 3506 1146}%
\special{pa 3476 1130}%
\special{pa 3438 1124}%
\special{pa 3406 1118}%
\special{pa 3388 1106}%
\special{pa 3396 1086}%
\special{pa 3418 1058}%
\special{pa 3446 1026}%
\special{pa 3474 994}%
\special{pa 3490 964}%
\special{pa 3486 944}%
\special{pa 3462 934}%
\special{pa 3426 926}%
\special{pa 3388 918}%
\special{pa 3358 900}%
\special{pa 3346 872}%
\special{pa 3354 838}%
\special{pa 3380 816}%
\special{pa 3418 808}%
\special{pa 3450 806}%
\special{pa 3462 794}%
\special{pa 3450 768}%
\special{pa 3422 740}%
\special{pa 3392 714}%
\special{pa 3364 698}%
\special{pa 3338 688}%
\special{pa 3308 678}%
\special{pa 3272 668}%
\special{pa 3232 652}%
\special{pa 3200 638}%
\special{pa 3192 630}%
\special{pa 3212 628}%
\special{pa 3250 630}%
\special{pa 3292 628}%
\special{pa 3326 616}%
\special{pa 3336 590}%
\special{pa 3328 556}%
\special{pa 3306 526}%
\special{pa 3282 508}%
\special{pa 3252 498}%
\special{pa 3216 490}%
\special{pa 3176 474}%
\special{pa 3150 458}%
\special{pa 3164 442}%
\special{pa 3196 426}%
\special{pa 3200 404}%
\special{pa 3172 378}%
\special{pa 3142 352}%
\special{pa 3118 326}%
\special{pa 3102 300}%
\special{pa 3090 280}%
\special{sp}%
% VECTOR 2 0 3 0
% 2 410 1880 2490 1880
% 
\special{pn 8}%
\special{pa 410 1880}%
\special{pa 2490 1880}%
\special{fp}%
\special{sh 1}%
\special{pa 2490 1880}%
\special{pa 2424 1860}%
\special{pa 2438 1880}%
\special{pa 2424 1900}%
\special{pa 2490 1880}%
\special{fp}%
% STR 2 0 3 0
% 3 1210 340 1210 440 2 0
% $T$
\put(12.1000,-4.4000){\makebox(0,0)[lb]{$T$}}%
% STR 2 0 3 0
% 3 1370 1950 1370 2050 2 0
% $0$
\put(13.7000,-20.5000){\makebox(0,0)[lb]{$0$}}%
% STR 2 0 3 0
% 3 2510 1930 2510 2030 2 0
% $x$
\put(25.1000,-20.3000){\makebox(0,0)[lb]{$x$}}%
% STR 2 0 3 0
% 3 1470 270 1470 370 2 0
% $t$
\put(14.7000,-3.7000){\makebox(0,0)[lb]{$t$}}%
% LINE 1 2 3 0
% 2 390 480 2460 480
% 
\special{pn 13}%
\special{pa 390 480}%
\special{pa 2460 480}%
\special{dt 0.045}%
% VECTOR 2 0 3 0
% 2 1400 1880 1400 280
% 
\special{pn 8}%
\special{pa 1400 1880}%
\special{pa 1400 280}%
\special{fp}%
\special{sh 1}%
\special{pa 1400 280}%
\special{pa 1380 348}%
\special{pa 1400 334}%
\special{pa 1420 348}%
\special{pa 1400 280}%
\special{fp}%
% SPLINE 2 0 3 0
% 30 1400 1880 1450 1800 1490 1740 1450 1660 1490 1640 1520 1570 1490 1500 1520 1470 1520 1420 1480 1360 1490 1340 1540 1330 1570 1290 1590 1250 1610 1210 1560 1140 1520 1080 1520 1000 1570 1000 1590 950 1560 910 1600 880 1600 800 1570 760 1620 710 1640 670 1600 620 1620 570 1680 540 1800 480
% 
\special{pn 8}%
\special{pa 1400 1880}%
\special{pa 1414 1852}%
\special{pa 1430 1826}%
\special{pa 1452 1800}%
\special{pa 1476 1774}%
\special{pa 1490 1746}%
\special{pa 1480 1716}%
\special{pa 1456 1684}%
\special{pa 1450 1660}%
\special{pa 1478 1646}%
\special{pa 1508 1628}%
\special{pa 1522 1596}%
\special{pa 1518 1562}%
\special{pa 1498 1532}%
\special{pa 1490 1504}%
\special{pa 1508 1482}%
\special{pa 1526 1454}%
\special{pa 1522 1422}%
\special{pa 1498 1394}%
\special{pa 1482 1370}%
\special{pa 1490 1340}%
\special{pa 1522 1334}%
\special{pa 1550 1324}%
\special{pa 1568 1296}%
\special{pa 1582 1266}%
\special{pa 1598 1238}%
\special{pa 1610 1208}%
\special{pa 1600 1182}%
\special{pa 1576 1156}%
\special{pa 1552 1132}%
\special{pa 1534 1108}%
\special{pa 1520 1076}%
\special{pa 1510 1038}%
\special{pa 1514 1006}%
\special{pa 1538 1000}%
\special{pa 1572 1000}%
\special{pa 1592 974}%
\special{pa 1584 938}%
\special{pa 1562 914}%
\special{pa 1576 896}%
\special{pa 1606 874}%
\special{pa 1616 842}%
\special{pa 1606 808}%
\special{pa 1582 780}%
\special{pa 1570 756}%
\special{pa 1590 736}%
\special{pa 1620 712}%
\special{pa 1640 682}%
\special{pa 1632 656}%
\special{pa 1608 632}%
\special{pa 1600 604}%
\special{pa 1618 574}%
\special{pa 1644 554}%
\special{pa 1674 542}%
\special{pa 1704 532}%
\special{pa 1732 518}%
\special{pa 1762 502}%
\special{pa 1790 486}%
\special{pa 1800 480}%
\special{sp}%
% SPLINE 2 0 3 0
% 37 1000 480 1090 550 1170 550 1160 580 1260 640 1300 650 1240 670 1210 690 1280 710 1280 760 1240 770 1240 800 1310 840 1370 870 1400 890 1410 950 1460 1020 1460 1090 1450 1090 1400 1090 1380 1180 1390 1190 1440 1210 1410 1270 1390 1350 1300 1360 1250 1390 1290 1430 1250 1480 1310 1570 1280 1640 1290 1690 1340 1720 1340 1740 1340 1790 1370 1810 1400 1880
% 
\special{pn 8}%
\special{pa 1000 480}%
\special{pa 1018 510}%
\special{pa 1040 534}%
\special{pa 1070 548}%
\special{pa 1108 550}%
\special{pa 1150 542}%
\special{pa 1170 550}%
\special{pa 1160 578}%
\special{pa 1166 606}%
\special{pa 1190 624}%
\special{pa 1230 634}%
\special{pa 1274 642}%
\special{pa 1300 650}%
\special{pa 1282 658}%
\special{pa 1238 672}%
\special{pa 1210 690}%
\special{pa 1228 696}%
\special{pa 1268 702}%
\special{pa 1290 732}%
\special{pa 1278 762}%
\special{pa 1246 768}%
\special{pa 1238 796}%
\special{pa 1258 822}%
\special{pa 1288 834}%
\special{pa 1320 844}%
\special{pa 1348 858}%
\special{pa 1378 874}%
\special{pa 1402 894}%
\special{pa 1408 924}%
\special{pa 1412 956}%
\special{pa 1428 980}%
\special{pa 1450 1004}%
\special{pa 1468 1040}%
\special{pa 1472 1078}%
\special{pa 1450 1090}%
\special{pa 1416 1088}%
\special{pa 1390 1100}%
\special{pa 1374 1130}%
\special{pa 1374 1168}%
\special{pa 1394 1192}%
\special{pa 1428 1204}%
\special{pa 1442 1220}%
\special{pa 1424 1248}%
\special{pa 1406 1284}%
\special{pa 1400 1320}%
\special{pa 1394 1346}%
\special{pa 1374 1358}%
\special{pa 1342 1358}%
\special{pa 1302 1360}%
\special{pa 1264 1372}%
\special{pa 1252 1394}%
\special{pa 1278 1416}%
\special{pa 1288 1440}%
\special{pa 1264 1464}%
\special{pa 1248 1490}%
\special{pa 1264 1516}%
\special{pa 1292 1542}%
\special{pa 1310 1570}%
\special{pa 1304 1598}%
\special{pa 1286 1628}%
\special{pa 1278 1662}%
\special{pa 1290 1690}%
\special{pa 1322 1706}%
\special{pa 1342 1726}%
\special{pa 1338 1760}%
\special{pa 1340 1790}%
\special{pa 1366 1808}%
\special{pa 1388 1832}%
\special{pa 1398 1862}%
\special{pa 1400 1880}%
\special{sp}%
% SPLINE 2 0 3 0
% 42 600 480 620 510 640 520 670 530 660 580 620 580 660 630 730 660 760 680 720 730 770 780 800 810 740 860 710 900 710 950 780 950 830 980 850 1030 800 1090 770 1120 780 1150 790 1180 880 1180 840 1200 840 1250 900 1330 840 1360 870 1410 880 1450 880 1520 930 1540 1020 1600 1020 1640 960 1640 1010 1690 1090 1690 1160 1740 1200 1760 1250 1790 1260 1840 1300 1850 1400 1880
% 
\special{pn 8}%
\special{pa 600 480}%
\special{pa 618 508}%
\special{pa 646 522}%
\special{pa 674 536}%
\special{pa 670 570}%
\special{pa 640 582}%
\special{pa 618 584}%
\special{pa 638 612}%
\special{pa 672 638}%
\special{pa 700 650}%
\special{pa 730 660}%
\special{pa 758 676}%
\special{pa 748 700}%
\special{pa 722 726}%
\special{pa 726 750}%
\special{pa 758 772}%
\special{pa 790 794}%
\special{pa 800 816}%
\special{pa 778 838}%
\special{pa 744 858}%
\special{pa 720 880}%
\special{pa 708 910}%
\special{pa 708 944}%
\special{pa 728 956}%
\special{pa 764 952}%
\special{pa 798 954}%
\special{pa 824 974}%
\special{pa 844 1002}%
\special{pa 850 1032}%
\special{pa 838 1058}%
\special{pa 812 1082}%
\special{pa 782 1104}%
\special{pa 772 1128}%
\special{pa 782 1160}%
\special{pa 798 1186}%
\special{pa 840 1188}%
\special{pa 876 1182}%
\special{pa 870 1184}%
\special{pa 838 1202}%
\special{pa 834 1234}%
\special{pa 854 1268}%
\special{pa 884 1296}%
\special{pa 902 1320}%
\special{pa 892 1336}%
\special{pa 856 1348}%
\special{pa 840 1368}%
\special{pa 858 1394}%
\special{pa 878 1426}%
\special{pa 880 1458}%
\special{pa 876 1492}%
\special{pa 880 1520}%
\special{pa 906 1536}%
\special{pa 942 1544}%
\special{pa 974 1554}%
\special{pa 1000 1570}%
\special{pa 1020 1598}%
\special{pa 1024 1634}%
\special{pa 998 1644}%
\special{pa 964 1640}%
\special{pa 964 1656}%
\special{pa 996 1682}%
\special{pa 1032 1696}%
\special{pa 1064 1692}%
\special{pa 1094 1692}%
\special{pa 1120 1706}%
\special{pa 1146 1728}%
\special{pa 1172 1748}%
\special{pa 1202 1760}%
\special{pa 1234 1774}%
\special{pa 1252 1796}%
\special{pa 1256 1830}%
\special{pa 1276 1848}%
\special{pa 1312 1852}%
\special{pa 1344 1858}%
\special{pa 1374 1870}%
\special{pa 1400 1880}%
\special{sp}%
% SPLINE 2 0 3 0
% 27 1400 1880 1340 1820 1300 1790 1300 1740 1210 1680 1230 1630 1120 1590 1180 1510 1120 1460 1130 1430 1090 1370 1160 1330 1080 1230 1220 1170 1110 1090 1200 1040 1020 910 1130 850 1020 770 1060 760 1010 710 1080 670 840 560 890 520 860 480 860 480 860 480
% 
\special{pn 8}%
\special{pa 1400 1880}%
\special{pa 1382 1856}%
\special{pa 1358 1834}%
\special{pa 1330 1814}%
\special{pa 1304 1796}%
\special{pa 1300 1766}%
\special{pa 1298 1734}%
\special{pa 1268 1716}%
\special{pa 1232 1702}%
\special{pa 1210 1684}%
\special{pa 1222 1656}%
\special{pa 1230 1626}%
\special{pa 1204 1612}%
\special{pa 1162 1604}%
\special{pa 1126 1598}%
\special{pa 1124 1580}%
\special{pa 1150 1554}%
\special{pa 1178 1526}%
\special{pa 1174 1502}%
\special{pa 1142 1484}%
\special{pa 1120 1462}%
\special{pa 1130 1432}%
\special{pa 1116 1402}%
\special{pa 1092 1376}%
\special{pa 1102 1358}%
\special{pa 1138 1344}%
\special{pa 1162 1328}%
\special{pa 1152 1304}%
\special{pa 1122 1276}%
\special{pa 1090 1250}%
\special{pa 1082 1228}%
\special{pa 1102 1210}%
\special{pa 1142 1198}%
\special{pa 1184 1188}%
\special{pa 1216 1178}%
\special{pa 1218 1164}%
\special{pa 1190 1146}%
\special{pa 1152 1126}%
\special{pa 1120 1106}%
\special{pa 1112 1088}%
\special{pa 1136 1072}%
\special{pa 1172 1056}%
\special{pa 1200 1042}%
\special{pa 1200 1024}%
\special{pa 1174 1006}%
\special{pa 1134 986}%
\special{pa 1090 966}%
\special{pa 1050 946}%
\special{pa 1024 928}%
\special{pa 1020 910}%
\special{pa 1046 896}%
\special{pa 1086 882}%
\special{pa 1120 866}%
\special{pa 1130 848}%
\special{pa 1104 828}%
\special{pa 1062 806}%
\special{pa 1026 786}%
\special{pa 1018 772}%
\special{pa 1050 764}%
\special{pa 1056 752}%
\special{pa 1024 726}%
\special{pa 1012 704}%
\special{pa 1040 686}%
\special{pa 1076 672}%
\special{pa 1088 660}%
\special{pa 1070 648}%
\special{pa 1032 636}%
\special{pa 982 624}%
\special{pa 928 612}%
\special{pa 880 598}%
\special{pa 848 584}%
\special{pa 836 568}%
\special{pa 854 550}%
\special{pa 884 530}%
\special{pa 886 506}%
\special{pa 860 480}%
\special{sp}%
% SPLINE 2 0 3 0
% 26 3980 1880 4200 1780 4220 1780 4180 1710 4290 1680 4380 1610 4260 1540 4340 1490 4410 1440 4340 1370 4530 1320 4520 1250 4550 1190 4410 1070 4490 1000 4480 930 4550 860 4530 790 4650 740 4470 640 4680 560 4680 480 4620 410 4700 340 4700 300 4750 290
% 
\special{pn 8}%
\special{pa 3980 1880}%
\special{pa 4006 1860}%
\special{pa 4030 1840}%
\special{pa 4058 1822}%
\special{pa 4084 1806}%
\special{pa 4114 1792}%
\special{pa 4146 1784}%
\special{pa 4180 1780}%
\special{pa 4216 1782}%
\special{pa 4220 1766}%
\special{pa 4194 1732}%
\special{pa 4180 1704}%
\special{pa 4194 1690}%
\special{pa 4232 1686}%
\special{pa 4278 1682}%
\special{pa 4324 1670}%
\special{pa 4360 1652}%
\special{pa 4380 1628}%
\special{pa 4378 1606}%
\special{pa 4354 1588}%
\special{pa 4316 1572}%
\special{pa 4282 1558}%
\special{pa 4262 1544}%
\special{pa 4268 1526}%
\special{pa 4296 1510}%
\special{pa 4336 1492}%
\special{pa 4376 1476}%
\special{pa 4404 1458}%
\special{pa 4408 1436}%
\special{pa 4386 1412}%
\special{pa 4356 1388}%
\special{pa 4340 1368}%
\special{pa 4350 1358}%
\special{pa 4384 1354}%
\special{pa 4428 1352}%
\special{pa 4474 1348}%
\special{pa 4510 1340}%
\special{pa 4530 1322}%
\special{pa 4528 1292}%
\special{pa 4520 1258}%
\special{pa 4530 1228}%
\special{pa 4548 1200}%
\special{pa 4546 1174}%
\special{pa 4522 1150}%
\special{pa 4486 1130}%
\special{pa 4448 1110}%
\special{pa 4420 1090}%
\special{pa 4410 1072}%
\special{pa 4426 1054}%
\special{pa 4458 1036}%
\special{pa 4486 1012}%
\special{pa 4490 982}%
\special{pa 4480 950}%
\special{pa 4486 922}%
\special{pa 4512 900}%
\special{pa 4540 878}%
\special{pa 4550 852}%
\special{pa 4538 820}%
\special{pa 4530 790}%
\special{pa 4548 770}%
\special{pa 4586 758}%
\special{pa 4624 750}%
\special{pa 4650 742}%
\special{pa 4648 728}%
\special{pa 4626 714}%
\special{pa 4588 696}%
\special{pa 4546 678}%
\special{pa 4506 662}%
\special{pa 4478 646}%
\special{pa 4470 636}%
\special{pa 4482 628}%
\special{pa 4510 622}%
\special{pa 4550 614}%
\special{pa 4592 606}%
\special{pa 4634 592}%
\special{pa 4668 574}%
\special{pa 4688 546}%
\special{pa 4694 514}%
\special{pa 4680 482}%
\special{pa 4654 452}%
\special{pa 4628 428}%
\special{pa 4622 408}%
\special{pa 4642 392}%
\special{pa 4674 372}%
\special{pa 4698 348}%
\special{pa 4700 314}%
\special{pa 4710 292}%
\special{pa 4748 290}%
\special{pa 4750 290}%
\special{sp}%
% SPLINE 2 0 3 0
% 29 1400 1880 1660 1770 1720 1740 1620 1680 1700 1630 1820 1580 1710 1520 1840 1450 1820 1400 1730 1320 1850 1270 1970 1240 1960 1190 2000 1130 1850 1050 1780 1020 1870 990 1920 990 1960 960 1880 920 2000 870 2050 830 1910 740 1940 690 2020 660 2000 620 2050 560 2020 480 2020 480
% 
\special{pn 8}%
\special{pa 1400 1880}%
\special{pa 1424 1868}%
\special{pa 1450 1856}%
\special{pa 1474 1844}%
\special{pa 1500 1832}%
\special{pa 1530 1820}%
\special{pa 1560 1808}%
\special{pa 1592 1794}%
\special{pa 1628 1782}%
\special{pa 1668 1768}%
\special{pa 1704 1754}%
\special{pa 1720 1740}%
\special{pa 1702 1724}%
\special{pa 1664 1708}%
\special{pa 1630 1690}%
\special{pa 1620 1672}%
\special{pa 1638 1654}%
\special{pa 1676 1638}%
\special{pa 1724 1624}%
\special{pa 1768 1612}%
\special{pa 1804 1600}%
\special{pa 1822 1588}%
\special{pa 1814 1574}%
\special{pa 1782 1558}%
\special{pa 1744 1542}%
\special{pa 1716 1528}%
\special{pa 1714 1514}%
\special{pa 1738 1502}%
\special{pa 1778 1490}%
\special{pa 1816 1474}%
\special{pa 1840 1454}%
\special{pa 1836 1424}%
\special{pa 1814 1394}%
\special{pa 1784 1370}%
\special{pa 1754 1348}%
\special{pa 1734 1330}%
\special{pa 1732 1312}%
\special{pa 1750 1298}%
\special{pa 1782 1284}%
\special{pa 1822 1274}%
\special{pa 1868 1270}%
\special{pa 1912 1266}%
\special{pa 1948 1262}%
\special{pa 1968 1250}%
\special{pa 1968 1226}%
\special{pa 1960 1192}%
\special{pa 1974 1162}%
\special{pa 1996 1138}%
\special{pa 2004 1116}%
\special{pa 1990 1098}%
\special{pa 1962 1084}%
\special{pa 1924 1070}%
\special{pa 1880 1058}%
\special{pa 1834 1046}%
\special{pa 1796 1034}%
\special{pa 1780 1020}%
\special{pa 1796 1006}%
\special{pa 1832 996}%
\special{pa 1872 990}%
\special{pa 1904 992}%
\special{pa 1936 984}%
\special{pa 1960 962}%
\special{pa 1944 944}%
\special{pa 1908 930}%
\special{pa 1880 920}%
\special{pa 1884 910}%
\special{pa 1914 898}%
\special{pa 1958 886}%
\special{pa 2002 870}%
\special{pa 2036 850}%
\special{pa 2050 830}%
\special{pa 2038 814}%
\special{pa 2006 800}%
\special{pa 1968 786}%
\special{pa 1932 768}%
\special{pa 1912 746}%
\special{pa 1914 716}%
\special{pa 1940 690}%
\special{pa 1980 678}%
\special{pa 2014 670}%
\special{pa 2018 652}%
\special{pa 2000 622}%
\special{pa 2012 596}%
\special{pa 2040 572}%
\special{pa 2054 548}%
\special{pa 2046 518}%
\special{pa 2024 486}%
\special{pa 2020 480}%
\special{sp}%
% SPLINE 2 0 3 0
% 56 3980 1880 3950 1790 3930 1770 3950 1720 3920 1670 3980 1610 3880 1580 3880 1530 3930 1490 3840 1440 3910 1370 3940 1360 3880 1330 3840 1320 3840 1290 3920 1280 3990 1260 4030 1240 4000 1210 3950 1210 3950 1180 4020 1160 4040 1130 3990 1100 3910 1040 3860 1000 3890 960 3930 930 3990 900 4020 880 3960 850 3910 830 3870 810 3830 780 3890 750 3890 720 3850 700 3810 690 3780 670 3840 670 3900 650 3920 600 3860 600 3770 580 3750 540 3840 510 3870 480 3810 460 3680 430 3670 410 3740 380 3780 360 3680 330 3630 320 3630 290 3630 290
% 
\special{pn 8}%
\special{pa 3980 1880}%
\special{pa 3978 1846}%
\special{pa 3970 1816}%
\special{pa 3950 1790}%
\special{pa 3930 1766}%
\special{pa 3946 1736}%
\special{pa 3946 1706}%
\special{pa 3924 1680}%
\special{pa 3928 1654}%
\special{pa 3960 1628}%
\special{pa 3980 1610}%
\special{pa 3964 1604}%
\special{pa 3926 1598}%
\special{pa 3888 1586}%
\special{pa 3872 1558}%
\special{pa 3884 1526}%
\special{pa 3916 1504}%
\special{pa 3930 1488}%
\special{pa 3908 1474}%
\special{pa 3870 1460}%
\special{pa 3842 1442}%
\special{pa 3842 1418}%
\special{pa 3866 1392}%
\special{pa 3904 1372}%
\special{pa 3938 1362}%
\special{pa 3930 1346}%
\special{pa 3890 1332}%
\special{pa 3852 1328}%
\special{pa 3838 1306}%
\special{pa 3852 1280}%
\special{pa 3884 1278}%
\special{pa 3922 1280}%
\special{pa 3952 1274}%
\special{pa 3982 1264}%
\special{pa 4018 1252}%
\special{pa 4028 1232}%
\special{pa 4000 1210}%
\special{pa 3962 1214}%
\special{pa 3948 1196}%
\special{pa 3962 1170}%
\special{pa 3994 1166}%
\special{pa 4028 1156}%
\special{pa 4040 1128}%
\special{pa 4016 1110}%
\special{pa 3980 1096}%
\special{pa 3956 1078}%
\special{pa 3932 1058}%
\special{pa 3904 1038}%
\special{pa 3874 1018}%
\special{pa 3860 996}%
\special{pa 3880 970}%
\special{pa 3906 948}%
\special{pa 3932 930}%
\special{pa 3958 916}%
\special{pa 3990 902}%
\special{pa 4018 884}%
\special{pa 4010 868}%
\special{pa 3972 854}%
\special{pa 3936 842}%
\special{pa 3910 830}%
\special{pa 3882 816}%
\special{pa 3848 798}%
\special{pa 3830 780}%
\special{pa 3858 768}%
\special{pa 3890 750}%
\special{pa 3890 718}%
\special{pa 3862 702}%
\special{pa 3830 696}%
\special{pa 3798 684}%
\special{pa 3782 668}%
\special{pa 3812 668}%
\special{pa 3856 670}%
\special{pa 3886 664}%
\special{pa 3910 638}%
\special{pa 3922 606}%
\special{pa 3906 594}%
\special{pa 3868 600}%
\special{pa 3830 604}%
\special{pa 3798 598}%
\special{pa 3772 582}%
\special{pa 3752 552}%
\special{pa 3756 530}%
\special{pa 3786 524}%
\special{pa 3824 516}%
\special{pa 3860 498}%
\special{pa 3868 476}%
\special{pa 3840 464}%
\special{pa 3798 460}%
\special{pa 3760 460}%
\special{pa 3728 458}%
\special{pa 3700 448}%
\special{pa 3678 428}%
\special{pa 3672 402}%
\special{pa 3700 388}%
\special{pa 3742 380}%
\special{pa 3776 366}%
\special{pa 3778 352}%
\special{pa 3752 340}%
\special{pa 3708 332}%
\special{pa 3664 330}%
\special{pa 3634 324}%
\special{pa 3630 296}%
\special{pa 3630 290}%
\special{sp}%
% SPLINE 2 0 3 0
% 33 1400 1880 1570 1850 1660 1820 1780 1810 1900 1770 1840 1710 1900 1670 2040 1630 1980 1590 1980 1560 2100 1510 2120 1470 2080 1430 2170 1360 2090 1290 2070 1240 2120 1210 2140 1160 2070 1090 2100 1070 2150 1020 2080 960 2080 930 2170 880 2160 820 2110 780 2180 730 2220 690 2160 640 2150 610 2200 580 2210 520 2140 480
% 
\special{pn 8}%
\special{pa 1400 1880}%
\special{pa 1432 1878}%
\special{pa 1464 1874}%
\special{pa 1496 1868}%
\special{pa 1526 1862}%
\special{pa 1558 1854}%
\special{pa 1590 1844}%
\special{pa 1620 1832}%
\special{pa 1650 1824}%
\special{pa 1678 1818}%
\special{pa 1708 1814}%
\special{pa 1738 1812}%
\special{pa 1776 1810}%
\special{pa 1820 1806}%
\special{pa 1862 1800}%
\special{pa 1892 1788}%
\special{pa 1902 1772}%
\special{pa 1882 1750}%
\special{pa 1852 1728}%
\special{pa 1840 1704}%
\special{pa 1860 1686}%
\special{pa 1902 1670}%
\special{pa 1950 1660}%
\special{pa 1994 1650}%
\special{pa 2030 1644}%
\special{pa 2044 1636}%
\special{pa 2030 1626}%
\special{pa 2000 1610}%
\special{pa 1978 1584}%
\special{pa 1984 1554}%
\special{pa 2008 1538}%
\special{pa 2040 1532}%
\special{pa 2076 1524}%
\special{pa 2106 1506}%
\special{pa 2122 1476}%
\special{pa 2102 1450}%
\special{pa 2080 1430}%
\special{pa 2092 1410}%
\special{pa 2128 1390}%
\special{pa 2162 1370}%
\special{pa 2170 1352}%
\special{pa 2152 1334}%
\special{pa 2118 1312}%
\special{pa 2086 1286}%
\special{pa 2070 1254}%
\special{pa 2080 1232}%
\special{pa 2110 1216}%
\special{pa 2136 1192}%
\special{pa 2140 1162}%
\special{pa 2116 1134}%
\special{pa 2084 1112}%
\special{pa 2070 1092}%
\special{pa 2092 1076}%
\special{pa 2126 1054}%
\special{pa 2148 1030}%
\special{pa 2144 1010}%
\special{pa 2116 992}%
\special{pa 2086 970}%
\special{pa 2078 938}%
\special{pa 2096 914}%
\special{pa 2128 902}%
\special{pa 2160 890}%
\special{pa 2176 866}%
\special{pa 2168 832}%
\special{pa 2142 806}%
\special{pa 2116 786}%
\special{pa 2114 770}%
\special{pa 2140 752}%
\special{pa 2178 732}%
\special{pa 2210 710}%
\special{pa 2220 688}%
\special{pa 2200 670}%
\special{pa 2168 648}%
\special{pa 2150 618}%
\special{pa 2166 598}%
\special{pa 2198 582}%
\special{pa 2214 554}%
\special{pa 2210 520}%
\special{pa 2190 498}%
\special{pa 2158 486}%
\special{pa 2140 480}%
\special{sp}%
\end{picture}%
\end{center}
\caption{Illustrations of the
noncolliding Brownian motions $\X (t), t \in (0, T]$ 
in the left picture, 
and $\Y (t), t \in (0, \infty)$ in the right picture,
both start from ${\bf 0}$}
\end{figure}
%%%%%%%%%%%%%%%%%%%%%%%%%%%%%%%%%%%%%%%%%%
%%%%%%%%%%%%%%%%%%%%%%%%%%%%%%%%%%%%%%%%%%

When the time $T$ becomes infinity, the process $\X (t)$ converges to 
the temporally homogeneous process $\Y(t)$,
the {\bf noncolliding Brownian motion in an infinite
time-interval} $t \in (0, \infty)$,
whose transition probability density function $p_N(t,\y|\x)$
is given by follows;

%%%%%%%%%%%%%%%%%%%%%%%%%%%%%%%%%%%%%%%%%%%%%%%%
\begin{eqnarray}
\label{eqn:p(x,y)}
p_N(t,\y|\x)
&=& \frac{h_N(\y)}{h_N(\x)} f_N(t,\y|\x),
\quad t > 0, \quad \x, \y \in \WA, 
\\
%\label{eqn:p(0,y)}
p_N(t,\y|\0)
&=& \frac{t^{-N^2/2}}{C_1(N)} h_N(\y)^2 
\exp\left(-\frac{|\y|^2}{2t}\right),
\quad t > 0, \quad \y \in \WA.
\nonumber
\end{eqnarray}
%%%%%%%%%%%%%%%%%%%%%%%%%%%%%%%%%%%%
The above formulas are derived form (\ref{eqn:g(x,y)})
and (\ref{eqn:g(0,y)}) by taking the 
limit $T \to \infty$ using (\ref{eqn:asym}).
See the right picture of Fig. 1,
which illustrates $\Y(t)$ starting from ${\bf 0}$.

The first formula of (\ref{eqn:p(x,y)}) implies that
$\Y (t)$ is the Doob $h$-transformation
of the absorbing Brownian motion in $\WA$,
whose transition probability density is given 
by $f_N(t, \y|\x)$ \cite{Gra99}. 
One-parameter family of interacting Brownian motions
on $\R$ satisfying the following system of 
stochastic differential equations
\begin{equation}
Y_i(t)=B_i(t)+ \frac{\beta}{2} 
\sum_{\substack{1 \leq j \leq N\\ j \ne i}}
\int_{0}^{t}
\frac{1}{Y_i(s)-Y_j(s)} ds,
\quad 
1 \leq i \leq N
\label{eqn:Dyson0}
\end{equation}
is called {\bf Dyson's Brownian motion model}
with parameter $\beta >0$ \cite{Dys62b,Mehta}.
It is readily seen from (\ref{eqn:p(x,y)})
that $\Y(t), t \in (0, \infty)$ solves 
the system of equations (\ref{eqn:Dyson0})
with $\beta=2$.

By comparing
(\ref{eqn:GT2}) with (\ref{eqn:g(x,y)}), (\ref{eqn:g(0,y)}),
we find that $\X(t)$ can be
regarded as a multi-dimensional extension of
the Brownian meander $X(t)$.
Similarly, 
by comparing
(\ref{eqn:Bessel1}) with (\ref{eqn:p(x,y)}),
$\Y(t)$ can be considered to be a multi-dimensional version of
the three-dimensional Bessel process $Y(t)$.
Moreover, a multi-dimensional extension of 
Imhof's relation (\ref{Imhof:original}) 
is derived from (\ref{eqn:g(x,y)}), (\ref{eqn:g(0,y)}) 
and (\ref{eqn:p(x,y)}) as \cite{KT03a}
\begin{equation}
P(\X (\cdot) \in dw)=\frac{C_1(N)}{C_2(N)}
\frac{T^{N(N-1)/4}}{h(w(T))} P(\Y (\cdot) \in dw).
\label{Imhof:BM}
\end{equation}

%%%%%%%%%%%%%%%%%%%%%%%%%%%%%%%%%%%%%%%%%%%%%%%%%%%%%%%
\subsection{Noncolliding generalized meander
and noncolliding Bessel process}
%%%%%%%%%%%%%%%%%%%%%%%%%%%%%%%%%%%%%%%%%%%%%%%%%%%%%%%

We introduce the subsets $\WC$ and $\WD$ of $\R^N$ 
defined by
\begin{eqnarray}
&&\WC=\{{\x} \in \R^N: 0<x_1<x_2< \, \cdots \, <x_N\},
\nonumber\\
&&\WD=\{{\x} \in \R^N: 0\leq |x_1|<x_2< \, \cdots \, <x_N\},
\nonumber
\end{eqnarray}
which are called the Wyle chambers of 
type $C_{N}$ and of type $D_{N}$, respectively.
The Bessel process $Y^{(\nu)}(t), t \geq 0$
has the origin as a transient point when $\nu \geq 0$, 
and it has the origin as a recurrent point
when $-1 < \nu < 0$.
Then the state spaces $\W$ of the noncolliding generalized meander
and the noncolliding Bessel process,
which will be introduced in this subsection,
are $\WC$ when $\nu\ge 0$,
and $\WD$ when $-1 < \nu <0$.
If the Bessel process is defined by 
the square root of the squared Bessel process,
it can be a multi-valued stochastic process.
Consider the squared Bessel process 
starting from a positive initial point.
When $\nu \geq 0$, the process stays positive
with probability one,
and its square root is determined uniquely,
which coincides with the Bessel process introduced in Section 2.
While, when $-1 < \nu < 0$, it hits the origin with probability one
and then the square root process becomes a bi-valued process 
after hitting.
The generalized meander and the leftmost particles
in the $N$ particle systems of the noncolliding generalized meander
and of the noncolliding Bessel process
are in the same situation.
The absolute value $|x_1|$ appearing in the definition of $\WD$ 
implies that bi-valued processes are allowed,
when $-1 < \nu < 0$.
See Fig 2.
However, we usually consider only nonnegative parts
of such bi-valued processes
just for simplicity of explanation.

The density function of an 
$N$-component generalized meander
at time $t$, which 
starts from ${\x}$ in $\W$ at time $s$ 
and stays in $\W$ up to time $t$, is given by
\begin{equation*}
f_N^{(\nu,\kappa)}(s,{\x};t,{\y})
=\det_{1 \leq i,j \leq N}
\Big(G^{(\nu,\kappa)}(s,x_i; t,y_j) \Big)
\end{equation*}
from the Karlin-McGregor formula.
The probability that the process stays in  $\W$ 
during the time-interval $(0,t]$ is given by
\begin{equation*}
\cNnk(t,{\x})
=\int_{\W}d{\y}f_N^{(\nu,\kappa)}(0, {\x}; t, {\y}).
\end{equation*}
Then the transition probability density function 
of the {\bf noncolliding generalized meander} 
${\bf X}^{(\nu,\kappa)}(t)
=(X^{(\nu,\kappa)}_1(t),X^{(\nu,\kappa)}_2(t),
 \dots ,X^{(\nu,\kappa)}_N(t))$
is given by
$$
g^{(\nu,\kappa)}_{N,T} (s,\x,t,\y) 
= \frac{\cNnk(T-t, {\y})}{\cNnk(T-s, {\x})}
f^{(\nu,\kappa)}_N(t-s,\x,\y), \quad
0 \leq s < t \leq T, \, \x, \y \in \W.
$$

%%%%%%%%%%%%%%%%%%%%%%%%%%%%%%%%%%%%%%%%%%
%%%%%%%%%%%%%%%%%%%%%%%%%%%%%%%%%%%%%%%%%%
\vskip 3mm
\begin{figure}[ch]
\begin{center}
%\input fig_ronsetu_2.1.tex
%WinTpicVersion3.08
\unitlength 0.1in
\begin{picture}( 20.6000, 16.7500)(  0.0000,-18.0500)
% VECTOR 2 0 3 0
% 2 400 1800 2060 1800
% 
\special{pn 8}%
\special{pa 400 1800}%
\special{pa 2060 1800}%
\special{fp}%
\special{sh 1}%
\special{pa 2060 1800}%
\special{pa 1994 1780}%
\special{pa 2008 1800}%
\special{pa 1994 1820}%
\special{pa 2060 1800}%
\special{fp}%
% SPLINE 2 0 3 0
% 33 400 1800 650 1740 710 1720 660 1620 720 1590 830 1530 780 1480 810 1440 890 1390 800 1340 810 1260 890 1240 970 1190 930 1150 1000 1080 970 1010 910 980 870 890 910 830 1020 810 1100 790 1090 750 1050 750 1040 690 1090 680 1110 660 1150 620 1090 570 1120 530 1110 480 1030 470 1090 420 1130 400
% 
\special{pn 8}%
\special{pa 400 1800}%
\special{pa 428 1792}%
\special{pa 454 1784}%
\special{pa 482 1776}%
\special{pa 512 1768}%
\special{pa 542 1760}%
\special{pa 576 1754}%
\special{pa 610 1746}%
\special{pa 650 1740}%
\special{pa 688 1734}%
\special{pa 710 1720}%
\special{pa 704 1696}%
\special{pa 682 1664}%
\special{pa 662 1634}%
\special{pa 666 1610}%
\special{pa 696 1598}%
\special{pa 738 1586}%
\special{pa 782 1572}%
\special{pa 816 1558}%
\special{pa 832 1540}%
\special{pa 820 1520}%
\special{pa 792 1498}%
\special{pa 780 1474}%
\special{pa 800 1448}%
\special{pa 836 1426}%
\special{pa 872 1408}%
\special{pa 890 1394}%
\special{pa 880 1382}%
\special{pa 846 1370}%
\special{pa 812 1352}%
\special{pa 792 1324}%
\special{pa 794 1288}%
\special{pa 812 1260}%
\special{pa 836 1248}%
\special{pa 870 1244}%
\special{pa 910 1236}%
\special{pa 948 1222}%
\special{pa 972 1202}%
\special{pa 960 1180}%
\special{pa 934 1156}%
\special{pa 936 1136}%
\special{pa 964 1114}%
\special{pa 994 1090}%
\special{pa 1004 1060}%
\special{pa 990 1030}%
\special{pa 964 1006}%
\special{pa 934 992}%
\special{pa 906 978}%
\special{pa 884 952}%
\special{pa 872 920}%
\special{pa 872 888}%
\special{pa 884 858}%
\special{pa 906 834}%
\special{pa 932 818}%
\special{pa 962 812}%
\special{pa 996 810}%
\special{pa 1034 812}%
\special{pa 1072 812}%
\special{pa 1098 798}%
\special{pa 1100 764}%
\special{pa 1076 750}%
\special{pa 1044 746}%
\special{pa 1032 712}%
\special{pa 1048 686}%
\special{pa 1082 682}%
\special{pa 1106 664}%
\special{pa 1136 640}%
\special{pa 1150 618}%
\special{pa 1126 598}%
\special{pa 1094 578}%
\special{pa 1098 556}%
\special{pa 1122 528}%
\special{pa 1122 494}%
\special{pa 1094 474}%
\special{pa 1056 472}%
\special{pa 1030 470}%
\special{pa 1042 454}%
\special{pa 1074 430}%
\special{pa 1106 412}%
\special{pa 1130 400}%
\special{sp}%
% SPLINE 2 0 3 0
% 33 400 1800 740 1770 990 1760 1110 1730 1280 1670 1230 1610 1160 1520 1170 1500 1250 1500 1310 1420 1230 1380 1330 1330 1390 1320 1470 1290 1480 1240 1420 1200 1500 1180 1500 1130 1330 1070 1360 1040 1420 980 1580 920 1660 900 1660 880 1520 830 1560 780 1670 780 1680 700 1600 640 1620 610 1720 550 1680 450 1710 400
% 
\special{pn 8}%
\special{pa 400 1800}%
\special{pa 432 1796}%
\special{pa 464 1792}%
\special{pa 494 1788}%
\special{pa 526 1784}%
\special{pa 558 1782}%
\special{pa 590 1778}%
\special{pa 622 1776}%
\special{pa 654 1774}%
\special{pa 686 1772}%
\special{pa 718 1770}%
\special{pa 752 1770}%
\special{pa 784 1770}%
\special{pa 818 1770}%
\special{pa 850 1772}%
\special{pa 884 1770}%
\special{pa 916 1770}%
\special{pa 946 1768}%
\special{pa 976 1764}%
\special{pa 1006 1758}%
\special{pa 1034 1750}%
\special{pa 1064 1742}%
\special{pa 1098 1734}%
\special{pa 1136 1726}%
\special{pa 1178 1720}%
\special{pa 1218 1712}%
\special{pa 1252 1702}%
\special{pa 1274 1688}%
\special{pa 1280 1670}%
\special{pa 1268 1646}%
\special{pa 1244 1622}%
\special{pa 1212 1596}%
\special{pa 1184 1574}%
\special{pa 1164 1548}%
\special{pa 1160 1520}%
\special{pa 1178 1496}%
\special{pa 1208 1496}%
\special{pa 1244 1502}%
\special{pa 1280 1488}%
\special{pa 1306 1460}%
\special{pa 1314 1428}%
\special{pa 1294 1406}%
\special{pa 1258 1392}%
\special{pa 1232 1382}%
\special{pa 1236 1368}%
\special{pa 1266 1352}%
\special{pa 1306 1338}%
\special{pa 1344 1328}%
\special{pa 1376 1322}%
\special{pa 1406 1318}%
\special{pa 1438 1312}%
\special{pa 1466 1294}%
\special{pa 1484 1266}%
\special{pa 1476 1234}%
\special{pa 1442 1212}%
\special{pa 1420 1200}%
\special{pa 1440 1198}%
\special{pa 1482 1190}%
\special{pa 1508 1170}%
\special{pa 1504 1136}%
\special{pa 1476 1110}%
\special{pa 1436 1098}%
\special{pa 1394 1092}%
\special{pa 1356 1088}%
\special{pa 1332 1082}%
\special{pa 1334 1066}%
\special{pa 1358 1044}%
\special{pa 1382 1020}%
\special{pa 1404 996}%
\special{pa 1426 976}%
\special{pa 1450 958}%
\special{pa 1476 942}%
\special{pa 1506 932}%
\special{pa 1540 924}%
\special{pa 1576 920}%
\special{pa 1616 920}%
\special{pa 1648 914}%
\special{pa 1662 890}%
\special{pa 1648 864}%
\special{pa 1614 854}%
\special{pa 1572 850}%
\special{pa 1536 844}%
\special{pa 1520 830}%
\special{pa 1532 802}%
\special{pa 1562 780}%
\special{pa 1596 778}%
\special{pa 1630 784}%
\special{pa 1660 784}%
\special{pa 1682 766}%
\special{pa 1690 732}%
\special{pa 1678 698}%
\special{pa 1648 676}%
\special{pa 1616 660}%
\special{pa 1600 642}%
\special{pa 1616 614}%
\special{pa 1648 594}%
\special{pa 1682 578}%
\special{pa 1710 566}%
\special{pa 1720 546}%
\special{pa 1710 520}%
\special{pa 1692 488}%
\special{pa 1680 456}%
\special{pa 1688 426}%
\special{pa 1710 400}%
\special{sp}%
% SPLINE 2 0 3 0
% 32 400 1800 810 1750 860 1730 840 1690 860 1670 1080 1620 1070 1550 1010 1530 1020 1450 980 1390 1070 1330 1180 1320 1180 1280 1160 1180 1130 1140 1200 1130 1230 1090 1110 1020 1110 940 1240 910 1240 870 1200 820 1260 780 1310 700 1260 670 1410 630 1550 570 1440 530 1470 490 1530 460 1370 430 1310 400
% 
\special{pn 8}%
\special{pa 400 1800}%
\special{pa 426 1798}%
\special{pa 452 1796}%
\special{pa 478 1792}%
\special{pa 504 1790}%
\special{pa 532 1786}%
\special{pa 562 1784}%
\special{pa 592 1780}%
\special{pa 624 1776}%
\special{pa 656 1772}%
\special{pa 692 1766}%
\special{pa 730 1762}%
\special{pa 772 1756}%
\special{pa 814 1750}%
\special{pa 852 1740}%
\special{pa 856 1720}%
\special{pa 840 1690}%
\special{pa 864 1668}%
\special{pa 900 1658}%
\special{pa 938 1654}%
\special{pa 976 1656}%
\special{pa 1010 1656}%
\special{pa 1040 1652}%
\special{pa 1066 1642}%
\special{pa 1082 1618}%
\special{pa 1086 1584}%
\special{pa 1074 1554}%
\special{pa 1042 1540}%
\special{pa 1012 1532}%
\special{pa 1010 1508}%
\special{pa 1020 1472}%
\special{pa 1016 1438}%
\special{pa 996 1412}%
\special{pa 980 1388}%
\special{pa 990 1364}%
\special{pa 1018 1344}%
\special{pa 1060 1332}%
\special{pa 1102 1330}%
\special{pa 1142 1334}%
\special{pa 1170 1332}%
\special{pa 1182 1314}%
\special{pa 1180 1280}%
\special{pa 1180 1246}%
\special{pa 1176 1216}%
\special{pa 1164 1186}%
\special{pa 1140 1156}%
\special{pa 1132 1138}%
\special{pa 1164 1136}%
\special{pa 1204 1130}%
\special{pa 1228 1104}%
\special{pa 1226 1080}%
\special{pa 1204 1064}%
\special{pa 1170 1054}%
\special{pa 1134 1040}%
\special{pa 1106 1016}%
\special{pa 1096 980}%
\special{pa 1104 948}%
\special{pa 1128 930}%
\special{pa 1164 928}%
\special{pa 1202 928}%
\special{pa 1232 920}%
\special{pa 1244 896}%
\special{pa 1236 862}%
\special{pa 1210 836}%
\special{pa 1202 816}%
\special{pa 1224 800}%
\special{pa 1260 780}%
\special{pa 1294 754}%
\special{pa 1314 726}%
\special{pa 1310 700}%
\special{pa 1282 680}%
\special{pa 1258 668}%
\special{pa 1260 658}%
\special{pa 1284 652}%
\special{pa 1326 644}%
\special{pa 1376 636}%
\special{pa 1430 626}%
\special{pa 1482 614}%
\special{pa 1524 600}%
\special{pa 1550 586}%
\special{pa 1556 574}%
\special{pa 1534 566}%
\special{pa 1498 558}%
\special{pa 1460 548}%
\special{pa 1440 532}%
\special{pa 1452 506}%
\special{pa 1484 482}%
\special{pa 1516 468}%
\special{pa 1532 458}%
\special{pa 1520 452}%
\special{pa 1490 448}%
\special{pa 1448 444}%
\special{pa 1402 438}%
\special{pa 1362 428}%
\special{pa 1332 414}%
\special{pa 1310 400}%
\special{sp}%
% SPLINE 2 0 3 0
% 26 400 1800 550 1710 450 1660 610 1590 610 1500 540 1470 670 1380 560 1310 520 1260 560 1230 640 1210 640 1170 790 1090 750 1020 650 990 610 920 680 880 720 860 700 810 700 740 830 690 760 620 820 580 760 520 930 460 910 390
% 
\special{pn 8}%
\special{pa 400 1800}%
\special{pa 450 1784}%
\special{pa 494 1768}%
\special{pa 530 1752}%
\special{pa 552 1734}%
\special{pa 556 1718}%
\special{pa 538 1702}%
\special{pa 504 1686}%
\special{pa 470 1672}%
\special{pa 450 1660}%
\special{pa 456 1652}%
\special{pa 482 1644}%
\special{pa 520 1636}%
\special{pa 562 1622}%
\special{pa 598 1602}%
\special{pa 622 1572}%
\special{pa 628 1538}%
\special{pa 616 1506}%
\special{pa 586 1488}%
\special{pa 552 1476}%
\special{pa 540 1466}%
\special{pa 556 1450}%
\special{pa 592 1432}%
\special{pa 632 1410}%
\special{pa 662 1390}%
\special{pa 672 1372}%
\special{pa 656 1356}%
\special{pa 624 1342}%
\special{pa 584 1324}%
\special{pa 548 1302}%
\special{pa 524 1276}%
\special{pa 524 1252}%
\special{pa 554 1232}%
\special{pa 592 1226}%
\special{pa 626 1222}%
\special{pa 642 1204}%
\special{pa 640 1172}%
\special{pa 656 1146}%
\special{pa 688 1132}%
\special{pa 726 1124}%
\special{pa 764 1114}%
\special{pa 788 1098}%
\special{pa 790 1072}%
\special{pa 774 1042}%
\special{pa 746 1018}%
\special{pa 716 1008}%
\special{pa 684 1002}%
\special{pa 654 992}%
\special{pa 626 968}%
\special{pa 612 938}%
\special{pa 614 912}%
\special{pa 638 894}%
\special{pa 674 882}%
\special{pa 710 870}%
\special{pa 720 850}%
\special{pa 704 818}%
\special{pa 690 784}%
\special{pa 694 752}%
\special{pa 718 728}%
\special{pa 756 716}%
\special{pa 798 708}%
\special{pa 826 700}%
\special{pa 828 686}%
\special{pa 800 664}%
\special{pa 768 640}%
\special{pa 762 618}%
\special{pa 792 600}%
\special{pa 820 582}%
\special{pa 810 560}%
\special{pa 778 536}%
\special{pa 758 516}%
\special{pa 770 504}%
\special{pa 802 498}%
\special{pa 846 492}%
\special{pa 890 482}%
\special{pa 924 468}%
\special{pa 934 446}%
\special{pa 922 414}%
\special{pa 910 390}%
\special{sp}%
% STR 2 0 3 0
% 3 420 430 420 530 2 0
% $T$
\put(4.2000,-5.3000){\makebox(0,0)[lb]{$T$}}%
% STR 2 0 3 0
% 3 290 1850 290 1950 2 0
% $0$
\put(2.9000,-19.5000){\makebox(0,0)[lb]{$0$}}%
% STR 2 0 3 0
% 3 1850 1850 1850 1950 2 0
% $x$
\put(18.5000,-19.5000){\makebox(0,0)[lb]{$x$}}%
% STR 2 0 3 0
% 3 470 200 470 300 2 0
% $t$
\put(4.7000,-3.0000){\makebox(0,0)[lb]{$t$}}%
% LINE 3 0 3 0
% 48 390 1110 270 1230 390 1170 270 1290 390 1230 270 1350 390 1290 270 1410 390 1350 270 1470 390 1410 270 1530 390 1470 270 1590 390 1530 270 1650 390 1590 270 1710 390 1650 270 1770 390 1710 310 1790 390 1050 270 1170 390 990 270 1110 390 930 270 1050 390 870 270 990 390 810 270 930 390 750 270 870 390 690 270 810 390 630 270 750 390 570 270 690 390 510 270 630 390 450 270 570 380 400 270 510 320 400 270 450
% 
\special{pn 4}%
\special{pa 390 1110}%
\special{pa 270 1230}%
\special{fp}%
\special{pa 390 1170}%
\special{pa 270 1290}%
\special{fp}%
\special{pa 390 1230}%
\special{pa 270 1350}%
\special{fp}%
\special{pa 390 1290}%
\special{pa 270 1410}%
\special{fp}%
\special{pa 390 1350}%
\special{pa 270 1470}%
\special{fp}%
\special{pa 390 1410}%
\special{pa 270 1530}%
\special{fp}%
\special{pa 390 1470}%
\special{pa 270 1590}%
\special{fp}%
\special{pa 390 1530}%
\special{pa 270 1650}%
\special{fp}%
\special{pa 390 1590}%
\special{pa 270 1710}%
\special{fp}%
\special{pa 390 1650}%
\special{pa 270 1770}%
\special{fp}%
\special{pa 390 1710}%
\special{pa 310 1790}%
\special{fp}%
\special{pa 390 1050}%
\special{pa 270 1170}%
\special{fp}%
\special{pa 390 990}%
\special{pa 270 1110}%
\special{fp}%
\special{pa 390 930}%
\special{pa 270 1050}%
\special{fp}%
\special{pa 390 870}%
\special{pa 270 990}%
\special{fp}%
\special{pa 390 810}%
\special{pa 270 930}%
\special{fp}%
\special{pa 390 750}%
\special{pa 270 870}%
\special{fp}%
\special{pa 390 690}%
\special{pa 270 810}%
\special{fp}%
\special{pa 390 630}%
\special{pa 270 750}%
\special{fp}%
\special{pa 390 570}%
\special{pa 270 690}%
\special{fp}%
\special{pa 390 510}%
\special{pa 270 630}%
\special{fp}%
\special{pa 390 450}%
\special{pa 270 570}%
\special{fp}%
\special{pa 380 400}%
\special{pa 270 510}%
\special{fp}%
\special{pa 320 400}%
\special{pa 270 450}%
\special{fp}%
% BOX 2 5 3 0
% 2 270 400 390 1790
% 
\special{pn 8}%
\special{pa 270 400}%
\special{pa 390 400}%
\special{pa 390 1790}%
\special{pa 270 1790}%
\special{pa 270 400}%
\special{ip}%
% LINE 1 2 3 0
% 2 200 400 1850 400
% 
\special{pn 13}%
\special{pa 200 400}%
\special{pa 1850 400}%
\special{dt 0.045}%
% LINE 2 0 3 0
% 2 0 1800 400 1800
% 
\special{pn 8}%
\special{pa 0 1800}%
\special{pa 400 1800}%
\special{fp}%
% VECTOR 2 0 3 0
% 2 400 1800 400 200
% 
\special{pn 8}%
\special{pa 400 1800}%
\special{pa 400 200}%
\special{fp}%
\special{sh 1}%
\special{pa 400 200}%
\special{pa 380 268}%
\special{pa 400 254}%
\special{pa 420 268}%
\special{pa 400 200}%
\special{fp}%
\end{picture}%
\hspace*{20mm}
%\input fig_ronsetu_2.2.tex
%WinTpicVersion3.08
\unitlength 0.1in
\begin{picture}( 20.9000, 16.8500)(  1.8000,-18.0500)
% VECTOR 2 0 3 0
% 2 610 1800 2270 1800
% 
\special{pn 8}%
\special{pa 610 1800}%
\special{pa 2270 1800}%
\special{fp}%
\special{sh 1}%
\special{pa 2270 1800}%
\special{pa 2204 1780}%
\special{pa 2218 1800}%
\special{pa 2204 1820}%
\special{pa 2270 1800}%
\special{fp}%
% STR 2 0 3 0
% 3 640 430 640 530 2 0
% $T$
\put(6.4000,-5.3000){\makebox(0,0)[lb]{$T$}}%
% STR 2 0 3 0
% 3 500 1850 500 1950 2 0
% $0$
\put(5.0000,-19.5000){\makebox(0,0)[lb]{$0$}}%
% STR 2 0 3 0
% 3 2060 1860 2060 1960 2 0
% $x$
\put(20.6000,-19.6000){\makebox(0,0)[lb]{$x$}}%
% STR 2 0 3 0
% 3 680 190 680 290 2 0
% $t$
\put(6.8000,-2.9000){\makebox(0,0)[lb]{$t$}}%
% LINE 1 2 3 0
% 6 190 390 2060 390 2060 390 2060 390 2060 390 2060 390
% 
\special{pn 13}%
\special{pa 190 390}%
\special{pa 2060 390}%
\special{dt 0.045}%
\special{pa 2060 390}%
\special{pa 2060 390}%
\special{dt 0.045}%
\special{pa 2060 390}%
\special{pa 2060 390}%
\special{dt 0.045}%
% LINE 2 0 3 0
% 2 210 1800 610 1800
% 
\special{pn 8}%
\special{pa 210 1800}%
\special{pa 610 1800}%
\special{fp}%
% VECTOR 2 0 3 0
% 2 610 1800 610 200
% 
\special{pn 8}%
\special{pa 610 1800}%
\special{pa 610 200}%
\special{fp}%
\special{sh 1}%
\special{pa 610 200}%
\special{pa 590 268}%
\special{pa 610 254}%
\special{pa 630 268}%
\special{pa 610 200}%
\special{fp}%
% SPLINE 2 0 3 0
% 43 610 1800 690 1770 700 1740 680 1710 740 1670 740 1640 710 1610 650 1590 610 1550 670 1540 690 1520 690 1490 720 1480 660 1420 610 1400 630 1360 670 1360 760 1340 800 1290 760 1270 790 1220 770 1180 700 1150 660 1150 610 1120 660 1100 720 1100 820 1070 820 1020 780 1000 880 920 920 860 920 800 840 750 830 710 880 700 900 640 810 600 830 550 890 530 980 490 980 450 1040 400
% 
\special{pn 8}%
\special{pa 610 1800}%
\special{pa 644 1796}%
\special{pa 674 1784}%
\special{pa 696 1762}%
\special{pa 694 1730}%
\special{pa 680 1706}%
\special{pa 708 1692}%
\special{pa 740 1672}%
\special{pa 740 1640}%
\special{pa 718 1616}%
\special{pa 690 1604}%
\special{pa 660 1594}%
\special{pa 626 1574}%
\special{pa 610 1552}%
\special{pa 634 1546}%
\special{pa 672 1540}%
\special{pa 690 1514}%
\special{pa 696 1486}%
\special{pa 724 1478}%
\special{pa 716 1454}%
\special{pa 678 1428}%
\special{pa 636 1414}%
\special{pa 612 1402}%
\special{pa 618 1372}%
\special{pa 646 1356}%
\special{pa 678 1362}%
\special{pa 710 1362}%
\special{pa 740 1352}%
\special{pa 774 1332}%
\special{pa 802 1304}%
\special{pa 790 1284}%
\special{pa 760 1270}%
\special{pa 772 1246}%
\special{pa 792 1216}%
\special{pa 776 1186}%
\special{pa 750 1164}%
\special{pa 720 1152}%
\special{pa 690 1150}%
\special{pa 656 1150}%
\special{pa 620 1134}%
\special{pa 614 1116}%
\special{pa 646 1102}%
\special{pa 682 1098}%
\special{pa 712 1100}%
\special{pa 746 1102}%
\special{pa 780 1098}%
\special{pa 810 1084}%
\special{pa 826 1054}%
\special{pa 822 1022}%
\special{pa 792 1006}%
\special{pa 776 994}%
\special{pa 792 976}%
\special{pa 826 956}%
\special{pa 864 932}%
\special{pa 894 908}%
\special{pa 912 882}%
\special{pa 922 852}%
\special{pa 926 818}%
\special{pa 914 792}%
\special{pa 886 778}%
\special{pa 856 764}%
\special{pa 832 736}%
\special{pa 834 708}%
\special{pa 866 704}%
\special{pa 896 688}%
\special{pa 904 652}%
\special{pa 886 628}%
\special{pa 850 618}%
\special{pa 818 610}%
\special{pa 810 586}%
\special{pa 826 554}%
\special{pa 852 538}%
\special{pa 886 532}%
\special{pa 922 526}%
\special{pa 958 516}%
\special{pa 978 496}%
\special{pa 980 466}%
\special{pa 986 436}%
\special{pa 1010 416}%
\special{pa 1040 400}%
\special{sp}%
% SPLINE 2 0 3 0
% 38 610 1800 770 1770 830 1700 860 1610 780 1580 840 1540 940 1510 920 1480 900 1410 940 1370 850 1320 860 1260 860 1240 1020 1200 1050 1110 1050 1090 1010 1080 970 1080 1010 1030 1050 1020 1110 1000 1070 980 1020 970 1010 930 1120 910 1160 890 1180 830 1060 770 1040 730 1110 730 1180 680 1180 640 1140 610 1200 570 1200 530 1130 480 1170 440 1240 400
% 
\special{pn 8}%
\special{pa 610 1800}%
\special{pa 646 1800}%
\special{pa 680 1798}%
\special{pa 712 1792}%
\special{pa 742 1784}%
\special{pa 768 1772}%
\special{pa 790 1754}%
\special{pa 810 1730}%
\special{pa 832 1700}%
\special{pa 854 1666}%
\special{pa 866 1634}%
\special{pa 860 1610}%
\special{pa 826 1596}%
\special{pa 792 1588}%
\special{pa 782 1574}%
\special{pa 804 1554}%
\special{pa 848 1538}%
\special{pa 892 1530}%
\special{pa 928 1526}%
\special{pa 940 1512}%
\special{pa 926 1488}%
\special{pa 904 1456}%
\special{pa 896 1426}%
\special{pa 912 1398}%
\special{pa 938 1376}%
\special{pa 932 1360}%
\special{pa 898 1348}%
\special{pa 862 1332}%
\special{pa 848 1308}%
\special{pa 858 1274}%
\special{pa 860 1242}%
\special{pa 874 1222}%
\special{pa 904 1214}%
\special{pa 942 1214}%
\special{pa 982 1212}%
\special{pa 1016 1204}%
\special{pa 1036 1182}%
\special{pa 1046 1152}%
\special{pa 1050 1116}%
\special{pa 1048 1086}%
\special{pa 1014 1080}%
\special{pa 976 1082}%
\special{pa 974 1066}%
\special{pa 1000 1038}%
\special{pa 1030 1024}%
\special{pa 1066 1018}%
\special{pa 1104 1008}%
\special{pa 1104 994}%
\special{pa 1070 980}%
\special{pa 1034 978}%
\special{pa 1014 958}%
\special{pa 1012 926}%
\special{pa 1032 912}%
\special{pa 1068 910}%
\special{pa 1108 912}%
\special{pa 1142 904}%
\special{pa 1166 884}%
\special{pa 1180 852}%
\special{pa 1178 824}%
\special{pa 1158 810}%
\special{pa 1126 800}%
\special{pa 1090 790}%
\special{pa 1058 768}%
\special{pa 1040 738}%
\special{pa 1054 724}%
\special{pa 1090 730}%
\special{pa 1128 728}%
\special{pa 1158 714}%
\special{pa 1178 686}%
\special{pa 1186 654}%
\special{pa 1164 628}%
\special{pa 1140 610}%
\special{pa 1158 596}%
\special{pa 1194 578}%
\special{pa 1206 546}%
\special{pa 1188 518}%
\special{pa 1154 500}%
\special{pa 1132 484}%
\special{pa 1142 462}%
\special{pa 1172 440}%
\special{pa 1202 422}%
\special{pa 1228 406}%
\special{pa 1240 400}%
\special{sp}%
% SPLINE 2 0 3 0
% 31 610 1800 870 1780 990 1750 950 1720 990 1690 1050 1660 1090 1560 1190 1510 1230 1470 1110 1420 1070 1360 1170 1330 1260 1270 1240 1240 1170 1130 1170 1070 1220 980 1390 950 1340 900 1370 860 1420 810 1400 740 1320 710 1360 670 1430 640 1450 600 1400 560 1360 520 1450 500 1480 490 1380 400
% 
\special{pn 8}%
\special{pa 610 1800}%
\special{pa 634 1798}%
\special{pa 658 1796}%
\special{pa 682 1794}%
\special{pa 710 1792}%
\special{pa 742 1790}%
\special{pa 778 1788}%
\special{pa 818 1784}%
\special{pa 864 1780}%
\special{pa 914 1778}%
\special{pa 960 1772}%
\special{pa 990 1764}%
\special{pa 990 1750}%
\special{pa 962 1732}%
\special{pa 954 1712}%
\special{pa 986 1692}%
\special{pa 1020 1680}%
\special{pa 1046 1666}%
\special{pa 1062 1640}%
\special{pa 1072 1608}%
\special{pa 1082 1576}%
\special{pa 1098 1552}%
\special{pa 1124 1538}%
\special{pa 1154 1526}%
\special{pa 1188 1512}%
\special{pa 1218 1490}%
\special{pa 1230 1468}%
\special{pa 1218 1454}%
\special{pa 1186 1446}%
\special{pa 1146 1436}%
\special{pa 1108 1418}%
\special{pa 1080 1392}%
\special{pa 1070 1364}%
\special{pa 1082 1348}%
\special{pa 1112 1340}%
\special{pa 1152 1334}%
\special{pa 1194 1324}%
\special{pa 1232 1310}%
\special{pa 1256 1290}%
\special{pa 1258 1266}%
\special{pa 1240 1240}%
\special{pa 1216 1212}%
\special{pa 1196 1186}%
\special{pa 1180 1160}%
\special{pa 1170 1132}%
\special{pa 1168 1100}%
\special{pa 1170 1068}%
\special{pa 1178 1034}%
\special{pa 1192 1004}%
\special{pa 1216 982}%
\special{pa 1256 970}%
\special{pa 1302 966}%
\special{pa 1346 964}%
\special{pa 1380 962}%
\special{pa 1394 956}%
\special{pa 1380 942}%
\special{pa 1352 920}%
\special{pa 1342 892}%
\special{pa 1362 868}%
\special{pa 1390 846}%
\special{pa 1414 824}%
\special{pa 1424 794}%
\special{pa 1418 762}%
\special{pa 1392 736}%
\special{pa 1354 726}%
\special{pa 1324 716}%
\special{pa 1326 698}%
\special{pa 1354 674}%
\special{pa 1386 660}%
\special{pa 1416 650}%
\special{pa 1440 630}%
\special{pa 1450 600}%
\special{pa 1434 578}%
\special{pa 1400 560}%
\special{pa 1370 536}%
\special{pa 1362 516}%
\special{pa 1384 508}%
\special{pa 1424 504}%
\special{pa 1464 498}%
\special{pa 1486 486}%
\special{pa 1480 468}%
\special{pa 1456 446}%
\special{pa 1416 422}%
\special{pa 1380 400}%
\special{sp}%
% SPLINE 2 0 3 0
% 35 610 1800 920 1790 1050 1770 1120 1750 1080 1710 1170 1640 1310 1570 1400 1540 1310 1470 1350 1440 1450 1430 1480 1380 1430 1360 1630 1310 1490 1270 1370 1230 1320 1200 1490 1150 1420 1090 1460 1050 1520 990 1540 940 1530 880 1600 840 1680 810 1740 740 1680 740 1630 720 1620 640 1710 600 1710 540 1650 480 1580 470 1610 430 1760 410
% 
\special{pn 8}%
\special{pa 610 1800}%
\special{pa 644 1800}%
\special{pa 678 1802}%
\special{pa 712 1802}%
\special{pa 746 1802}%
\special{pa 778 1800}%
\special{pa 810 1800}%
\special{pa 842 1798}%
\special{pa 872 1796}%
\special{pa 902 1794}%
\special{pa 928 1790}%
\special{pa 956 1784}%
\special{pa 986 1780}%
\special{pa 1018 1774}%
\special{pa 1060 1770}%
\special{pa 1102 1764}%
\special{pa 1120 1752}%
\special{pa 1102 1732}%
\special{pa 1078 1706}%
\special{pa 1086 1684}%
\special{pa 1116 1666}%
\special{pa 1156 1648}%
\special{pa 1190 1630}%
\special{pa 1218 1614}%
\special{pa 1242 1600}%
\special{pa 1268 1586}%
\special{pa 1300 1574}%
\special{pa 1340 1564}%
\special{pa 1378 1554}%
\special{pa 1400 1544}%
\special{pa 1388 1530}%
\special{pa 1354 1510}%
\special{pa 1320 1490}%
\special{pa 1312 1466}%
\special{pa 1336 1446}%
\special{pa 1372 1438}%
\special{pa 1404 1440}%
\special{pa 1434 1438}%
\special{pa 1464 1420}%
\special{pa 1482 1390}%
\special{pa 1464 1370}%
\special{pa 1432 1360}%
\special{pa 1426 1354}%
\special{pa 1452 1346}%
\special{pa 1494 1340}%
\special{pa 1544 1332}%
\special{pa 1592 1324}%
\special{pa 1624 1314}%
\special{pa 1630 1306}%
\special{pa 1612 1296}%
\special{pa 1578 1288}%
\special{pa 1534 1280}%
\special{pa 1492 1270}%
\special{pa 1458 1262}%
\special{pa 1428 1254}%
\special{pa 1402 1244}%
\special{pa 1370 1230}%
\special{pa 1336 1214}%
\special{pa 1320 1198}%
\special{pa 1338 1188}%
\special{pa 1378 1180}%
\special{pa 1428 1174}%
\special{pa 1470 1166}%
\special{pa 1490 1156}%
\special{pa 1478 1140}%
\special{pa 1446 1118}%
\special{pa 1422 1096}%
\special{pa 1430 1072}%
\special{pa 1460 1050}%
\special{pa 1488 1028}%
\special{pa 1508 1006}%
\special{pa 1526 980}%
\special{pa 1540 950}%
\special{pa 1536 918}%
\special{pa 1530 886}%
\special{pa 1544 864}%
\special{pa 1574 848}%
\special{pa 1608 838}%
\special{pa 1638 830}%
\special{pa 1668 818}%
\special{pa 1700 796}%
\special{pa 1730 768}%
\special{pa 1742 744}%
\special{pa 1724 738}%
\special{pa 1688 740}%
\special{pa 1652 736}%
\special{pa 1628 716}%
\special{pa 1614 684}%
\special{pa 1616 650}%
\special{pa 1636 630}%
\special{pa 1670 620}%
\special{pa 1700 608}%
\special{pa 1716 586}%
\special{pa 1714 552}%
\special{pa 1702 518}%
\special{pa 1682 494}%
\special{pa 1650 480}%
\special{pa 1610 478}%
\special{pa 1582 474}%
\special{pa 1590 452}%
\special{pa 1618 424}%
\special{pa 1650 408}%
\special{pa 1682 402}%
\special{pa 1714 402}%
\special{pa 1746 408}%
\special{pa 1760 410}%
\special{sp}%
% SPLINE 2 1 3 0
% 43 610 1800 530 1770 520 1740 540 1710 480 1670 480 1640 510 1610 570 1590 610 1550 550 1540 530 1520 530 1490 500 1480 560 1420 610 1400 590 1360 550 1360 460 1340 420 1290 460 1270 430 1220 450 1180 520 1150 560 1150 610 1120 560 1100 500 1100 400 1070 400 1020 440 1000 340 920 300 860 300 800 380 750 390 710 340 700 320 640 410 600 390 550 330 530 240 490 240 450 180 400
% 
\special{pn 8}%
\special{pa 610 1800}%
\special{pa 578 1796}%
\special{pa 548 1784}%
\special{pa 524 1762}%
\special{pa 526 1730}%
\special{pa 540 1706}%
\special{pa 514 1692}%
\special{pa 482 1672}%
\special{pa 480 1640}%
\special{pa 502 1616}%
\special{pa 530 1604}%
\special{pa 562 1594}%
\special{pa 594 1574}%
\special{pa 610 1552}%
\special{pa 588 1546}%
\special{pa 548 1540}%
\special{pa 530 1514}%
\special{pa 526 1486}%
\special{pa 496 1478}%
\special{pa 506 1454}%
\special{pa 542 1428}%
\special{pa 584 1414}%
\special{pa 610 1402}%
\special{pa 602 1372}%
\special{pa 574 1356}%
\special{pa 542 1362}%
\special{pa 512 1362}%
\special{pa 482 1352}%
\special{pa 446 1332}%
\special{pa 420 1304}%
\special{pa 430 1284}%
\special{pa 460 1270}%
\special{pa 448 1246}%
\special{pa 430 1216}%
\special{pa 444 1186}%
\special{pa 472 1164}%
\special{pa 502 1152}%
\special{pa 532 1150}%
\special{pa 566 1150}%
\special{pa 602 1134}%
\special{pa 608 1116}%
\special{pa 576 1102}%
\special{pa 538 1098}%
\special{pa 508 1100}%
\special{pa 476 1102}%
\special{pa 440 1098}%
\special{pa 412 1084}%
\special{pa 394 1054}%
\special{pa 398 1022}%
\special{pa 428 1006}%
\special{pa 444 994}%
\special{pa 430 976}%
\special{pa 396 956}%
\special{pa 356 932}%
\special{pa 326 908}%
\special{pa 308 882}%
\special{pa 298 852}%
\special{pa 296 818}%
\special{pa 308 792}%
\special{pa 334 778}%
\special{pa 366 764}%
\special{pa 390 736}%
\special{pa 388 708}%
\special{pa 356 704}%
\special{pa 326 688}%
\special{pa 316 652}%
\special{pa 334 628}%
\special{pa 370 618}%
\special{pa 402 610}%
\special{pa 410 586}%
\special{pa 394 554}%
\special{pa 368 538}%
\special{pa 336 532}%
\special{pa 298 526}%
\special{pa 264 516}%
\special{pa 242 496}%
\special{pa 240 466}%
\special{pa 234 436}%
\special{pa 210 416}%
\special{pa 180 400}%
\special{sp 0.070}%
\end{picture}%
\end{center}
\caption{Illustrations of 
the noncolliding generalized meanders 
${\bf X}^{(\nu,\kappa)}(t)$ with $\nu\ge 0$ 
in the left picture, and $-1< \nu< 0$ in the right picture.
When $-1< \nu < 0$, a bi-valued process is assigned to
describe the motion of the leftmost particle.}
\end{figure}
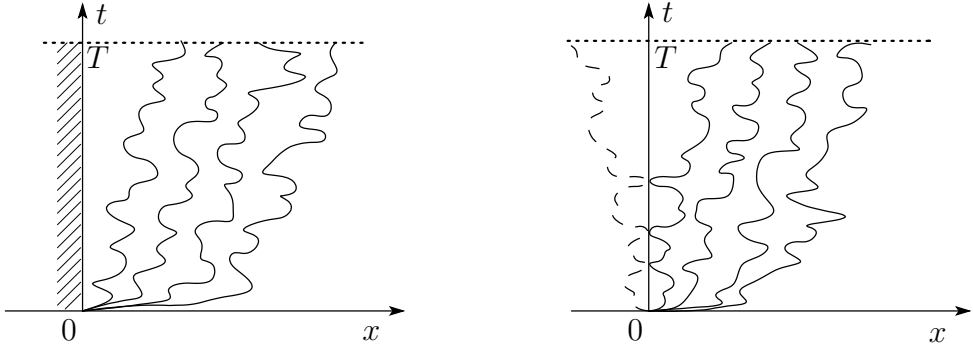
%%%%%%%%%%%%%%%%%%%%%%%%%%%%%%%%%%%%%%%%%%

Consider the noncolliding generalized meander,
when all $N$ particles start from the origin. 
Define $f_N^{(\nu)}(t,{\y}|{\x}) \equiv f_N^{(\nu,0)}(0,{\x}; t, {\y})$.
Then we see
\begin{equation*}
f_N^{(\nu)}(t,{\y}|{\x})
=t^{-N}\prod_{i=1}^N \left( \frac{y_i^{\nu+1}}{x_i^\nu}\right)
\exp\left\{-\frac{1}{2t}\sum_{i=1}^N(x_i^2+y_i^2)\right\}
\det_{1 \leq i,j\leq N}
\left(I_\nu\left(\frac{x_i y_j}{t}\right)\right)
\end{equation*}
and using it 
$g^{(\nu, \kappa)}_{N,T}(s, {\x}; t, {\y})$ is rewritten as
\begin{equation}
g^{(\nu, \kappa)}_{N,T}(s, {\x}; t, {\y}) 
= \frac{\cNt^{(\nu, \kappa)}(T-t, {\y})}
{\widetilde{\mathcal{N}}_{N}^{(\nu, \kappa)}(T-s, {\x})}
f^{(\nu)}(t-s, {\y}|{\x}),
\label{eqn:gNnk1}
\end{equation}
where 
$\widetilde{\mathcal{N}}_{N}^{(\nu, \kappa)}(t, {\x})
= \int_{\W} d {\y} \
f_{N}^{(\nu)}(t, {\y}|{\x}) 
\prod_{i=1}^{N} y_{i}^{-\kappa}$.
Again performing bilinear expansion 
with respect to the Schur functions,
we have obtained the following asymptotics \cite{KT04}:
\begin{eqnarray}
f_{N}^{(\nu)}(t, \y|\x) &\sim&
\frac{t^{-N(N+1+2\nu)/2}}{C^{(\nu)}(N)}
h^{(0)}_{N} \left( \frac{\x}{\sqrt{t}} \right) 
h^{(2\nu+1)}_{N}(\y)
\exp \left( - \frac{|\y|^2}{2t} \right),
\nonumber\\
\widetilde{\mathcal{N}}_{N}^{(\nu,\kappa)}(t,\x) 
&\sim&\frac{t^{-N \kappa/2}C^{(\nu, \kappa)}_N}{C^{(\nu)}(N)}
h^{(0)}_N \left(\frac{\x}{\sqrt{t}} \right),
\qquad \frac{|\x|}{\sqrt{t}} \to 0,
\nonumber
\end{eqnarray}
where 
$C^{(\nu)}(N)=2^{N(N+\nu-1)} \prod_{i=1}^{N} \Gamma(i)\Gamma(i+\nu)$, 
$C^{(\nu, \kappa)}(N)= 2^{N(N+2\nu-\kappa-1)/2} \pi^{-N/2}$
$\prod_{i=1}^{N} \{ \Gamma (i/2)
\Gamma((i+2\nu+1-\kappa)/2)\}$, and
$$
h_{N}^{(\alpha)}(\a)=\prod_{1 \leq i<j \leq N}
(a_{j}^2-a_{i}^2) \prod_{k=1}^{N} a_{k}^{\alpha}.
$$
By the above estimates 
the transition probability density function of 
the noncolliding generalized meander,
when all $N$ particles start
from the origin is determined as
%%%%%%%%%%%%%%%%%%%%%%%%%%%%%%%%%%%%%%%%%%%%%%%%%%
%%%%%%%%%%%%%%%%%%%%%%%%%%%%%%%%%%%%%%%%%%%%%%%%%%
\begin{eqnarray}
&& g_{N,T}^{(\nu, \kappa)}(0, \0; t, {\y}) 
=\frac{T^{N(N+\kappa-1)/2} t^{-N(N+\nu)}}{C^{(\nu, \kappa)}(N)}
\widetilde{{\mathcal{N}}}_{N}^{(\nu, \kappa)}(T-t, {\y})
h_{N}^{(2\nu+1)}(\y)
e^{-|{\y}|^2/2t}, 
\label{eqn:gNnk0}\\
&& \hskip 7cm t \in (0, T], \quad \y \in \W.
\nonumber
\end{eqnarray}
%%%%%%%%%%%%%%%%%%%%%%%%%%%%%%%%%%%%%%%%%%
%%%%%%%%%%%%%%%%%%%%%%%%%%%%%%%%%%%%%%%%%%
\vskip 3mm
\begin{figure}[ch]
\begin{center}
%\input fig_ronsetu_3.1.tex
%WinTpicVersion3.08
\unitlength 0.1in
\begin{picture}( 20.6000, 16.7500)(  0.0000,-18.0500)
% VECTOR 2 0 3 0
% 2 400 1800 2060 1800
% 
\special{pn 8}%
\special{pa 400 1800}%
\special{pa 2060 1800}%
\special{fp}%
\special{sh 1}%
\special{pa 2060 1800}%
\special{pa 1994 1780}%
\special{pa 2008 1800}%
\special{pa 1994 1820}%
\special{pa 2060 1800}%
\special{fp}%
% SPLINE 2 0 3 0
% 36 400 1800 650 1740 710 1720 660 1620 720 1590 830 1530 780 1480 810 1440 890 1390 800 1340 810 1260 890 1240 970 1190 930 1150 1000 1080 970 1010 910 980 870 890 910 830 1020 810 1100 790 1090 750 1050 750 1040 690 1090 680 1110 660 1150 620 1090 570 1120 530 1110 480 1030 470 1090 420 1120 300 1170 270 1210 200 1210 200
% 
\special{pn 8}%
\special{pa 400 1800}%
\special{pa 428 1792}%
\special{pa 454 1784}%
\special{pa 482 1776}%
\special{pa 512 1768}%
\special{pa 542 1760}%
\special{pa 576 1754}%
\special{pa 610 1746}%
\special{pa 650 1740}%
\special{pa 688 1734}%
\special{pa 710 1720}%
\special{pa 704 1696}%
\special{pa 682 1664}%
\special{pa 662 1634}%
\special{pa 666 1610}%
\special{pa 696 1598}%
\special{pa 738 1586}%
\special{pa 782 1572}%
\special{pa 816 1558}%
\special{pa 832 1540}%
\special{pa 820 1520}%
\special{pa 792 1498}%
\special{pa 780 1474}%
\special{pa 800 1448}%
\special{pa 836 1426}%
\special{pa 872 1408}%
\special{pa 890 1394}%
\special{pa 880 1382}%
\special{pa 846 1370}%
\special{pa 812 1352}%
\special{pa 792 1324}%
\special{pa 794 1288}%
\special{pa 812 1260}%
\special{pa 836 1248}%
\special{pa 870 1244}%
\special{pa 910 1236}%
\special{pa 948 1222}%
\special{pa 972 1202}%
\special{pa 960 1180}%
\special{pa 934 1156}%
\special{pa 936 1136}%
\special{pa 964 1114}%
\special{pa 994 1090}%
\special{pa 1004 1060}%
\special{pa 990 1030}%
\special{pa 964 1006}%
\special{pa 934 992}%
\special{pa 906 978}%
\special{pa 884 952}%
\special{pa 872 920}%
\special{pa 872 888}%
\special{pa 884 858}%
\special{pa 906 834}%
\special{pa 932 818}%
\special{pa 962 812}%
\special{pa 996 810}%
\special{pa 1034 812}%
\special{pa 1072 812}%
\special{pa 1098 798}%
\special{pa 1100 764}%
\special{pa 1076 750}%
\special{pa 1044 746}%
\special{pa 1032 712}%
\special{pa 1048 686}%
\special{pa 1082 682}%
\special{pa 1106 664}%
\special{pa 1136 640}%
\special{pa 1150 618}%
\special{pa 1126 598}%
\special{pa 1094 578}%
\special{pa 1098 556}%
\special{pa 1122 528}%
\special{pa 1122 494}%
\special{pa 1094 474}%
\special{pa 1054 472}%
\special{pa 1030 470}%
\special{pa 1044 458}%
\special{pa 1076 434}%
\special{pa 1100 404}%
\special{pa 1108 370}%
\special{pa 1108 336}%
\special{pa 1116 308}%
\special{pa 1136 290}%
\special{pa 1166 274}%
\special{pa 1188 252}%
\special{pa 1202 222}%
\special{pa 1210 200}%
\special{sp}%
% SPLINE 2 0 3 0
% 36 400 1800 740 1770 990 1760 1110 1730 1280 1670 1230 1610 1160 1520 1170 1500 1250 1500 1310 1420 1230 1380 1330 1330 1390 1320 1470 1290 1480 1240 1420 1200 1500 1180 1500 1130 1330 1070 1360 1040 1420 980 1580 920 1660 900 1660 880 1520 830 1560 780 1670 780 1680 700 1600 640 1620 610 1720 550 1680 450 1820 370 1740 310 1780 250 1810 200
% 
\special{pn 8}%
\special{pa 400 1800}%
\special{pa 432 1796}%
\special{pa 464 1792}%
\special{pa 494 1788}%
\special{pa 526 1784}%
\special{pa 558 1782}%
\special{pa 590 1778}%
\special{pa 622 1776}%
\special{pa 654 1774}%
\special{pa 686 1772}%
\special{pa 718 1770}%
\special{pa 752 1770}%
\special{pa 784 1770}%
\special{pa 818 1770}%
\special{pa 850 1772}%
\special{pa 884 1770}%
\special{pa 916 1770}%
\special{pa 946 1768}%
\special{pa 976 1764}%
\special{pa 1006 1758}%
\special{pa 1034 1750}%
\special{pa 1064 1742}%
\special{pa 1098 1734}%
\special{pa 1136 1726}%
\special{pa 1178 1720}%
\special{pa 1218 1712}%
\special{pa 1252 1702}%
\special{pa 1274 1688}%
\special{pa 1280 1670}%
\special{pa 1268 1646}%
\special{pa 1244 1622}%
\special{pa 1212 1596}%
\special{pa 1184 1574}%
\special{pa 1164 1548}%
\special{pa 1160 1520}%
\special{pa 1178 1496}%
\special{pa 1208 1496}%
\special{pa 1244 1502}%
\special{pa 1280 1488}%
\special{pa 1306 1460}%
\special{pa 1314 1428}%
\special{pa 1294 1406}%
\special{pa 1258 1392}%
\special{pa 1232 1382}%
\special{pa 1236 1368}%
\special{pa 1266 1352}%
\special{pa 1306 1338}%
\special{pa 1344 1328}%
\special{pa 1376 1322}%
\special{pa 1406 1318}%
\special{pa 1438 1312}%
\special{pa 1466 1294}%
\special{pa 1484 1266}%
\special{pa 1476 1234}%
\special{pa 1442 1212}%
\special{pa 1420 1200}%
\special{pa 1440 1198}%
\special{pa 1482 1190}%
\special{pa 1508 1170}%
\special{pa 1504 1136}%
\special{pa 1476 1110}%
\special{pa 1436 1098}%
\special{pa 1394 1092}%
\special{pa 1356 1088}%
\special{pa 1332 1082}%
\special{pa 1334 1066}%
\special{pa 1358 1044}%
\special{pa 1382 1020}%
\special{pa 1404 996}%
\special{pa 1426 976}%
\special{pa 1450 958}%
\special{pa 1476 942}%
\special{pa 1506 932}%
\special{pa 1540 924}%
\special{pa 1576 920}%
\special{pa 1616 920}%
\special{pa 1648 914}%
\special{pa 1662 890}%
\special{pa 1648 864}%
\special{pa 1614 854}%
\special{pa 1572 850}%
\special{pa 1536 844}%
\special{pa 1520 830}%
\special{pa 1532 802}%
\special{pa 1562 780}%
\special{pa 1596 778}%
\special{pa 1630 784}%
\special{pa 1660 784}%
\special{pa 1682 766}%
\special{pa 1690 732}%
\special{pa 1678 698}%
\special{pa 1648 676}%
\special{pa 1616 660}%
\special{pa 1600 642}%
\special{pa 1616 614}%
\special{pa 1648 594}%
\special{pa 1684 580}%
\special{pa 1712 566}%
\special{pa 1720 546}%
\special{pa 1708 518}%
\special{pa 1688 486}%
\special{pa 1680 456}%
\special{pa 1696 430}%
\special{pa 1730 412}%
\special{pa 1772 398}%
\special{pa 1806 386}%
\special{pa 1822 374}%
\special{pa 1806 358}%
\special{pa 1772 340}%
\special{pa 1744 318}%
\special{pa 1744 294}%
\special{pa 1766 266}%
\special{pa 1790 240}%
\special{pa 1804 212}%
\special{pa 1810 200}%
\special{sp}%
% SPLINE 2 0 3 0
% 36 400 1800 810 1750 860 1730 840 1690 860 1670 1080 1620 1070 1550 1010 1530 1020 1450 980 1390 1070 1330 1180 1320 1180 1280 1160 1180 1130 1140 1200 1130 1230 1090 1110 1020 1110 940 1240 910 1240 870 1200 820 1260 780 1310 700 1260 670 1410 630 1550 570 1440 530 1470 490 1530 460 1370 430 1420 360 1510 360 1440 290 1510 200 1510 200
% 
\special{pn 8}%
\special{pa 400 1800}%
\special{pa 426 1798}%
\special{pa 452 1796}%
\special{pa 478 1792}%
\special{pa 504 1790}%
\special{pa 532 1786}%
\special{pa 562 1784}%
\special{pa 592 1780}%
\special{pa 624 1776}%
\special{pa 656 1772}%
\special{pa 692 1766}%
\special{pa 730 1762}%
\special{pa 772 1756}%
\special{pa 814 1750}%
\special{pa 852 1740}%
\special{pa 856 1720}%
\special{pa 840 1690}%
\special{pa 864 1668}%
\special{pa 900 1658}%
\special{pa 938 1654}%
\special{pa 976 1656}%
\special{pa 1010 1656}%
\special{pa 1040 1652}%
\special{pa 1066 1642}%
\special{pa 1082 1618}%
\special{pa 1086 1584}%
\special{pa 1074 1554}%
\special{pa 1042 1540}%
\special{pa 1012 1532}%
\special{pa 1010 1508}%
\special{pa 1020 1472}%
\special{pa 1016 1438}%
\special{pa 996 1412}%
\special{pa 980 1388}%
\special{pa 990 1364}%
\special{pa 1018 1344}%
\special{pa 1060 1332}%
\special{pa 1102 1330}%
\special{pa 1142 1334}%
\special{pa 1170 1332}%
\special{pa 1182 1314}%
\special{pa 1180 1280}%
\special{pa 1180 1246}%
\special{pa 1176 1216}%
\special{pa 1164 1186}%
\special{pa 1140 1156}%
\special{pa 1132 1138}%
\special{pa 1164 1136}%
\special{pa 1204 1130}%
\special{pa 1228 1104}%
\special{pa 1226 1080}%
\special{pa 1204 1064}%
\special{pa 1170 1054}%
\special{pa 1134 1040}%
\special{pa 1106 1016}%
\special{pa 1096 980}%
\special{pa 1104 948}%
\special{pa 1128 930}%
\special{pa 1164 928}%
\special{pa 1202 928}%
\special{pa 1232 920}%
\special{pa 1244 896}%
\special{pa 1236 862}%
\special{pa 1210 836}%
\special{pa 1202 816}%
\special{pa 1224 800}%
\special{pa 1260 780}%
\special{pa 1294 754}%
\special{pa 1314 726}%
\special{pa 1310 700}%
\special{pa 1282 680}%
\special{pa 1258 668}%
\special{pa 1260 658}%
\special{pa 1284 652}%
\special{pa 1326 644}%
\special{pa 1376 636}%
\special{pa 1430 626}%
\special{pa 1482 614}%
\special{pa 1522 600}%
\special{pa 1550 586}%
\special{pa 1556 574}%
\special{pa 1534 566}%
\special{pa 1498 558}%
\special{pa 1460 548}%
\special{pa 1440 532}%
\special{pa 1452 506}%
\special{pa 1484 482}%
\special{pa 1518 466}%
\special{pa 1530 458}%
\special{pa 1512 458}%
\special{pa 1474 456}%
\special{pa 1428 454}%
\special{pa 1388 444}%
\special{pa 1368 426}%
\special{pa 1376 396}%
\special{pa 1402 368}%
\special{pa 1442 356}%
\special{pa 1484 358}%
\special{pa 1508 362}%
\special{pa 1504 352}%
\special{pa 1476 330}%
\special{pa 1448 302}%
\special{pa 1438 274}%
\special{pa 1450 248}%
\special{pa 1478 222}%
\special{pa 1510 200}%
\special{sp}%
% SPLINE 2 0 3 0
% 29 400 1800 550 1710 450 1660 610 1590 610 1500 540 1470 670 1380 560 1310 520 1260 560 1230 640 1210 640 1170 790 1090 750 1020 650 990 610 920 680 880 720 860 700 810 700 740 830 690 760 620 820 580 760 520 930 460 880 330 920 250 910 200 910 200
% 
\special{pn 8}%
\special{pa 400 1800}%
\special{pa 450 1784}%
\special{pa 494 1768}%
\special{pa 530 1752}%
\special{pa 552 1734}%
\special{pa 556 1718}%
\special{pa 538 1702}%
\special{pa 504 1686}%
\special{pa 470 1672}%
\special{pa 450 1660}%
\special{pa 456 1652}%
\special{pa 482 1644}%
\special{pa 520 1636}%
\special{pa 562 1622}%
\special{pa 598 1602}%
\special{pa 622 1572}%
\special{pa 628 1538}%
\special{pa 616 1506}%
\special{pa 586 1488}%
\special{pa 552 1476}%
\special{pa 540 1466}%
\special{pa 556 1450}%
\special{pa 592 1432}%
\special{pa 632 1410}%
\special{pa 662 1390}%
\special{pa 672 1372}%
\special{pa 656 1356}%
\special{pa 624 1342}%
\special{pa 584 1324}%
\special{pa 548 1302}%
\special{pa 524 1276}%
\special{pa 524 1252}%
\special{pa 554 1232}%
\special{pa 592 1226}%
\special{pa 626 1222}%
\special{pa 642 1204}%
\special{pa 640 1172}%
\special{pa 656 1146}%
\special{pa 688 1132}%
\special{pa 726 1124}%
\special{pa 764 1114}%
\special{pa 788 1098}%
\special{pa 790 1072}%
\special{pa 774 1042}%
\special{pa 746 1018}%
\special{pa 716 1008}%
\special{pa 684 1002}%
\special{pa 654 992}%
\special{pa 626 968}%
\special{pa 612 938}%
\special{pa 614 912}%
\special{pa 638 894}%
\special{pa 674 882}%
\special{pa 710 870}%
\special{pa 720 850}%
\special{pa 704 818}%
\special{pa 690 784}%
\special{pa 694 752}%
\special{pa 718 728}%
\special{pa 756 716}%
\special{pa 798 708}%
\special{pa 826 700}%
\special{pa 828 686}%
\special{pa 800 664}%
\special{pa 768 640}%
\special{pa 762 618}%
\special{pa 792 600}%
\special{pa 820 582}%
\special{pa 810 560}%
\special{pa 778 536}%
\special{pa 758 516}%
\special{pa 770 504}%
\special{pa 802 496}%
\special{pa 848 488}%
\special{pa 892 480}%
\special{pa 924 468}%
\special{pa 934 448}%
\special{pa 922 420}%
\special{pa 902 388}%
\special{pa 884 354}%
\special{pa 882 324}%
\special{pa 896 294}%
\special{pa 914 266}%
\special{pa 920 236}%
\special{pa 912 206}%
\special{pa 910 200}%
\special{sp}%
% STR 2 0 3 0
% 3 290 1850 290 1950 2 0
% $0$
\put(2.9000,-19.5000){\makebox(0,0)[lb]{$0$}}%
% STR 2 0 3 0
% 3 1850 1850 1850 1950 2 0
% $x$
\put(18.5000,-19.5000){\makebox(0,0)[lb]{$x$}}%
% STR 2 0 3 0
% 3 470 200 470 300 2 0
% $t$
\put(4.7000,-3.0000){\makebox(0,0)[lb]{$t$}}%
% BOX 2 5 3 0
% 2 260 210 400 1790
% 
\special{pn 8}%
\special{pa 260 210}%
\special{pa 400 210}%
\special{pa 400 1790}%
\special{pa 260 1790}%
\special{pa 260 210}%
\special{ip}%
% LINE 2 0 3 0
% 2 0 1800 400 1800
% 
\special{pn 8}%
\special{pa 0 1800}%
\special{pa 400 1800}%
\special{fp}%
% VECTOR 2 0 3 0
% 2 400 1800 400 200
% 
\special{pn 8}%
\special{pa 400 1800}%
\special{pa 400 200}%
\special{fp}%
\special{sh 1}%
\special{pa 400 200}%
\special{pa 380 268}%
\special{pa 400 254}%
\special{pa 420 268}%
\special{pa 400 200}%
\special{fp}%
% LINE 3 0 3 0
% 56 400 500 260 640 400 560 260 700 400 620 260 760 400 680 260 820 400 740 260 880 400 800 260 940 400 860 260 1000 400 920 260 1060 400 980 260 1120 400 1040 260 1180 400 1100 260 1240 400 1160 260 1300 400 1220 260 1360 400 1280 260 1420 400 1340 260 1480 400 1400 260 1540 400 1460 260 1600 400 1520 260 1660 400 1580 260 1720 400 1640 260 1780 400 1700 310 1790 400 1760 370 1790 400 440 260 580 400 380 260 520 400 320 260 460 400 260 260 400 380 220 260 340 330 210 260 280
% 
\special{pn 4}%
\special{pa 400 500}%
\special{pa 260 640}%
\special{fp}%
\special{pa 400 560}%
\special{pa 260 700}%
\special{fp}%
\special{pa 400 620}%
\special{pa 260 760}%
\special{fp}%
\special{pa 400 680}%
\special{pa 260 820}%
\special{fp}%
\special{pa 400 740}%
\special{pa 260 880}%
\special{fp}%
\special{pa 400 800}%
\special{pa 260 940}%
\special{fp}%
\special{pa 400 860}%
\special{pa 260 1000}%
\special{fp}%
\special{pa 400 920}%
\special{pa 260 1060}%
\special{fp}%
\special{pa 400 980}%
\special{pa 260 1120}%
\special{fp}%
\special{pa 400 1040}%
\special{pa 260 1180}%
\special{fp}%
\special{pa 400 1100}%
\special{pa 260 1240}%
\special{fp}%
\special{pa 400 1160}%
\special{pa 260 1300}%
\special{fp}%
\special{pa 400 1220}%
\special{pa 260 1360}%
\special{fp}%
\special{pa 400 1280}%
\special{pa 260 1420}%
\special{fp}%
\special{pa 400 1340}%
\special{pa 260 1480}%
\special{fp}%
\special{pa 400 1400}%
\special{pa 260 1540}%
\special{fp}%
\special{pa 400 1460}%
\special{pa 260 1600}%
\special{fp}%
\special{pa 400 1520}%
\special{pa 260 1660}%
\special{fp}%
\special{pa 400 1580}%
\special{pa 260 1720}%
\special{fp}%
\special{pa 400 1640}%
\special{pa 260 1780}%
\special{fp}%
\special{pa 400 1700}%
\special{pa 310 1790}%
\special{fp}%
\special{pa 400 1760}%
\special{pa 370 1790}%
\special{fp}%
\special{pa 400 440}%
\special{pa 260 580}%
\special{fp}%
\special{pa 400 380}%
\special{pa 260 520}%
\special{fp}%
\special{pa 400 320}%
\special{pa 260 460}%
\special{fp}%
\special{pa 400 260}%
\special{pa 260 400}%
\special{fp}%
\special{pa 380 220}%
\special{pa 260 340}%
\special{fp}%
\special{pa 330 210}%
\special{pa 260 280}%
\special{fp}%
\end{picture}%
\hspace*{20mm}
%\input fig_ronsetu_3.2.tex
%WinTpicVersion3.08
\unitlength 0.1in
\begin{picture}( 21.1700, 16.8500)(  2.0000,-18.6500)
% VECTOR 2 0 3 0
% 2 657 1860 2317 1860
% 
\special{pn 8}%
\special{pa 658 1860}%
\special{pa 2318 1860}%
\special{fp}%
\special{sh 1}%
\special{pa 2318 1860}%
\special{pa 2250 1840}%
\special{pa 2264 1860}%
\special{pa 2250 1880}%
\special{pa 2318 1860}%
\special{fp}%
% STR 2 0 3 0
% 3 547 1910 547 2010 2 0
% $0$
\put(5.4700,-20.1000){\makebox(0,0)[lb]{$0$}}%
% STR 2 0 3 0
% 3 2107 1920 2107 2020 2 0
% $x$
\put(21.0700,-20.2000){\makebox(0,0)[lb]{$x$}}%
% STR 2 0 3 0
% 3 727 250 727 350 2 0
% $t$
\put(7.2700,-3.5000){\makebox(0,0)[lb]{$t$}}%
% LINE 2 0 3 0
% 2 257 1860 657 1860
% 
\special{pn 8}%
\special{pa 258 1860}%
\special{pa 658 1860}%
\special{fp}%
% VECTOR 2 0 3 0
% 2 657 1860 657 260
% 
\special{pn 8}%
\special{pa 658 1860}%
\special{pa 658 260}%
\special{fp}%
\special{sh 1}%
\special{pa 658 260}%
\special{pa 638 328}%
\special{pa 658 314}%
\special{pa 678 328}%
\special{pa 658 260}%
\special{fp}%
% SPLINE 2 0 3 0
% 46 657 1860 737 1830 747 1800 727 1770 787 1730 787 1700 757 1670 697 1650 657 1610 717 1600 737 1580 737 1550 767 1540 707 1480 657 1460 677 1420 717 1420 807 1400 847 1350 807 1330 837 1280 817 1240 747 1210 707 1210 657 1180 707 1160 767 1160 867 1130 867 1080 827 1060 927 980 967 920 967 860 887 810 877 770 927 760 947 700 857 660 877 610 937 590 1027 550 1027 510 1027 360 1107 330 1097 260 1097 260
% 
\special{pn 8}%
\special{pa 658 1860}%
\special{pa 690 1856}%
\special{pa 720 1844}%
\special{pa 744 1822}%
\special{pa 742 1790}%
\special{pa 728 1766}%
\special{pa 754 1752}%
\special{pa 786 1732}%
\special{pa 788 1700}%
\special{pa 766 1676}%
\special{pa 738 1664}%
\special{pa 706 1654}%
\special{pa 674 1634}%
\special{pa 658 1612}%
\special{pa 680 1606}%
\special{pa 720 1600}%
\special{pa 738 1574}%
\special{pa 742 1546}%
\special{pa 772 1538}%
\special{pa 762 1514}%
\special{pa 726 1488}%
\special{pa 684 1474}%
\special{pa 658 1462}%
\special{pa 666 1432}%
\special{pa 694 1416}%
\special{pa 726 1422}%
\special{pa 756 1422}%
\special{pa 786 1412}%
\special{pa 822 1392}%
\special{pa 848 1364}%
\special{pa 838 1344}%
\special{pa 808 1330}%
\special{pa 820 1306}%
\special{pa 838 1276}%
\special{pa 824 1246}%
\special{pa 796 1224}%
\special{pa 766 1212}%
\special{pa 736 1210}%
\special{pa 702 1210}%
\special{pa 666 1194}%
\special{pa 660 1176}%
\special{pa 692 1162}%
\special{pa 730 1158}%
\special{pa 760 1160}%
\special{pa 792 1162}%
\special{pa 828 1158}%
\special{pa 856 1144}%
\special{pa 874 1114}%
\special{pa 870 1082}%
\special{pa 840 1066}%
\special{pa 824 1054}%
\special{pa 838 1036}%
\special{pa 872 1016}%
\special{pa 912 992}%
\special{pa 942 968}%
\special{pa 960 942}%
\special{pa 970 912}%
\special{pa 972 878}%
\special{pa 960 852}%
\special{pa 934 838}%
\special{pa 902 824}%
\special{pa 878 796}%
\special{pa 880 768}%
\special{pa 912 764}%
\special{pa 942 748}%
\special{pa 952 712}%
\special{pa 934 688}%
\special{pa 898 678}%
\special{pa 866 670}%
\special{pa 856 646}%
\special{pa 874 614}%
\special{pa 900 598}%
\special{pa 932 592}%
\special{pa 968 586}%
\special{pa 1002 574}%
\special{pa 1024 556}%
\special{pa 1030 526}%
\special{pa 1024 492}%
\special{pa 1014 454}%
\special{pa 1008 416}%
\special{pa 1010 386}%
\special{pa 1024 364}%
\special{pa 1054 352}%
\special{pa 1088 342}%
\special{pa 1110 326}%
\special{pa 1110 298}%
\special{pa 1098 262}%
\special{pa 1098 260}%
\special{sp}%
% SPLINE 2 0 3 0
% 40 657 1860 817 1830 877 1760 907 1670 827 1640 887 1600 987 1570 967 1540 947 1470 987 1430 897 1380 907 1320 907 1300 1067 1260 1097 1170 1097 1150 1057 1140 1017 1140 1057 1090 1097 1080 1157 1060 1117 1040 1067 1030 1057 990 1167 970 1207 950 1227 890 1107 830 1087 790 1157 790 1227 740 1227 700 1187 670 1247 630 1247 590 1177 540 1217 500 1217 360 1307 330 1357 260
% 
\special{pn 8}%
\special{pa 658 1860}%
\special{pa 692 1860}%
\special{pa 728 1858}%
\special{pa 760 1852}%
\special{pa 790 1844}%
\special{pa 816 1832}%
\special{pa 838 1814}%
\special{pa 858 1790}%
\special{pa 878 1760}%
\special{pa 900 1726}%
\special{pa 914 1694}%
\special{pa 906 1670}%
\special{pa 874 1656}%
\special{pa 838 1648}%
\special{pa 828 1634}%
\special{pa 852 1614}%
\special{pa 894 1598}%
\special{pa 938 1590}%
\special{pa 974 1586}%
\special{pa 988 1572}%
\special{pa 974 1548}%
\special{pa 952 1516}%
\special{pa 942 1486}%
\special{pa 958 1458}%
\special{pa 984 1436}%
\special{pa 980 1420}%
\special{pa 944 1408}%
\special{pa 908 1392}%
\special{pa 896 1368}%
\special{pa 904 1334}%
\special{pa 908 1302}%
\special{pa 920 1282}%
\special{pa 950 1274}%
\special{pa 990 1274}%
\special{pa 1030 1272}%
\special{pa 1062 1264}%
\special{pa 1082 1242}%
\special{pa 1092 1212}%
\special{pa 1096 1176}%
\special{pa 1094 1146}%
\special{pa 1060 1140}%
\special{pa 1024 1142}%
\special{pa 1020 1126}%
\special{pa 1048 1098}%
\special{pa 1078 1084}%
\special{pa 1114 1078}%
\special{pa 1150 1068}%
\special{pa 1152 1054}%
\special{pa 1116 1040}%
\special{pa 1082 1038}%
\special{pa 1060 1018}%
\special{pa 1058 986}%
\special{pa 1080 972}%
\special{pa 1116 970}%
\special{pa 1154 972}%
\special{pa 1188 964}%
\special{pa 1212 944}%
\special{pa 1228 912}%
\special{pa 1224 884}%
\special{pa 1204 870}%
\special{pa 1172 860}%
\special{pa 1138 850}%
\special{pa 1106 828}%
\special{pa 1088 798}%
\special{pa 1100 784}%
\special{pa 1136 790}%
\special{pa 1176 788}%
\special{pa 1206 774}%
\special{pa 1224 746}%
\special{pa 1232 714}%
\special{pa 1212 688}%
\special{pa 1188 670}%
\special{pa 1206 656}%
\special{pa 1240 638}%
\special{pa 1254 606}%
\special{pa 1234 578}%
\special{pa 1200 558}%
\special{pa 1178 544}%
\special{pa 1192 526}%
\special{pa 1218 500}%
\special{pa 1228 466}%
\special{pa 1222 430}%
\special{pa 1214 396}%
\special{pa 1214 368}%
\special{pa 1232 352}%
\special{pa 1264 344}%
\special{pa 1298 334}%
\special{pa 1326 318}%
\special{pa 1344 292}%
\special{pa 1358 260}%
\special{sp}%
% SPLINE 2 0 3 0
% 35 657 1860 917 1840 1037 1810 997 1780 1037 1750 1097 1720 1137 1620 1237 1570 1277 1530 1157 1480 1117 1420 1217 1390 1307 1330 1287 1300 1217 1190 1217 1130 1267 1040 1437 1010 1387 960 1417 920 1467 870 1447 800 1367 770 1407 730 1477 700 1497 660 1447 620 1407 580 1497 560 1527 550 1407 400 1467 390 1517 350 1617 300 1627 260
% 
\special{pn 8}%
\special{pa 658 1860}%
\special{pa 680 1858}%
\special{pa 704 1856}%
\special{pa 730 1854}%
\special{pa 758 1852}%
\special{pa 788 1850}%
\special{pa 824 1848}%
\special{pa 864 1844}%
\special{pa 910 1840}%
\special{pa 962 1838}%
\special{pa 1008 1832}%
\special{pa 1038 1824}%
\special{pa 1038 1810}%
\special{pa 1008 1792}%
\special{pa 1000 1772}%
\special{pa 1032 1752}%
\special{pa 1068 1740}%
\special{pa 1092 1726}%
\special{pa 1108 1700}%
\special{pa 1118 1668}%
\special{pa 1128 1636}%
\special{pa 1146 1612}%
\special{pa 1170 1598}%
\special{pa 1202 1586}%
\special{pa 1234 1572}%
\special{pa 1264 1550}%
\special{pa 1278 1528}%
\special{pa 1266 1514}%
\special{pa 1234 1506}%
\special{pa 1194 1496}%
\special{pa 1154 1478}%
\special{pa 1128 1452}%
\special{pa 1118 1424}%
\special{pa 1128 1408}%
\special{pa 1158 1400}%
\special{pa 1198 1394}%
\special{pa 1242 1384}%
\special{pa 1280 1370}%
\special{pa 1304 1350}%
\special{pa 1306 1326}%
\special{pa 1286 1300}%
\special{pa 1264 1272}%
\special{pa 1244 1246}%
\special{pa 1228 1220}%
\special{pa 1218 1192}%
\special{pa 1216 1160}%
\special{pa 1218 1128}%
\special{pa 1224 1094}%
\special{pa 1238 1064}%
\special{pa 1264 1042}%
\special{pa 1302 1030}%
\special{pa 1348 1026}%
\special{pa 1392 1024}%
\special{pa 1426 1022}%
\special{pa 1440 1016}%
\special{pa 1428 1002}%
\special{pa 1398 980}%
\special{pa 1388 952}%
\special{pa 1410 928}%
\special{pa 1438 906}%
\special{pa 1460 884}%
\special{pa 1470 854}%
\special{pa 1464 822}%
\special{pa 1440 796}%
\special{pa 1402 786}%
\special{pa 1372 776}%
\special{pa 1372 758}%
\special{pa 1400 734}%
\special{pa 1434 720}%
\special{pa 1462 710}%
\special{pa 1488 690}%
\special{pa 1498 660}%
\special{pa 1480 638}%
\special{pa 1446 620}%
\special{pa 1416 596}%
\special{pa 1408 576}%
\special{pa 1432 568}%
\special{pa 1472 564}%
\special{pa 1512 558}%
\special{pa 1532 544}%
\special{pa 1522 520}%
\special{pa 1492 490}%
\special{pa 1456 458}%
\special{pa 1422 430}%
\special{pa 1404 408}%
\special{pa 1416 398}%
\special{pa 1448 394}%
\special{pa 1482 384}%
\special{pa 1504 362}%
\special{pa 1532 342}%
\special{pa 1564 332}%
\special{pa 1594 320}%
\special{pa 1618 300}%
\special{pa 1626 270}%
\special{pa 1628 260}%
\special{sp}%
% SPLINE 2 0 3 0
% 39 657 1860 967 1850 1097 1830 1167 1810 1127 1770 1217 1700 1357 1630 1447 1600 1357 1530 1397 1500 1497 1490 1527 1440 1477 1420 1677 1370 1537 1330 1417 1290 1367 1260 1537 1210 1467 1150 1507 1110 1567 1050 1587 1000 1577 940 1647 900 1727 870 1787 800 1727 800 1677 780 1667 700 1757 660 1757 600 1697 540 1627 530 1657 490 1847 440 1797 390 1857 350 1807 310 1917 260
% 
\special{pn 8}%
\special{pa 658 1860}%
\special{pa 692 1860}%
\special{pa 726 1862}%
\special{pa 760 1862}%
\special{pa 794 1862}%
\special{pa 826 1860}%
\special{pa 858 1860}%
\special{pa 890 1858}%
\special{pa 920 1856}%
\special{pa 948 1854}%
\special{pa 976 1850}%
\special{pa 1004 1844}%
\special{pa 1032 1840}%
\special{pa 1066 1834}%
\special{pa 1106 1830}%
\special{pa 1148 1824}%
\special{pa 1168 1812}%
\special{pa 1148 1792}%
\special{pa 1126 1766}%
\special{pa 1132 1744}%
\special{pa 1162 1726}%
\special{pa 1202 1708}%
\special{pa 1238 1690}%
\special{pa 1264 1674}%
\special{pa 1290 1660}%
\special{pa 1314 1646}%
\special{pa 1346 1634}%
\special{pa 1386 1624}%
\special{pa 1426 1614}%
\special{pa 1448 1604}%
\special{pa 1436 1590}%
\special{pa 1402 1570}%
\special{pa 1368 1550}%
\special{pa 1358 1526}%
\special{pa 1384 1506}%
\special{pa 1418 1498}%
\special{pa 1450 1500}%
\special{pa 1482 1498}%
\special{pa 1512 1480}%
\special{pa 1530 1450}%
\special{pa 1512 1430}%
\special{pa 1478 1420}%
\special{pa 1474 1414}%
\special{pa 1498 1406}%
\special{pa 1542 1400}%
\special{pa 1592 1392}%
\special{pa 1638 1384}%
\special{pa 1670 1374}%
\special{pa 1678 1366}%
\special{pa 1660 1356}%
\special{pa 1624 1348}%
\special{pa 1580 1340}%
\special{pa 1538 1330}%
\special{pa 1504 1322}%
\special{pa 1476 1314}%
\special{pa 1448 1304}%
\special{pa 1416 1290}%
\special{pa 1382 1274}%
\special{pa 1368 1258}%
\special{pa 1386 1248}%
\special{pa 1426 1240}%
\special{pa 1474 1234}%
\special{pa 1516 1226}%
\special{pa 1538 1216}%
\special{pa 1526 1200}%
\special{pa 1492 1178}%
\special{pa 1468 1156}%
\special{pa 1478 1132}%
\special{pa 1508 1110}%
\special{pa 1534 1088}%
\special{pa 1556 1066}%
\special{pa 1574 1040}%
\special{pa 1586 1010}%
\special{pa 1584 978}%
\special{pa 1576 946}%
\special{pa 1590 924}%
\special{pa 1622 908}%
\special{pa 1656 898}%
\special{pa 1686 890}%
\special{pa 1714 878}%
\special{pa 1748 856}%
\special{pa 1776 828}%
\special{pa 1790 804}%
\special{pa 1772 798}%
\special{pa 1734 800}%
\special{pa 1700 796}%
\special{pa 1674 776}%
\special{pa 1660 744}%
\special{pa 1662 710}%
\special{pa 1682 690}%
\special{pa 1716 680}%
\special{pa 1748 668}%
\special{pa 1764 646}%
\special{pa 1760 612}%
\special{pa 1748 578}%
\special{pa 1728 554}%
\special{pa 1696 540}%
\special{pa 1656 538}%
\special{pa 1630 534}%
\special{pa 1634 512}%
\special{pa 1668 484}%
\special{pa 1712 468}%
\special{pa 1760 462}%
\special{pa 1804 460}%
\special{pa 1838 458}%
\special{pa 1854 452}%
\special{pa 1846 438}%
\special{pa 1814 416}%
\special{pa 1798 392}%
\special{pa 1820 374}%
\special{pa 1854 356}%
\special{pa 1848 338}%
\special{pa 1816 318}%
\special{pa 1802 300}%
\special{pa 1820 284}%
\special{pa 1860 272}%
\special{pa 1910 262}%
\special{pa 1918 260}%
\special{sp}%
% SPLINE 2 1 3 0
% 46 655 1860 575 1830 565 1800 585 1770 525 1730 525 1700 555 1670 615 1650 655 1610 595 1600 575 1580 575 1550 545 1540 605 1480 655 1460 635 1420 595 1420 505 1400 465 1350 505 1330 475 1280 495 1240 565 1210 605 1210 655 1180 605 1160 545 1160 445 1130 445 1080 485 1060 385 980 345 920 345 860 425 810 435 770 385 760 365 700 455 660 435 610 375 590 285 550 285 510 285 360 205 330 215 260 215 260
% 
\special{pn 8}%
\special{pa 656 1860}%
\special{pa 622 1856}%
\special{pa 592 1844}%
\special{pa 570 1822}%
\special{pa 572 1790}%
\special{pa 586 1766}%
\special{pa 558 1752}%
\special{pa 526 1732}%
\special{pa 526 1700}%
\special{pa 548 1676}%
\special{pa 576 1664}%
\special{pa 606 1654}%
\special{pa 640 1634}%
\special{pa 656 1612}%
\special{pa 632 1606}%
\special{pa 594 1600}%
\special{pa 576 1574}%
\special{pa 570 1546}%
\special{pa 542 1538}%
\special{pa 550 1514}%
\special{pa 588 1488}%
\special{pa 630 1474}%
\special{pa 654 1462}%
\special{pa 648 1432}%
\special{pa 620 1416}%
\special{pa 588 1422}%
\special{pa 556 1422}%
\special{pa 526 1412}%
\special{pa 492 1392}%
\special{pa 464 1364}%
\special{pa 476 1344}%
\special{pa 506 1330}%
\special{pa 494 1306}%
\special{pa 474 1276}%
\special{pa 490 1246}%
\special{pa 516 1224}%
\special{pa 546 1212}%
\special{pa 576 1210}%
\special{pa 610 1210}%
\special{pa 646 1194}%
\special{pa 652 1176}%
\special{pa 620 1162}%
\special{pa 584 1158}%
\special{pa 554 1160}%
\special{pa 520 1162}%
\special{pa 486 1158}%
\special{pa 456 1144}%
\special{pa 440 1114}%
\special{pa 444 1082}%
\special{pa 474 1066}%
\special{pa 490 1054}%
\special{pa 474 1036}%
\special{pa 440 1016}%
\special{pa 402 992}%
\special{pa 372 968}%
\special{pa 354 942}%
\special{pa 344 912}%
\special{pa 340 878}%
\special{pa 352 852}%
\special{pa 380 838}%
\special{pa 410 824}%
\special{pa 434 796}%
\special{pa 432 768}%
\special{pa 400 764}%
\special{pa 370 748}%
\special{pa 362 712}%
\special{pa 380 688}%
\special{pa 416 678}%
\special{pa 448 670}%
\special{pa 456 646}%
\special{pa 440 614}%
\special{pa 412 598}%
\special{pa 380 592}%
\special{pa 344 586}%
\special{pa 310 574}%
\special{pa 288 556}%
\special{pa 284 526}%
\special{pa 290 492}%
\special{pa 298 454}%
\special{pa 304 416}%
\special{pa 304 386}%
\special{pa 290 364}%
\special{pa 260 352}%
\special{pa 226 342}%
\special{pa 204 326}%
\special{pa 202 298}%
\special{pa 214 262}%
\special{pa 216 260}%
\special{sp 0.070}%
\end{picture}%
\end{center}
\caption{Illustrations of the
noncolliding Bessel processes
$\Y^{(\nu)}(t)$ with $\nu\ge 0$ 
in the left picture, and $-1< \nu< 0$ in the right picture.
When $-1< \nu < 0$, a bi-valued process is assigned to
describe the motion of the leftmost particle.}
\end{figure}
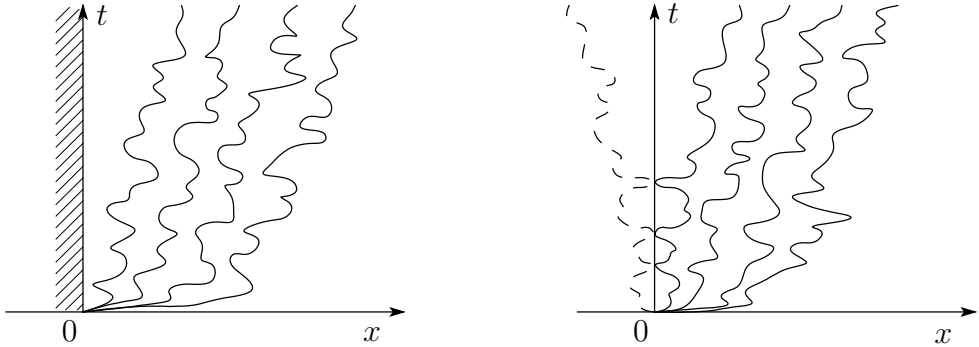
%%%%%%%%%%%%%%%%%%%%%%%%%%%%%%%%%%%%%%%%%%

When the time $T$ becomes infinity, 
the process ${\bf X}^{(\nu,\kappa)}(t)$ converges to 
the temporally homogeneous process $\Y^{(\nu)}(t)$,
whose transition probability density function is given by
\begin{eqnarray}
p^{(\nu)}_N (t, \y|\x)
&=& \frac{h_{N}^{(0)}(\y)}{h_{N}^{(0)}(\x)}
f_{N}^{(\nu)}(t, \y|\x), 
\quad t >0, \quad \x, \y \in \W,
\label{eqn:pNnk0}\\
p^{(\nu)}_N (t, \y| {\bf 0})
&=& \frac{ t^{-N(N+\nu)}}{C^{(\nu)}(N)}
h_{N}^{(\nu+1/2)}(\y)^2
\exp \left( -\frac{|\y|^2}{2t} \right), 
\ t>0, \ \y\in\W.
\nonumber
\end{eqnarray}
See Fig 3.
Since the parameter $\kappa$ controls 
the distribution of the process
when $t \to T$, 
it is irrelevant for
the process $\Y^{(\nu)}(t)$ in which $T \to \infty$.
The process $\Y^{(\nu)}(t)$ is 
the {\bf noncolliding $2(\nu+1)$-dimensional Bessel process},
which is temporally homogeneous 
and solves the following system of 
stochastic differential equations when $\nu \ge -1/2$F
\begin{eqnarray}
\label{eqn:Ynu}
&& Y_i^{(\nu)}(t) = B_i(t)
+\int_0^t \frac{\nu+1/2}{Y_i(s)}ds
+\sum_{\substack{1 \leq j \leq N\\ j \ne i}}
\int_0^t \frac{2Y_i^{(\nu)}(s)}{Y_i^{(\nu)}(s)^2-Y_j^{(\nu)}(s)^2}ds,
\\
&& \hskip 4.5cm 
t \in (0, \infty), \, 1 \leq i \leq N,
\nonumber
\end{eqnarray}
where we impose a reflecting wall at the origin 
when $\nu = -1/2$.
Although solution of the system of equations is not 
necessarily unique in general,
$\Y^{(\nu)}(t)$ can be defined as a unique solution 
such that all coordinates are positive \cite{Mc60}.
By comparing (\ref{eqn:gNnk1}), (\ref{eqn:gNnk0}) 
with (\ref{eqn:pNnk0}),
we have the following equality
$$
P(\X^{(\nu,\kappa)}(\cdot) \in dw) 
=\frac{C^{(\nu)}(N)}{C^{(\nu,\kappa)}(N)}
\frac{T^{N(N+\kappa-1)/2}}{h_N^{(\kappa)}(w(T))}
P(\Y^{(\nu)}(\cdot) \in dw),
$$
which is an extension of Imhof's relation
for the noncolliding Bessel process $\Y^{(\nu)}(t)$ and 
the noncolliding generalized meander $\X^{(\nu, \kappa)}(t)$ 
\cite{KT04}.

%%%%%%%%%%%%%%%%%%%%%%%%%%%%%%%%%%%%%%%%%%%%%%%%%%%%%%%%%%%%
%%%%%%%%%%Section 4%%%%%%%%%%%%%%%%%%%%%%%%%%%%%%%%%%%%%%%%%
%%%%%%%%%%%%%%%%%%%%%%%%%%%%%%%%%%%%%%%%%%%%%%%%%%%%%%%%%%%%
\section{Matrix-valued processes}
%%%%%%%%%%%%%%%%%%%%%%%%%%%%%%%%%%%%%%%%%%%%%%%%%%%%%%%%%%%%
%%%%%%%%%%%%%%%%%%%%%%%%%%%%%%%%%%%%%%%%%%%%%%%%%%%%%%%%%%%%
%%%%%%%%%%%%%%%%%%%%%%%%%%%%%%%%%%%%%%%%%%%%%%%%%%%%%%%%%%%%

%%%%%%%%%%%%%%%%%%%%%%%%%%%%%%%%%%%%%%%%%%%%%%%%%%%%%%
\subsection{Generalized Bru's theorem}
%%%%%%%%%%%%%%%%%%%%%%%%%%%%%%%%%%%%%%%%%%%%%%%%%%%%%%%%%%%%

We denote the space of $N \times N$ Hermitian matrices
by $\mathcal{H}(N)$ and
the spaces of $N \times N$ real symmetric matrices 
by $\mathcal{S}(N)$. For a matrix $A$
we indicate by $^{t}\! A$ the transposed matrix of $A$,
by $\overline{A}$ the complex conjugate of $A$, 
and  by $A^{*}\equiv \ ^{t}\overline{A}$
the adjoint matrix of $A$, respectively.
We denote the unit matrix of size $N \times N$ by $I_N$.
Bru \cite{Bru89, Bru91} studied the eigenvalue processes of
Wishart process, which is an $\H$-valued process,
and derived the stochastic differential equations
for the eigenvalue processes.
The result is generalized to the case that
each element of matrix-valued process,
$\xi_{ij}(t)$, $1\leq i,j \leq N$, 
is a complex-valued continuous
semi-martingale \cite{KT03b, KT04}.
In this section we state this generalized version 
of Bru's theorem and give its applications.

Let 
$\vlambda (t) = (\lambda_1(t),\lambda_2(t),\dots, \lambda_N(t))$ 
be the vector, whose coordinates are eigenvalues of 
$\H$-valued process, $\Xi (t)=(\xi_{ij}(t))_{1\leq i,j \leq N}$, 
with $\lambda_1(t) \leq \lambda_2(t) \leq \cdots \leq \lambda_N(t)$.
Then let $U(t)=(u_{ij}(t))_{1\le i,j\le N}$ 
be a unitary-matrix-valued
process, which diagonalizes $\Xi(t)$, 
$$
U(t)^{*} \Xi (t) U(t) 
= \Lambda (t) = {\rm diag} \Big( \lambda_1(t), \lambda_2(t),
\dots, \lambda_N(t) \Big).
$$
We put
$$
\Gamma_{ij,k\ell}(t)dt =\Big(U(t)^{*} 
d \Xi (t) U(t) \Big)_{ij}
\Big(U(t)^{*} d \Xi (t) U(t) \Big)_{k\ell},
$$
and the bounded variation part of $(U(t)^{*} d \Xi (t) U(t))_{ii}$
is written as $d \Upsilon_i(t)$.
Then we introduce the Markov times
\begin{eqnarray}
&& \sigma_{ij}=\inf\{ t\ge 0 : \lambda_i(t)\not=\lambda_j(t)\},
\nonumber\\
&&
\tau_{ij}= \inf\{t > \sigma_{ij}: \lambda_i(t)=\lambda_j(t)\},
\quad 
\tau = \displaystyle{\min_{1\le i<j \le N}\tau_{ij}}.
\nonumber
\end{eqnarray}
\vskip 3mm
%%%%%%%%%%%theorem%%%%%%%%%%%%%%%%%%%%%%%%%%%%%%%%%%%%%%
%%%%%%%%%%%theorem%%%%%%%%%%%%%%%%%%%%%%%%%%%%%%%%%%%%%%
%%%%%%%%%%%theorem%%%%%%%%%%%%%%%%%%%%%%%%%%%%%%%%%%%%%%
\noindent{\large\bf [Generalized Bru's theorem]}
{\rm (\cite{Bru89, Bru91, KT03b, KT04})}
\; 
Let $\xi_{ij}(t)$, $1\leq i,j \leq N$ be 
complex-valued continuous semi-martingales.
Then the eigenvalue process $\vlambda (t)$ of $\Xi(t)$
solves the following system of stochastic differential equations:
$$
\label{eqn:l=M+J}
d\lambda_i(t) = dM_i(t) + dJ_i(t), 
\quad t \in (0,\tau),
\quad 1 \leq i \leq N.
$$
where $\M(t)=(M_1(t),M_2(t),\dots,M_N(t))$ is
the martingale with 
$d M_i(t) dM_j(t)=\Gamma_{ii,jj}(t)dt$,
and 
$\J(t)=(J_1(t),J_2(t),\dots,J_N(t))$ is
the process with bounded variation given by
$$
d J_i(t)=
\sum_{\substack{1 \leq j \leq N\\ j \ne i}}
\frac{1}{\lambda_i(t) -\lambda_j(t)}
{\bf 1}_{\{\lambda_i(t) \not= \lambda_j(t)\}}
\Gamma_{ij,ji}(t)dt
+ d\Upsilon_i(t).
$$
Here ${\bf 1}_{\{\omega\}}$ denotes an indicator function of
a condition $\omega$. 
%%%%%%%%%%%theorem%%%%%%%%%%%%%%%%%%%%%%%%%%%%%%%%%%%%%%
%%%%%%%%%%%theorem%%%%%%%%%%%%%%%%%%%%%%%%%%%%%%%%%%%%%%
%%%%%%%%%%%theorem%%%%%%%%%%%%%%%%%%%%%%%%%%%%%%%%%%%%%%
\vskip 3mm

Most of the examples shown in the next subsection are
systems such that all particles start
from the origin and rapidly separate from each other
to avoid collision.
In these systems $\sigma_{ij}=0$, $1\le i, j \le N$,
and $\tau=\infty$, that is,
the repulsive forces among particles
are strong enough to prevent any collision.
For instance, in the stochastic differential equations
of Dyson's Brownian motion models (\ref{eqn:Dyson0}),
the repulsive force becomes stronger as 
the parameter $\beta$ becomes larger, 
and it is shown that
$\tau < \infty$, if $0< \beta < 1$, and
$\tau =\infty$, if $\beta \ge 1$ \cite{RS93}.
This corresponds to the fact that
the Bessel process is transient, 
if the dimension $d\ge 2$ ($\nu \ge 0$), and 
it is recurrent, if $0< d < 2$ ($-1< \nu < 0$).

%%%%%%%%%%%theorem%%%%%%%%%%%%%%%%%%%%%%%%%%%%%%%%%%%%%%
%%%%%%%%%%%%%%%%%%%%%%%%%%%%%%%%%%%%%%%%%%%%%%%%

\subsection{Examples}

In this subsection we give examples 
of eigenvalue processes obtained by 
the generalized Bru's theorem.

Let $\nu \in \N_0 \equiv \{0,1,2, \dots\}$, 
$B_{ij}(t)$, $\widetilde{B}_{ij}(t)$, 
$1\le i \leq N+\nu, 1 \leq j \leq N$
be independent one-dimensional Brownian motions, 
and
$s(t)$ and $a(t)$ be $N\times N$ matrices
whose element are given by
$$
s_{ij}(t)
=
\left\{
\begin{array}{ll}
\displaystyle{\frac{1}{\sqrt{2}} B_{ij}(t)}, 
	& \ \mbox{if } i < j,  \\
%& \\
B_{ii}(t), 
	& \ \mbox{if } i=j, \\
%& \\
\displaystyle{\frac{1}{\sqrt{2}} B_{ji}(t)}, 
	& \ \mbox{if } i > j, \\
\end{array}\right. 
\qquad
a_{ij}(t)
=
\left\{
\begin{array}{ll}
\displaystyle{\frac{1}{\sqrt{2}} 
\widetilde{B}_{ij}(t)}, 
	& \ \mbox{if } i < j, \\
%& \\
0, & \ \mbox{if } i=j, \\
%& \\
\displaystyle{-\frac{1}{\sqrt{2}} 
\widetilde{B}_{ji}(t)},
	& \ \mbox{if } i > j, \\
\end{array}\right.
$$
respectively.

%%%%%%%%%%%%%%%%%%%%%%%%%%%%%%%%%%%%%%%%%%%%%%%%%%%%%%%

%\begin{description}
\vskip 0.3cm
\noindent{\bf (i)} \quad {\bf GUE process}

Consider the $\H$-valued process defined by
$\Xi^{\rm GUE}(t) =s(t)+\sqrt{-1}a(t)$,
$t \in [0, \infty)$.
For any fixed $t\in [0,\infty)$, 
$\Xi^{\rm GUE} (t)$ is the $\H$-valued random variable
whose probability density function with respect to 
the volume element $\dH$ of $\H$ is
\begin{equation*}
\mu^{\rm GUE}(H , t)
= \frac{t^{-N^2/2}}{c_1(N)}
\exp \left( - \frac{1}{2t} {\rm Tr}H^2 \right),
\quad H\in \H, 
\end{equation*}
where ${\rm Tr}A$ represents the trace of a matrix $A$,
and $c_1(N)= 2^{N/2} \pi^{N^2/2}$.
We denote the group of $N \times N$
unitary matrices by ${\bf U}(N)$.
The probability $\mu^{\rm GUE}(H, t)\dH $ 
is invariant under any unitary transformation
 $H\to U^{*}HU$ for any $U\in{\bf U}(N)$.
In the random matrix theory,
such a statistical ensemble of $\H$-valued random variables
is called the
{\bf Gaussian unitary ensemble, GUE} \cite{Mehta,N05}. 
The probability density of eigenvalues of GUE
is given by
$$
g^{\rm GUE}(\x , t)
= \frac{t^{-N/2}}{C_1(N)} 
h_N\left( \frac{\x}{\sqrt{t}} \right)^2
\exp \left( - \frac{|\x|^2}{2t}\right)
$$
for $\x=(x_1, x_2, \dots, x_N) \in \WA$ \cite{Mehta,N05}. 
Here $C_1(N)$ is the same constant as $C_1(N)$ in (\ref{eqn:asym}).
By applying the generalized Bru's theorem to 
$\Xi^{\rm GUE}(t)$, we see that
$\vlambda(t),  t \in (0, \infty)$ solves the system of
stochastic differential equations of
Dyson's Brownian motion model (\ref{eqn:Dyson0})  
with $\beta=2$.
In Section 3.2 it was shown that
the noncolliding Brownian motion $\Y(t), t \in (0, \infty)$
solves the same equation.
Then the equivalence in distribution 
of the noncolliding Brownian motion $\Y(t)$
and the eigenvalue process $\vlambda(t)$ of $\Xi^{\rm GUE}(t)$
is established.

%%%%%%%%%%%%%%%%%%%%%%%%%%%%%%%%%%%%%%%%%%%%%%%%%%%%%%%%%
\vskip 0.3cm
\noindent{\bf (ii)} \quad {\bf GOE process}

Consider the $\S$-valued process defined by
$\Xi^{\rm GOE}(t)=s(t)$, $t \in [0, \infty)$.
For any fixed $t\in [0,\infty)$,
$\Xi^{\rm GOE} (t)$ is the $\S$-valued random variable
whose probability density function with respect to 
the volume element $\dS$ of $\S$ is given by
\begin{equation*}
\mu^{\rm GOE}(S,t)
= \frac{t^{-N(N+1)/4}}{c_2(N)}
\exp \left( - \frac{1}{2t} {\rm Tr}S^2 \right),
\quad S\in \S,
\end{equation*}
where $c_2(N)= 2^{N/2} \pi^{N(N+1)/4}$.
We denote the group of $N\times N$ real symmetric matrices
by ${\bf O}(N)$.
The probability $\mu^{\rm GOE}(S, t)\dS$ is invariant
under any orthogonal transformation $S \to ^{t}\!VSV$  
for any $V\in{\bf O}(N)$.
Such a statistical ensemble of $\S$-valued random variables
is called the {\bf Gaussian orthogonal ensemble, GOE}.
The probability density of eigenvalues of GOE 
is given by
\begin{equation}
g^{\rm GOE}(\x , t)
= \frac{t^{-N/2}}{C_2(N)} 
h_{N}\left( \frac{\x}{\sqrt{t}} \right)
\exp \left( - \frac{|\x|^2}{2t}\right)
\label{eqn:gGOE}
\end{equation}
for $\x=(x_1, x_2, \dots, x_N) \in \WA$ \cite{Mehta,N05}. 
Here $C_2(N)$ is the same constant as $C_2(N)$ in (\ref{eqn:asym}).
By applying the generalized Bru's theorem to $\Xi^{\rm GOE}(t)$,
we see that $\vlambda(t), t \in (0, \infty)$ solves
the system of
stochastic differential equations of
Dyson's Brownian motion model (\ref{eqn:Dyson0})  
with $\beta=1$.

%%%%%%%%%%%%%%%%%%%%%%%%%%%%%%%%%%%%%%%%%%%%%%%%%%%%%%%%%%%
\vskip 0.3cm
\noindent{\bf (iii)} \quad {\bf Laguerre process}

Let $\nu \in \N_0$.
We denote the space of
$(N+\nu) \times N$ complex matrices 
by $\mathcal{M}(N+\nu,N; \C)$.
Consider the
$\mathcal{M}(N+\nu,N ;\C)$-valued process defined by
$L(t)=(B_{ij}(t)
+\sqrt{-1}\widetilde{B}_{ij}(t))_{1\leq i\leq N+\nu,1\leq j\leq N}$.
For any fixed $t\in [0,\infty)$, 
$L(t)$ is the $\mathcal{M}(N+\nu,N; \C)$-valued random variable
whose density is given by
$$
\mu_\nu^{\rm chGUE}(L , t)
= \frac{t^{-N(N+\nu)}}{c_3(N)}
\exp \left( - \frac{1}{2t} {\rm Tr}L^{*} L \right), 
\quad L \in \mathcal{M}(N+\nu,N; \C),
$$
where $c_3(N)= (2\pi)^{N(N+\nu)}$.
Such a statistical ensemble is called
the {\bf chiral Gaussian unitary ensemble, chGUE} \cite{N05}. 
And the $\H$-valued process defined by 
$\Xi^{\rm L}(t)= L(t)^{*} L(t)$, 
$t \in [0, \infty)$
is called the {\bf Laguerre process} \cite{KO01}.
The matrix $\Xi^{\rm L}(t)$ is nonnegative definite
and has only nonnegative eigenvalues.
Applying the generalized Bru's theorem,
we see that the eigenvalue process $\vlambda(t)$ solves
the following system of stochastic differential equations:
$$
\lambda_{i}(t)= 2 \int_{0}^{t} \sqrt{\lambda_{i}(s)} dB_{i}(s)
+ 2(N+\nu)t
+2 \sum_{\substack{1 \leq j \leq N\\ j \ne i}}
\int_{0}^{t}
\frac{\lambda_{i}(s)+\lambda_{j}(s)}{\lambda_{i}(s)-\lambda_{j}(s)}ds,
\quad 1 \leq i \leq N,
$$
and $\tau=\infty$.
Moreover, by Ito's formula 
we can prove that 
$\vkappa(t)=(\kappa_1(t),\dots, \kappa_N(t)) \equiv
(\sqrt{\lambda_1(t)},\dots, \sqrt{\lambda_N(t)})$
solves the system of stochastic differential 
equations (\ref{eqn:Ynu}),
which implies the equivalence 
of $\vkappa(t)$ with the noncolliding 
$2(\nu+1)$-dimensional Bessel process $\Y^{(\nu)}(t)$.
Since $\nu$ is nonnegative integer
for the chGUE and the Laguerre process, the dimension 
$2(\nu+1)$ of the corresponding noncolliding Bessel process
is positive and even \cite{KO01}.

%%%%%%%%%%%%%%%%%%%%%%%%%%%%%%%%%%%%%%%%%%%%%%%%%%%%%%%%
\vskip 0.3cm
\noindent{\bf (iv)} \quad {\bf Wishart process}

Let $\nu \in \N_0$ and denote 
the space of $(N+\nu) \times N$ real matrices 
by $\mathcal{M}(N+\nu,N; \R)$.
Consider the $\mathcal{M}(N+\nu,N; \R)$-valued process
$W(t)=(B_{ij}(t))_{1 \leq i \leq N+\nu, 1 \leq j \leq N}$.
For any fixed $t\in [0,\infty)$, $W(t)$ is the
$\mathcal{M}(N+\nu,N; \R)$-valued random variable, 
whose probability density function is given by
\begin{equation*}
\mu_\nu^{\rm chGOE}(W , t)
= 
\frac{t^{-N(N+\nu)/2}}{c_4(N)}
\exp \left( - \frac{1}{2t} {\rm Tr}\ ^{t}W W \right),
\quad W \in \mathcal{M}(N+\nu,N; \R), 
\end{equation*}
where $c_4(N)= (2\pi)^{N(N+\nu)/2}$.
Such a statistical ensemble is called 
the {\bf chiral Gaussian orthogonal ensemble, chGOE} \cite{N05}. 
The $\S$-valued process defined by
$\Xi^{\rm W}(t)= {^{t}W(t)} W(t)$, $t \in [0, \infty)$
is called the {\bf Wishart process} \cite{Bru91}.
We can show that the eigenvalue process $\vlambda(t)$ of the 
Wishart process $\Xi^{\rm W}(t)$ solves the following system of
stochastic differential equations
\begin{equation*}
\lambda_{i}(t)= 2 \int_{0}^{t} \sqrt{\lambda_{i}(s)} dB_{i}(s)
+ (N+\nu)t
+\sum_{\substack{1 \leq j \leq N \\ j \not=i}} 
\int_{0}^{t}
\frac{\lambda_{i}(s)+\lambda_{j}(s)}
{\lambda_{i}(s)-\lambda_{j}(s)} ds,
\quad 1 \leq i \leq N
\end{equation*}
and $\tau=\infty$. 

\vskip 3mm

We are able to apply the generalized Bru's theorem to 
Hermitian-matrix-valued processes with additional symmetries.
We introduce the matrix $\sigma_{0}$ 
and the Pauli spin matrices $\sigma_{i}, i=1,2,3$ by
\begin{eqnarray}
&& \sigma_{0}=I_{2}= \left(
\begin{matrix}1 & 0 \cr 0 & 1\end{matrix} 
\right),
\;
\sigma_{1}=\left(
\begin{matrix} 0 & 1 \cr 1 & 0 \end{matrix}
\right), \;
\sigma_{2}=\left(
\begin{matrix} 0 & -\sqrt{-1} \cr \sqrt{-1} & 0 \end{matrix} 
\right), \;
\sigma_{3}=\left(
\begin{matrix} 1 & 0 \cr 0 & -1 \end{matrix} 
\right).
\nonumber
\end{eqnarray}
Suppose that $N \geq 2$ and define $2N \times 2N$ matrices
$\Sigma_{\rho}=I_{N} \otimes \sigma_{\rho}$,
$\rho=0,1,2,3$. By definition $\Sigma_{0}=I_{2N}$. 
Let
$s^{\rho}(t)=(s_{ij}^{\rho}(t))_{1 \leq i, j \leq N}$,
$a^{\rho}(t)=(a_{ij}^{\rho}(t))_{1 \leq i, j \leq N}$,
$0\leq \rho \leq 3$ be 
independent copies of $s(t)$, $a(t)$.
By using them the 
${\mathcal{H}}(2N)$-process $\Xi(t)$ is divided into 
$2 \times 4=8$ terms:
$$
\Xi(t)= \sum_{\rho=0}^{3} \Bigg\{
(s^{\rho}(t) \otimes \sigma_{\rho})+
(\sqrt{-1} a^{\rho}(t) \otimes \sigma_{\rho}) \Bigg\}.
$$
We consider the four Hermitian-matrix-valued
processes represented by the following four terms, 
$$
\Xi_{\theta \varepsilon}(t)=\sum_{\rho=0}^{3}
( \xi_{\theta \varepsilon}^{\rho}(t) \otimes \sigma_{\rho} ), \qquad
\theta=1,2, \quad \varepsilon=\pm,
$$
where
\begin{eqnarray}
&& 
(\xi_{\theta +}^{\rho}(t))=
\left\{
\begin{array}{ll}
s^{\rho}(t),
& \ \mbox{if }\theta =1, \rho \not=3 \mbox{ or }\theta=2,\rho=0,
\\
\sqrt{-1} a^{\rho}(t), 
& \ \mbox{if } \theta =1, \rho=3 \mbox{ or } \theta=2, \rho\not=0,
\end{array}\right.
\nonumber
\\
&& 
(\xi_{\theta -}^{\rho}(t))=
\left\{
\begin{array}{ll}
\sqrt{-1} a^{\rho}(t), 
	& \ \mbox{if } \theta = 1, \rho \not=3 
	\mbox{ or }\theta=2, \rho=0,
\\
s^{\rho}(t), 
	& \ \mbox{if } \theta = 1, \rho = 3 
	\mbox{ or } \theta=2, \rho\not=0 .
   \end{array}\right.
\nonumber
\end{eqnarray}
Putting
$\mathcal{H}_{\theta \varepsilon}(2N)= 
\{H\in \mathcal{H}: ^{t}\!\!\Sigma_\theta =\varepsilon \Sigma_\theta H \}$, 
we see that
$\Xi_{\theta \varepsilon}(t)$ takes values in
$\mathcal{H}_{\theta \varepsilon}(2N)$.
Due to the symmetries of matrices, 
eigenvalues have the following properties;

\noindent{\rm (i)} \;
When $\varepsilon=+$, they are pairwise degenerated;
$
\vlambda=(\omega_{1}, \omega_{1}, \omega_{2}, \omega_{2},
\dots, \omega_{N}, \omega_{N}).
$

\noindent{\rm (ii)} \;
When $\varepsilon=-$, they are in the form
$
\vlambda=(\omega_{1}, -\omega_{1}, \omega_{2}, -\omega_{2},
\dots, \omega_{N}, -\omega_{N}).
$

\vskip 0.3cm
%%%%%%%%%%%%%%%%%%%%%%%%%%%%%%%%%%%%%%%%%%%%%%%%%%%%%%%%%%%%%%%%
\noindent{\bf (v)} \quad {\bf GSE process}

In the case of $(\theta, \varepsilon)=(2, +)$,
a matrix $\Xi \in \mathcal{H}_{2+}(2N)$ is said to be 
a {\bf self-dual} Hermitian matrix,
if it has the symmetry  $^{t}\Xi \Sigma_{2}=\Sigma_{2} \Xi$ 
in addition to Hermitian property.
The matrix can be diagonalized by a unitary-symplectic matrix.
For any fixed $t\in (0,\infty)$ 
the statistical ensemble of $\Xi_{2+}(t)$ is invariant under
any unitary-symplectic transformation, and is called
{\bf Gaussian symplectic ensemble, GSE}.
The probability density function of eigenvalues of GSE is given by
\cite{Mehta,N05},
$$
g^{\rm GSE}(\x; t)=
\frac{t^{-N/2}}{C_3(N)}
h_{N}\left( \frac{\x}{\sqrt{t}} \right)^{4}
\exp \left( - \frac{|\x|^2}{2t} \right), \quad
\x \in {\bf W}_{N}^{A},
$$
where
$C_3(N)=(2\pi)^{N/2}\prod_{i=1}^N \Gamma (2i)$. 
The eigenvalues are pairwise degenerated
and represented as
$
\vlambda=(\omega_{1}, \omega_{1}, \omega_{2}, \omega_{2},
\dots, \omega_{N}, \omega_{N}).
$
This is known as the Kramers doublets
in quantum mechanics. 
Applying the generalized Bru's theorem,
we see that the distinct eigenvalues $\omega_i, 1 \leq i \leq N$ 
solves the system of equations
of Dyson's Brownian motion model (\ref{eqn:Dyson0}) with $\beta =4$, 
$$
\omega_{i}(t) = B_{i}(t) +
\sum_{\substack{1 \leq j \leq N \\ j \not=i} }
\int_0^t \frac{2}{\omega_{i}(s)-\omega_{j}(s)} ds, \quad
1 \leq i \leq N.
$$
For a pair of degenerated eigenvalues 
$\lambda_{2i-1}=\lambda_{2i}=\omega_i$  ,
$\sigma_{2i-1,2i}=\infty$, $1\le i \le N$.
All other pairs
separately move and never coincide with each other, that is,
$\tau=\infty$.

\vskip 0.3cm
%%%%%%%%%%%%%%%%%%%%%%%%%%%%%%%%%%%%%%%%%%%%%%%%%%%%%%%%%%%%%%%%%%%
\noindent{\bf (vi)} \quad {\bf Matrix-valued process of class C }

In the case of $(\theta, \varepsilon)=(2, -)$,
a matrix $\Xi \in \mathcal{H}_{2-}(2N)$ has
the symmetry 
$^{t}\Xi(t) \Sigma_{2}= - \Sigma_{2} \Xi(t)$ 
in addition to Hermitian property.
We denote by ${\rm sp}(2N; \C)$
the Lie algebra of complex symplectic group 
represented by $2N \times 2N$ matrices. Then 
$\mathcal{H}_{2-}(2N) \simeq {\rm sp}(2N; \C) \cap \mathcal{H}(2N)$.
For fixed $t\in (0,\infty)$
the statistical ensemble of $\Xi_{2-}(t)$
coincides with the random matrix ensemble called {\bf class C}
introduced by Altland and Zirnbauer \cite{AZ97}.
The eigenvalues of a matrix in this class
are in the form
$
\vlambda=(\omega_{1}, -\omega_{1}, \omega_{2}, -\omega_{2},
\dots, \omega_{N}, -\omega_{N}).
$
We denote the increasing sequence of nonnegative eigenvalues
by $\vomega(t)=(\omega_1(t),\omega_2(t),\dots,\omega_N(t))$.
By the generalized Bru's theorem, 
we see that $\vomega(t)$ solves the following
system of equations
\begin{eqnarray}
&&\omega_{i}(t) = B_{i}(t) 
+ \int_0^t \frac{1}{\omega_{i}(s)}ds
+ \sum_{\substack{ 1 \leq j \leq N \\ j \not=i} }
\int_0^t \left\{ \frac{1}{\omega_{i}(s)-\omega_{j}(s)}
+  \frac{1}{\omega_{i}(s)+\omega_{j}(s)} \right\}ds, 
\nonumber\\
&&\qquad\qquad\qquad\qquad\qquad\qquad\qquad\qquad\qquad\qquad
\qquad\qquad\qquad
1 \leq i \leq N.
\nonumber
\end{eqnarray}
It is also verified that
$\tau=\infty$ and 
$\vomega(t) \in {\bf W}^{\rm C}_{N}$, $
^{\forall} t \in (0, \infty)$ with probability one.
Comparing the above system of equations with (\ref{eqn:Ynu}),
we can conclude the equivalence
of $\vomega(t)$ and the noncolliding three-dimensional
Bessel process $\Y^{(1/2)}(t)$, in distribution \cite{KT04}.
Moreover, it coincides with the noncolliding Brownian motion
under the condition that it never hits the wall at the origin
\cite{KT04}.

%%%%%%%%%%%%%%%%%%%%%%%%%%%%%%%%%%%%%%%%%%%%%%%%%%%%%%%%%%
\vskip 0.3cm
%%%%%%%%%%%%%%%%%%%%%%%%%%%%%%%%%%%%%%%%%%%%%%%%%%%%%%%%%%%%%
\noindent{\bf (vii)} \quad {\bf Matrix-valued process of class D}

In the case of $(\theta, \varepsilon)=(1, -)$,
a matrix $\Xi \in \mathcal{H}_{1-}(2N)$ has
the symmetry 
$^{t}\Xi(t) \Sigma_{1} = - \Sigma_{1} \Xi(t)$
in addition to Hermitian property.
We denote by ${\rm so}(2N; \C)$
the complexification of Lie algebra of special orthogonal group
represented by $2N \times 2N$ matrices.
Then 
$\mathcal{H}_{1-}(2N) \simeq {\rm so}(2N; \C) \cap \mathcal{H}(2N)$. 
For any fixed $t\in (0,\infty)$
the statistical ensemble of $\Xi_{1-}(t)$ 
coincides with the random matrix ensemble called {\bf class D}
introduced by Altland and Zirnbauer \cite{AZ97}.
The eigenvalues of a matrix in this class
are also in the form
$
\vlambda=(\omega_{1}, -\omega_{1}, \omega_{2}, -\omega_{2},
\dots, \omega_{N}, -\omega_{N}).
$
We denote the increasing sequence of nonnegative eigenvalues
by $\vomega(t)=(\omega_1(t),\omega_2(t),\dots,\omega_N(t))$.
By the generalized Bru's theorem, 
we see that $\vomega(t)$ solves the following
system of equations
$$
\omega_{i}(t) = B_{i}(t) +
\sum_{\substack{ 1 \leq j \leq N \\ j \not=i} }
\int_0^t \left\{ \frac{1}{\omega_{i}(s)-\omega_{j}(s)}
+  \frac{1}{\omega_{i}(s)+\omega_{j}(s)}\right\} ds, \quad
1 \leq i \leq N.
\label{eqn:eq1-}
$$
It is also verified that
$\tau=\infty$ and $\vomega(t) \in {\bf W}^{\rm D}_{N}$, 
$^{\forall} t \in (0, \infty)$ with probability one.
Comparing the above system of equations with (\ref{eqn:Ynu}),
we can conclude that
$\vomega(t)$ 
is equivalent in distribution with
the noncolliding one-dimensional Bessel 
process $\Y^{(-1/2)}(t)$ \cite{KT04}.
Since one-dimensional Bessel process is identified with
a reflecting Brownian motion,
$\vomega(t)$ can be also regarded as the noncolliding 
reflecting Brownian motion
\cite{KT04}.

%%%%%%%%%%%%%%%%%%%%%%%%%%%%%%%%%%%%%%%%%%%%%%%%%%%%%%%%%%%%
%%%%%%%%%%Section 5%%%%%%%%%%%%%%%%%%%%%%%%%%%%%%%%%%%%%%%%%
%%%%%%%%%%%%%%%%%%%%%%%%%%%%%%%%%%%%%%%%%%%%%%%%%%%%%%%%%%%%
\section{Determinantal processes}
%%%%%%%%%%%%%%%%%%%%%%%%%%%%%%%%%%%%%%%%%%%%%%%%%%%%%%%%%%%%
%%%%%%%%%%%%%%%%%%%%%%%%%%%%%%%%%%%%%%%%%%%%%%%%%%%%%%%%%%%%
%%%%%%%%%%%%%%%%%%%%%%%%%%%%%%%%%%%%%%%%%%%%%%%%%%%%%%%%%%%%

%%%%%%%%%%%%%%%%%%%%%%%%%%%%%%%%%%%%%%%%%%%%%%%%%%
\subsection{Fredholm determinant}
%%%%%%%%%%%%%%%%%%%%%%%%%%%%%%%%%%%%%%%%%%%%%%%%%%

Let $\mathcal{X}$ be the space of countable subsets of $\R$ 
without accumulation points.
For $\x=(x_1,x_2,\dots, x_N)\in \R^N$ 
we write $\{ \x \}$ for an element $\{ x_1, x_2, \dots, x_N \}$
of $\mathcal{X}$.
For $\x_N \in \R^N$ and
$N' \in \{ 1,2,\dots,N\}$,
we write $\x_{N'}$ for 
$(x_1, x_2, \dots, x_{N'})\in \R^{N'}$.
For the temporally homogeneous noncolliding Brownian motion $\Y(t)$,
the $\mathcal{X}$-valued process $\xi^N (t) = \{ \Y(t) \}$ has
the transition probability density function
$$
\widetilde{p}_{N}(s, \{\x\} ; t, \{\y\})=
\left\{
\begin{array}{ll}
p_{N}(t-s, \y |\x),
& \mbox{if } s>0,  \ \x, \y \in \WA, 
\\
p_{N}(t, \y|{\bf 0}),
& \mbox{if } s=0, \ \x ={\bf 0}, \ \y \in \WA, 
\\
0, 
& \mbox{otherwise.}
\end{array}\right.
\nonumber
$$
We call the process $\xi^N (t)$ also the 
temporally homogeneous noncolliding Brownian motion
in this paper.
The $\mathcal{X}$-valued noncolliding Bessel process 
$\xi^{N,\nu}(t)$ can be defined as well.
For a sequence of time $0 < t_{1} < \cdots < t_{M} = T$ 
and a sequence of positive integers less than or equal to $N$, 
$\{N_m \}_{m=1}^{M}$, 
the {\bf multi-time correlation function} of
$\xi^N(\cdot)$ at
$(t_m, \{ \x^{(m)}_{N_{m}}\}), m=1,2, \dots, M$, 
is given by
\begin{eqnarray}
&& \cR \left(t_{1},\x^{(1)}_{N_1}; t_2,\x^{(2)}_{N_2}; 
\dots; t_{M},\x^{(M)}_{N_{M}} \right) 
\nonumber\\
&& \, =
\int  \limits _{\prod_{m=1}^{M} \R^{N-N_{m}}}
\prod_{m=1}^{M}\frac{1}{(N-N_{m})!}
\prod_{i=N_{m}+1}^{N} dx_{i}^{(m)}
\prod_{\ell=0}^{M-1} 
\widetilde{p}_{N}
(t_{\ell}, \{\x^{(\ell)}_N\}; t_{\ell+1}, \{\x^{(\ell+1)}_N\}),
\nonumber
\end{eqnarray}
where we put $t_0 = 0$, $\x_N^{(0)} = {\bf 0}$.
Let $C_{0}(\R)$ be the set of all real continuous functions
with compact supports.
For $\f=(f_{1}, f_{2}, \dots, f_{M}) \in C_{0}(\R)^{M}$,
$\vtheta=(\theta_{1}, \theta_{2}, \dots, \theta_{M})\in \R^{M}$,
we put
$\chi_{m}(x)=e^{\theta_{m} f_{m}(x)}-1$, $1 \leq m \leq M$
and $\vchi =(\chi_1, \chi_2, \dots, \chi_M)$.
The {\bf multi-time moment generating function} 
$$
\Psi_{N} (\vchi; \vtheta) \equiv 
E \left[\exp \left\{ \sum_{m=1}^{M}
\theta_{m} \sum_{i_{m}=1}^{N} 
f_{m}(X_{i_{m}}(t_{m})) \right\} \right]
$$
of
$\xi^{N}(t), t \in [0,T]$ 
can be expanded by means of the multi-time correlation functions
as follows:
\begin{eqnarray}
&& \sum_{N_{1}=0}^{N} \sum_{N_{2}=0}^{N} \cdots
\sum_{N_{M}=0}^{N}
\prod_{m=1}^{M}\frac{1}{N_m !}
\int_{\R^{N_{1}}} \prod_{i=1}^{N_1} d x_{i}^{(1)}
\int_{\R^{N_{2}}} \prod_{i=1}^{N_2} d x_{i}^{(2)}
 \cdots
\int_{\R^{N_{M}}} 
\prod_{i=1}^{N_{M}} d x_{i}^{(M)} \nonumber\\
&& \times \prod_{m=1}^{M} \prod_{i=1}^{N_{m}} 
\chi_{m} \Big(x_{i}^{(m)} \Big)
\rho_{N} \left(t_1, \x_{N_{1}}^{(1)} ; 
t_{2}, \x_{N_{2}}^{(2)} ;
\dots ; t_{M}, \x_{N_{M}}^{(M)} \right).
\nonumber
\end{eqnarray}

An $\mathcal{X}$-valued process $\xi(t)$ is called
a {\bf determinantal process}, if 
its multi-time moment generating function is written 
by {\bf Fredholm determinant} as
\begin{equation}
\Psi( \vchi; \vtheta)=
\Det \Big[ \delta_{m,n} \delta (x-y)
+{\bf K}(t_m, x; t_n, y)
\chi_n (y) \Big]
\label{eqn:Det3}
\end{equation}
with a locally integrable function ${\bf K}$.
We call the function ${\bf K}$ the correlation kernel of the process.
By definition of Fredholm determinant, 
the multi-time correlation function of $\xi(t)$
is then given by
$$
\rho_{N} 
\left(t_1, \x^{(1)}_{N_1}; t_2, \x^{(2)}_{N_2}; \dots; 
t_{M}, \x^{(M)}_{N_{M}} \right) 
=\det \left[
\mbA \Big( \x_{N_{0}}^{(0)}, \x_{N_{1}}^{(1)}, \dots,
\x_{N_{M}}^{(M)} \Big) \right]
$$
with a $\sum_{m=1}^{M} N_{m} \times \sum_{m=1}^{M} N_{m}$ matrix
$$
\mbA \Big( \x_{N_{1}}^{(1)}, \x_{N_{2}}^{(2)}, \dots,
\x_{N_{M}}^{(M)} \Big)
= \left( \mbK(t_m, x_{i}^{(m)}; t_n, x_{j}^{(n)} )
\right)_{1 \leq i \leq N_{m}, 1 \leq j \leq N_{n}
1 \leq m, n \leq M}.
$$

The noncolliding Brownian motion $\xi^N(t)$ 
is the determinantal process with the correlation kernel $\mbK_{N}$:
$$
\mbK_N(s,x;t,y) = \left\{
   \begin{array}{ll}
\displaystyle{
\frac{1}{\sqrt{2s}}
\sum_{n=0}^{N-1}
\left(\frac{t}{s}\right)^{n/2}
\varphi_{n}\left(\frac{x}{\sqrt{2s}}\right)
\varphi_{n}\left(\frac{y}{\sqrt{2t}}\right)},
& \quad \mbox{if $s \leq t$,}
\\
\displaystyle{-\frac{1}{\sqrt{2s}}
\sum_{n=N}^{\infty} 
\left(\frac{t}{s}\right)^{n/2}
\varphi_{n}\left(\frac{x}{\sqrt{2s}}\right)
\varphi_{n}\left(\frac{y}{\sqrt{2t}}\right)},
& \quad \mbox{if $s > t$,}
   \end{array} \right. 
$$
where $\varphi_n(x)= 
\{\sqrt{\pi} 2^{n} n!\}^{-1/2} H_n (x)e^{-x^2/2},
n=0,1,2, \dots,$ 
are the orthonormal functions on $\R$  
associated with the Hermite polynomials $H_n(x)$ \cite{KNT04,KT07b}.
The noncolliding Bessel process $\xi^{N,\nu}(t)$
is the determinantal process with the correlation kernel 
${\bf K}^{(\nu)}_{N}$:
$$
{\bf K}^{(\nu)}_{N}(s,x;t,y)= \left\{
   \begin{array}{ll}
\displaystyle{
\frac{\sqrt{xy}}{s}\sum_{n=0}^{N-1}
\left(\frac{t}{s}\right)^n
\varphi_n^\nu \left(\frac{x^2}{2s}\right)
\varphi_n^\nu \left(\frac{y^2}{2t}\right),}
& \quad \mbox{if $s \leq t$,}
\\
\displaystyle{
-\frac{\sqrt{xy}}{s}\sum_{n=N}^{\infty}
\left(\frac{t}{s}\right)^n
\varphi_n^\nu \left(\frac{x^2}{2s}\right)
\varphi_n^\nu \left(\frac{y^2}{2t}\right),}
& \quad \mbox{if $s > t$,}
\end{array} \right. 
$$
where $\varphi_n^\nu (x)
= \sqrt{\Gamma(n+1)/\Gamma(\nu+n+1)}
x^{\nu/2}L_n^{\nu}(x)e^{-x/2},
n=0,1,2, \dots,$ 
are the orthonormal functions on
$\R_{+}=\{x \in \R : x \geq 0 \}$
associated with the Laguerre polynomials $L_n^{\nu}(x)$
with parameter $\nu>-1$ \cite{KT07a}.
See Table 1, which summarizes the correspondence
between noncolliding diffusion processes with
finite number of particles,
statistical ensembles of random matrices (RM),
and orthogonal polynomials used to represent
correlation kernels.

%%%%%%%%%%%%%%%%%%%%%%%%%%%%%%%%%%%%%%%%%%%%%%%%%%%%
%%%%%%%%%%% Table 1%%%%%%%%%%%%%%%%%%%%%%%%%%%%%%%%%%%%%%%%%
%%%%%%%%%%%%%%%%%%%%%%%%%%%%%%%%%%%%%%%%%%%%%%%%%%%%
\begin{table}
\begin{center}
\begin{small}
\begin{tabular}{|c|c||c|c||c|}
\hline
diffusions &Wyle chambers &matrix-valued pr. 
&RM &orth. poly.
\cr\hline
\hline
Brownian motion & $A_{N-1}$  & GUE  &GUE & $H_n$
\cr
& & $\mathcal{H}(N)$& &
\cr\hline
\hline
even-dim. Bessel pr. & $C_{N}$  & Laguerre &chGUE 
& $L_n^{\nu}$, $\nu \in \N_0$
\cr\hline
3-dim. Bessel pr. & $C_{N}$  & class C &class C & $L_{n}^{1/2}$
\cr
 (absorbing BM) & & ${\rm sp}(2N; {\bf C}) \cap \mathcal{H}(2N)$
& &
\cr\hline
1-dim. Bessel pr. & $D_N$  &  class D  &class D & $L_{n}^{-1/2}$
\cr
(reflecting BM)& &${\rm so}(2N; {\bf C}) \cap \mathcal{H}(2N)$ & &
\cr\hline 
\end{tabular}
\end{small}
\end{center}
\caption{Noncolliding diffusion processes, 
random matrix ensembles, and orthogonal polynomials}
\end{table}
%%%%%%%%%%%%%%%%%%%%%%%%%%%%%%%%%%%%%%%%%%%%%%%%%%%%
%%%%%%%%%%% Table 1%%%%%%%%%%%%%%%%%%%%%%%%%%%%%%%%%%%%%%%%%
%%%%%%%%%%%%%%%%%%%%%%%%%%%%%%%%%%%%%%%%%%%%%%%%%%%%
%\hskip 1cm

\subsection{Scaling limits}

When the number of the diffusion process $N$ goes to infinity,
the asymptotic behaviors of the noncolliding processes
$\xi^N$ and $\xi^{N,\nu}$ are determined by 
the asymptotics of their correlation kernels
${\bf K}_N, {\bf K}_{N}^{(\nu)}$ in $N \to \infty$. 
Suppose that the correlation kernel converges
under an appropriate scaling limit.
Then the multi-time moment generating function 
$\Psi_{N} (\vchi; \vtheta)$ and
the multi-time correlation functions 
$\rho_N$ converge,
and then the process converges
in the sense of finite dimensional distributions.
In the following, we discuss the {\bf bulk scaling limit}
and the {\bf soft-edge scaling limit}
for the noncolliding Brownian motion $\Y(t)$
and the {\bf hard-edge scaling limit}
for the noncolliding Bessel process $\Y^{(\nu)}(t)$ 
\cite{NF98,PS02,Joh03,TW04,KNT04,AM05,KT07a,KT07b}.

\begin{enumerate}
\item{{\bf [Bulk scaling limit]}} \quad
As $N\to\infty$,
the process $\xi^N(N+t)$ converges to
the infinite-dimensional determinantal process, 
whose correlation kernel $\mathcal{K}^{{\rm sin}}$
is expressed by using trigonometric functions, 
$$
\mathcal{K}^{{\rm sin}}(s,x;t,y) = \left\{
   \begin{array}{ll}
\displaystyle{
\frac{1}{\pi} \int_{0}^{1} d u \,
e^{(t-s)u^2/2} \cos (u(x-y)),
}
& \mbox{if } s < t,
\\
& \\
%& \\
\displaystyle{
\frac{\sin(x-y)}{\pi(x-y)},
}
& \mbox{if } s=t,
\\
& \\
\displaystyle{
-\frac{1}{\pi} \int_{1}^{\infty} d u \,
e^{(t-s)u^2/2} \cos (u(x-y)),
}
& \mbox{if } s > t.
   \end{array} \right. 
$$

\item{{\bf [Soft-edge scaling limit]}} \quad
Define the scaled process
$\theta_{a(N,t)}\xi^N(N^{1/3}+t)
\equiv \{Y_1(N^{1/3}+t)-a(N,t), Y_2(N^{1/3}+t)-a(N,t),
\dots, Y_N(N^{1/3}+t)-a(N,t)\}$
with $a(N,t)=2N^{2/3}+N^{1/3}t-t^2/4$.
As $N\to\infty$, it
converges to the infinite-dimensional determinantal process, 
whose correlation kernel $\mathcal{K}^{{\rm Ai}}$ 
is expressed by using the Airy function ${\rm Ai}(x)$, 
$$
\mathcal{K}^{{\rm Ai}}(s,x; t,y) = \left\{
   \begin{array}{ll}
\displaystyle{
\int_{-\infty}^{0} du \,
e^{(t-s)u/2} {\rm Ai}(x-u) {\rm Ai}(y-u),
}
&\mbox{if } s \leq t,  \\
& \\
\displaystyle{
-\int_{0}^{\infty} du \,
e^{(t-s)u/2} {\rm Ai}(x-u) {\rm Ai}(y-u),
}
& \mbox{if } s > t. 
   \end{array} \right. 
$$

\item{{\bf [Hard-edge scaling limit]}} \quad
As $N\to\infty$, $\xi^{N,\nu}(N+t)$ 
converges to the infinite-dimensional determinantal process, 
whose correlation kernel $\mathcal{K}^{(\nu)}$ 
is expressed by using the Bessel function $J_{\nu}(x)$, 
$$
\mathcal{K}^{(\nu)}(s, x ; t, y)=\left\{
   \begin{array}{ll}
\displaystyle{
\sqrt{xy}\int_0^2 du \ e^{(t-s)u^2/2}
J_{\nu}(ux)u J_{\nu}(uy),}
 & \mbox{if } s < t,   \\
& \\
\displaystyle{
\frac{2\sqrt{xy}\{
J_{\nu}(2x)y J_{\nu}'(2y)
-J_{\nu}(2y)xJ_{\nu}'(2x)\}}{x^2-y^2 }},
 & \mbox{if } s=t,   \\
& \\
\displaystyle{
-\sqrt{xy}\int_2^{\infty} du \ e^{(t-s)u^2/2}
J_{\nu}(ux)u J_{\nu}(uy),}
 & \mbox{if } s> t. \quad  \\
   \end{array}\right. 
$$
\end{enumerate}

%%%%%%%%%%%%%%%%%%%%%%%%%%%%%%%%%%%%%%%%%%%%%%%%%%%%%%
The above three infinite particle systems 
are all temporally homogeneous.
The system obtained by the bulk scaling limit 
is spatially homogeneous, 
and the other systems
obtained by the soft- and hard-edge scaling limits
are spatially inhomogeneous (see Table 2). 
These infinite particle systems are reversible and
their equilibrium measures are determinantal point processes
\cite{KT07b,KT10}.
Osada \cite{Osa96,Osa04} constructed diffusion processes 
whose equilibrium measures are determinantal point processes, 
by the Dirichlet form technique.
Although the coincidence of Osada's processes and the above 
processes is expected \cite {KT07b}, it has not been proved yet.
If the coincidence were proved, it would be concluded  that
the infinite particle systems obtained by the bulk scaling limit 
and the hard-edge scaling limit solve 
the stochastic differential equations
(\ref{eqn:Dyson0}) and (\ref{eqn:Ynu}) with $N=\infty$ 
(see \cite{Osa10}).
Nonequilibrium dynamics of determinantal processes 
with infinite numbers of particles have been studied, 
which show the relaxation processes to the stationary 
determinantal processes with the correlation kernels 
$\mathcal{K}^{{\rm sin}}$ and $\mathcal{K}^{{\rm Ai}}$  \cite{KT09,KT10}.
There the theory of distributions of
zeros and orders of growth of entire functions \cite{Levin96} 
is applied to analyze the determinantal structures of 
noncolliding diffusion processes with infinite numbers of particles.

%%%%%%%%%%%%%%%%%%%%%%%%%%%%%%%%%%%%%%%%%%%%%%%%%%%
\subsection{Tracy-Widom distribution}
%%%%%%%%%%%%%%%%%%%%%%%%%%%%%%%%%%%%%%%%%%%%%%%%%%%

Consider the motion of the rightmost particle in
the temporally homogeneous noncolliding Brownian motion 
$\xi^N(t)=\{ \Y(t) \}$.
For a fixed time $t>0$, 
from (\ref{eqn:Det3}) with
$M=1, t_1=t, \theta_1=1$, $\chi_1(x)=-{\bf 1}_{\{x > \alpha \}}$
the probability that
the position of the rightmost particle is less than $\alpha\in\R$ 
is given by the Fredholm determinant as
\begin{eqnarray}
P\left(\max_{1\le i \le N}Y_i(t) \leq \alpha\right)
&=&
E \left[\exp\left\{\sum_{i=1}^N 
\log ({\bf 1}_{\{Y_i(t) \leq \alpha\}})\right\}\right]
\nonumber\\
&=&\Det \Big[\delta(x-y)-{\bf K}_N(t, x; t, y)
{\bf 1}_{\{y > \alpha \}} \Big].
\nonumber
\end{eqnarray}
Then the distribution function $F_{\rm max}(\alpha)$ 
of the position of the rightmost particle
in the process obtained by the soft-edge scaling limit 
is given by
$$
F_{\rm max}(\alpha) \equiv 
\Det \Big[\delta(x-y)-\mathcal{K}^{{\rm Ai}}(t, x; t, y)
{\bf 1}_{\{y > \alpha \}} \Big],
$$
where the kernel $\mathcal{K}^{{\rm Ai}}$ was
defined in the previous subsection.
Tracy and Widom \cite{TW94}
represented the distribution function 
({\bf Tracy-Widom distribution}) as 
$$
F_{\rm max}(\alpha)= \exp\left(-\int_\alpha^\infty 
(x-\alpha)q(x)^2 dx \right)
$$
with the solution  $q(x)$ of {\bf Painlev\'{e} II}
(see for instance \cite{NY01}j
$$
P_{\rm II} \quad : \quad
\frac{d^2q(x)}{dx^2}=2q(x)^3+xq(x)
$$
satisfying the boundary condition $q(x)\sim {\rm Ai}(x)$, 
$x\to\infty$.
For determinantal process,
the distribution of position of the
right-nearest particle to the fixed point (for instance the origin)
is described by Fredholm determinant
as well as the rightmost particle.
In the bulk scaling limit \cite{JMMS80} and
the hard-edge scaling limit \cite{TW94B}
the distributions of right-nearest particles 
to fixed points are 
studied precisely, and are represented by solutions of
Painlev\'{e} V ($P_{\rm V}$) and
Painlev\'{e} III ($P_{\rm III}$),
respectively (see Table 2).

%%%%%%%%%%%%%%%%%%%%%%%%%%%%%%%%%%%%%%%%%%%%%%%%%%%%%%%%%%%%%
\hskip 1cm
%%%%%%%%%%%%%%%%

%%%%%%%%%%%%%%%%%%%%%%%%%%%%%%%%%%%%%%%%%%%%%%%%%%%%
%%%%%%%%%%% Table 2%%%%%%%%%%%%%%%%%%%%%%%%%%%%%%%%%%%%%%%%%
%%%%%%%%%%%%%%%%%%%%%%%%%%%%%%%%%%%%%%%%%%%%%%%%%%%%
\begin{table}
\begin{center}
\begin{tabular}{|c||c|c|c|}
\hline
kernels (finite system) & Hermite $H_n$ & Hermite $H_n$ &
Laguerre $L_n^{\nu}$ \cr
\hline
kernels (infinite system) & (bulk) &
(soft-edge ) & (hard-edge) \cr
   & trigonometric & Airy & Bessel \cr
   & $\sin, \cos$ & Ai & $J_{\nu}$ \cr
\hline
spatial homogeneity & homogeneous & inhomogeneous & inhomogeneous \cr
\hline
\hline
Painlev\'{e} equation & $P_{\rm V}$ & $P_{\rm II}$ & $P_{\rm III}$ \cr
\hline
\end{tabular}
\end{center}
\caption{Determinantal processes and Painlev\'{e} equations}
\end{table}
%%%%%%%%%%%%%%%%%%%%%%%%%%%%%%%%%%%%%%%%%%%%%%%%%%%%
%%%%%%%%%%% Table 2%%%%%%%%%%%%%%%%%%%%%%%%%%%%%%%%%%%%%%%%%
%%%%%%%%%%%%%%%%%%%%%%%%%%%%%%%%%%%%%%%%%%%%%%%%%%%%
%%%%%%%%%%%%%%%%%%%%%%%%%%%%%%%%%%%%%%%%%%%%%%%%%%%%%%%%%%%%%%%

%%%%%%%%%%%%%%%%%%%%%%%%%%%%%%%%%%%%%%%%%%%%%%%%%%%%%%%%%%%%
%%%%%%%%%%Section 6%%%%%%%%%%%%%%%%%%%%%%%%%%%%%%%%%%%%%%%%%
%%%%%%%%%%%%%%%%%%%%%%%%%%%%%%%%%%%%%%%%%%%%%%%%%%%%%%%%%%%%
\section{Temporally inhomogeneous processes}
%%%%%%%%%%%%%%%%%%%%%%%%%%%%%%%%%%%%%%%%%%%%%%%%%%%%%%%%%%%%
%%%%%%%%%%%%%%%%%%%%%%%%%%%%%%%%%%%%%%%%%%%%%%%%%%%%%%%%%%%%
%%%%%%%%%%%%%%%%%%%%%%%%%%%%%%%%%%%%%%%%%%%%%%%%%%%%%%%%%%%%

The noncolliding processes discussed in the previous section
are temporally homogeneous diffusion processes,
in which noncolliding conditions are imposed 
in the infinite time-intervals $(0,\infty)$.
On the other hand, in Section 3 we explained that
the noncolliding Brownian motion  $\X(t), t \in [0, T]$ is a
temporally inhomogeneous diffusion process 
with transition probability density function (\ref{eqn:g(x,y)})
with (\ref{eqn:g(0,y)}), 
if noncolliding conditions are imposed during 
a finite time-interval $(0,T]$.
Since $\mathcal{N}_N(0,\y)=1$, $\y \in \WA$ by definition,
the probability density function (\ref{eqn:g(0,y)})
of the process $\X(T)$ coincides 
with $g^{\rm GOE}$ given by (\ref{eqn:gGOE}).
In other words, 
the distribution of the noncolliding Brownian motion 
at the terminal time $T$ of the noncolliding time-interval
is equal to the eigenvalue distribution of GOE.
When $0<t <T$, however, the distribution of $\X(t)$ is
different from the eigenvalue distribution of GOE.
In particular, when $0<t \ll T$, we see that
it is close to the eigenvalue distribution of GUE
by the asymptotic behavior (\ref{eqn:asym}) of $\mathcal{N}_N(t,\y)$.
>From the above observations it is expected that
$\X(t), t \in [0, T]$ exhibits a
{\bf transition from GUE to GOE} as $t$ approaches $T$.

Remind that the off-diagonal elements of the GOE 
process $\Xi^{\rm GOE}(t)$ are
one-dimensional Brownian motions and 
those of GUE process $\Xi^{\rm GUE}(t)$ 
are complex Brownian motions.
Hence, an $\H$-valued process, whose eigenvalue process 
realizes $\X(t), t \in (0, T]$, 
should satisfy the condition that
each off-diagonal element behaves like a complex Brownian motion
for $0<t \ll T$,
and it becomes to behave like 
a real Brownian motion as $t \nearrow T$
\cite{PM83}.
We find that, if each imaginary part of
off-diagonal element is given 
by the Brownian bridge of duration $T$,
which was introduced in Section 2,
this condition is fulfilled.
Let $\beta_{ij}^{T}(t)$, $1 \leq i < j \leq N$ be independent 
Brownian bridges, which are assumed to be 
independent of the Brownian motions
$B_{ij}(t)$, $1 \leq i\leq j \leq N$ 
used in the definition of $s_{ij}(t)$.
Then we put
$$
a_{ij}^{T}(t)
=
\left\{
   \begin{array}{ll}
      \displaystyle{
      \frac{1}{\sqrt{2}} \beta_{ij}^{T}(t),
      } & 
   \mbox{if $i < j$ } \\
%& \\
    0, & \mbox{if $i=j$} \\
%& \\
     \displaystyle{
      -\frac{1}{\sqrt{2}} \beta_{ji}^{T}(t),
      } & 
   \mbox{if $i > j$ } \\
   \end{array}\right.
%\label{eqn:aT}
$$
and define an $\H$-valued process by
\begin{equation}
\Xi^{T}(t)= \Big( s_{ij}(t)+ \sqrt{-1} a_{ij}^{T}(t)
\Big)_{1 \leq i,j \leq N}, \quad
t \in [0,T].
\label{eqn:XiT}
\end{equation}
Then we see that the eigenvalue process $\vlambda^{T}(t)$
of $\Xi^{T}(t), t \in (0,T]$ is a temporally inhomogeneous 
diffusion process and
is equivalent in distribution
with the temporally inhomogeneous 
noncolliding Brownian motion $\X(t), t \in (0, T]$ 
with $\X(0)={\bf 0}$ \cite{KT03b}. 
The equivalence is proved by the fact that
the eigenvalue process $\vlambda^{T}(t)$ of $\Xi^{T}(t)$
and that of $\Xi^{\rm GUE}(t)$ satisfy
the generalized Imhof's relation (\ref{Imhof:BM}).
Remember that
the eigenvalue process of $\Xi^{\rm GUE}(t)$ 
is identified with 
the temporally homogeneous noncolliding Brownian motion $\Y(t)$
in distribution as shown in Section 4.

This result implies that the process
$\X(t), t \in (0, T]$ has two different representations, 
`the representation by a noncolliding Brownian motion with 
the transition probability density
(\ref{eqn:g(0,y)}) given by
the Karlin-McGregor formula',
and `the representation by
an eigenvalue process of 
$\Xi^{T}(t)$ given by (\ref{eqn:XiT})'.
This claim is a generalization of the result that
the three-dimensional Bessel process 
as well as the generalized meander
have two different representations, 
`the representation by 
conditional one-dimensional Brownian motions'
and `the representation by 
radial parts of three-dimensional diffusion processes'.

The $\H$-valued process $\Xi^{T}(t)$ is decomposed into 
an eigenvalue part $\Lambda(t)$ and a unitary matrix part $U(t)$. 
The latter representation of the process $\X^{T}(t), t \in (0, T]$
implies that it is obtained from $\Xi^{T}(t)$ by integrating 
its unitary matrix part $U(t)$.
By this observation 
the following identity is derived \cite{KT03b},
which is called
the {\bf Harish-Chandra integral formula} \cite{HC57}
or the {\bf Itzykson-Zuber integral formula} \cite{IZ80}.

\vskip 3mm

\noindent{\large\bf [Harish-Chandra integral formula]}
\; 
Let $dU$ be the Haar measure of the space ${\bf U}(N)$ normalized as
$\int_{{\bf U}(N)}dU=1$. 
For $\x=(x_{1}, x_{2}, \dots, x_{N}) \in \WA$ and
$\y=(y_{1}, y_{2}, \dots, y_{N}) \in \WA$, put
$\Lambda_{\x}={\rm diag}(x_{1}, x_{2}, \dots, x_{N})$
and 
$\Lambda_{\y}={\rm diag}(y_{1}, y_{2}, \dots, y_{N})$.
Then for any $\sigma \in \R$ the following identity holds:
%\begin{equation*}
$$
\int_{{\bf U}(N)} dU \,
e^{ - {\rm Tr}(\Lambda_{\tiny \x}-U^{*} \Lambda_{\tiny \y} U)^2
/(2 \sigma^{2})}
= \frac{C_{1}(N) \sigma^{N^2}}
{h_{N}(\x) h_{N}(\y)}
\det_{1 \leq i, j \leq N}
\Big( G(\sigma^2, y_j |x_{i}) \Big).
%\label{eqn:HC}
%\end{equation*}
$$

The above argument is also valid for 
the noncolliding generalized meander $\X^{(\nu,\kappa)}(t)$.
The matrix-valued process, whose elements are the complex-valued process
having Brownian motions as its real part 
and Brownian bridges as its imaginary part,
$$
M_T(t)=\Big(B_{ij}(t)
+\sqrt{-1}\beta^T_{ij}(t) \Big)_{1\le i \le N+\nu, 1\le j\le N},
$$
exhibits a {\bf transition from chGUE to chGOE}. 
Its eigenvalue process is a complex-valued process
and different from the noncolliding generalized meander.
In stead of the process $M_{T}(t)$, 
we consider the $\H$-valued process defined by
$$
\Xi_T^{\rm LW}(t)=M_T(t)^{*} M_T(t), \quad t\in [0,T]
$$
The eigenvalue process
$\vlambda^{\rm LW} (t)
=(\lambda_1^{\rm LW}(t),\lambda_2^{\rm LW}(t),
\dots,\lambda_N^{\rm LW}(t))$
of $\Xi_T^{\rm LW}(t)$
is the stochastic process with $N$ nonnegative coordinates.
We put 
$\kappa^{\rm LW}_i(t)=\sqrt{\lambda_i^{\rm LW}(t)}$, $1\le i \le N$
and consider 
$\vkappa^{\rm LW}(t)=(\kappa_1^{\rm LW}(t), \kappa_2^{\rm LW}(t), 
\dots, \kappa_N^{\rm LW}(t))$.
Then the process $\vkappa^{\rm LW}(t)$ 
is temporally inhomogeneous and
is identified with the noncolliding generalized meander 
$\X^{(\nu,\kappa)}(t), t \in (0, T]$ 
with $\nu \in \N_0, \kappa=\nu+1$ and 
$\X^{(\nu, \kappa)}(0)={\bf 0}$ \cite{KT04}. 
As another example, by setting $(\nu,\kappa)=(1/2,1)$,
we can construct an $\H$-valued process, 
whose eigenvalue process is a noncolliding generalized meander
exhibiting a {\bf transition from class C to class CI}.
(For the definitions of the random matrix ensembles
called class CI and class DIII-odd/even mentioned below, see
\cite{Zir96,AZ97,Iva01,CM04}.)

For a system of $2N$ independent Brownian motions, 
we impose the condition that
pairs of $(2i-1)$-th and $2i$-th particles
meet at the terminal time $T$ for $1 \leq i \leq N$, 
in addition to the noncolliding condition in the
time-interval $(0,T)$.
Then we can show that the system realizes
the eigenvalue process of the matrix-valued process, 
which exhibits a {\bf transition from GUE to GSE}.

Consider the noncolliding generalized meanders with
$(\nu, \kappa)=(\nu, \nu+1), \nu \in \N_0$, 
$(\nu, \kappa)=(1/2, 0)$ and $(\nu, \kappa)=(-1/2, 0)$,
with the above mentioned additional condition at $t=T$.
We can prove that they realize the eigenvalue processes
of the matrix-valued processes,
which shows transitions {\bf from chGUE to chGSE}, 
{\bf from D to class DIII-odd} and
{\bf from class D to class DIII-even}, respectively
\cite{KT04}. 
(Here chGSE indicates the random matrix ensemble called the
{\bf chiral Gaussian symplectic ensemble} \cite{N05}.)

Several temporally inhomogeneous processes 
have two different representations, 
`the representation by noncolliding diffusion processes', and
`the representation by 
eigenvalue processes of matrix-value processes' 
(see Table 3).
By identifying these two representations, 
Harish-Chandra (Itzykson-Zuber) formulas
are derived for matrices with a variety of 
symmetries \cite{KT04}.

Recently, a family of stochastic processes, 
whose multi-time moment generating functions are
represented by Fredholm Pfaffians \cite{Rai00,KT07a},
has been intensively studied \cite{Sos04,BR05}.
We call such a stochastic process 
a {\bf Pfaffian process}.
In general, the multi-time $N$ point correlation 
function of a Pfaffian process is described by
a Pfaffian of $2N\times 2N$ matrix (see Table 3). 
Since Pfaffians of $2N\times 2N$ matrices
are reduced to determinants of $N\times N$ matrices
in special cases,
Pfaffian process is regarded as a generalization
of determinantal process.
Dyson introduced stochastic processes, 
whose $N$ point correlation functions
are represented by $N\times N$ quaternion determinants,
and they have been studied since then \cite{Dys70,FNH99,N05}.
These processes are also members of Pfaffian processes, because
quaternion determinants can be expressed by Pfaffians.
We showed that the temporally inhomogeneous version of
noncolliding Brownian motion $\X(t)$
and the noncolliding generalized meanders $\X^{(\nu, \kappa)}$ 
are Pfaffian processes \cite{KNT04, KT07a}.
By evaluating the asymptotics of Pfaffians in $N\to\infty$,
we can prove the existence of infinite-dimensional Pfaffian 
processes in appropriate scaling limits.
They describe temporally inhomogeneous infinite particle systems.
For the noncolliding generalized meander $\X^{(\nu, \kappa)}$,
the general form of correlation kernel is
described by using 
{\bf Riemann-Liouville differintegrals} \cite{KT07a}.

%%%%%%%%%%%%%%%%%%%%%%%%%%%%%%%%%%%%%%%%%%%%%%%%%%%%%%%
\hskip 1cm
\begin{table}
\begin{center}
\begin{tabular}{|c||c|c|}
\hline
homogeneity & homogeneous & inhomogeneous \cr
\hline
\hline
1 dim. diffusion & Brownian motion & Brownian bridge \cr
& Bessel process & Bessel bridge \cr
& & generalized meander \cr
\hline
\hline
matrix-valued pr. & GUE& GUE-to-GOE, GUE-to-GSE \cr
& chGUE & chGUE-to-chGOE, chGUE-to-chGSE 
\cr
& class C & class C-to-class CI %($\nu=1/2, \kappa=1$) 
\cr
& class D & class D-to-class DIII-odd %($\nu=1/2, \kappa=0$)
\cr
& & class D-to-class DIII-even %($\nu=-1/2, \kappa=0$) 
\cr
\hline
\hline
process & determinantal pr. & Pfaffian process \cr
\hline
corr. func. & determinant & Pfaffian \cr
\hline
moment gen. func. & Fredholm det. &
Fredholm Pfaffian \cr
\hline
\end{tabular}
\end{center}
\caption{Temporally homogeneous and inhomogeneous processes}
\end{table}
%%%%%%%%%%%%%%%%%%%%%%%%%%%%%%%%%%%%%%%%%%%%%%%%%%%%%%%%%%%

%%%%%%%%%%%%%%%%%%%%%%%%%%%%%%%%%%%%%%%%%%%%%%%%%%%%%%%%%%%%
%%%%%%%%%%Section 7%%%%%%%%%%%%%%%%%%%%%%%%%%%%%%%%%%%%%%%%%
%%%%%%%%%%%%%%%%%%%%%%%%%%%%%%%%%%%%%%%%%%%%%%%%%%%%%%%%%%%%
\section{Miscellanea}
%%%%%%%%%%%%%%%%%%%%%%%%%%%%%%%%%%%%%%%%%%%%%%%%%%%%%%%%%%%%
%%%%%%%%%%%%%%%%%%%%%%%%%%%%%%%%%%%%%%%%%%%%%%%%%%%%%%%%%%%%

\noindent {\rm 1.} \quad
In this paper, we have discussed noncolliding systems
of one-dimensional diffusion processes
in the unbounded domains, $\R$ and $\R_{+}$,
which are related to Gaussian ensembles of random matrices.
We can also consider noncolliding systems in bounded domains.
In particular, the systems on a circle have been studied 
and the relation with the statistical ensembles
of random unitary matrices called {\bf circular ensembles}
are reported \cite{Dys62a, HW96, Mehta, NF03}.
For the systems on finite intervals,
the transition probability density functions are described by
using the Jacobi polynomials and 
the systems are related to the random matrix 
model called
{\bf MANOVA (multivariate analysis of variance) model} \cite{DE02}.

\vskip 0.3cm
\noindent {\rm 2.} \quad
Dyson's Brownian motion models, which solve 
the system of equations (\ref{eqn:Dyson0}),
form a family of processes with a parameter $\beta >0$. 
In this paper by applying the generalized Bru's theorem 
we have clarified the correspondence between 
the eigenvalue processes 
associated with the random matrix ensembles GOE, GUE and GSE, 
and the Dyson's Brownian motion models with $\beta=1,2$ and $4$.
In particular, it was shown that, when $\beta=2$, 
the process is also realized by
the noncolliding Brownian motion.
Recently, 
a family of random matrix ensemble with 
a parameter $\beta >0$ is proposed,
in which the eigenvalue distribution is give by
$$
g^{\beta}(\x)=
\frac{1}{C_\beta (N)}
h_{N}(\x)^{\beta}
\exp \left( - \frac{|\x|^2}{2} \right), \quad
\x \in {\bf W}_{N}^{A}, 
$$
where $C_\beta (N)$ is the normalization constant.
This random matrix ensemble is called the {\bf Gaussian beta ensemble},
whose elements are
tridiagonal matrices such that
diagonal elements are independent Gaussian random variables
and 
$(k,k+1)$-elements and $(k+1,k)$-elements, $1\le k \le N-1$, 
are independent random variables with $\chi$-square distribution
with degree-of-freedom $(N-k)\beta$ \cite{DE02}.

\vskip 0.3cm
\noindent {\rm 3.} \quad
The ensembles of random matrices, whose elements are 
independent complex Gaussian random variables, 
is called the {\bf Ginibre ensemble} \cite{G65}.
The eigenvalues of matrices in this ensemble
are complex in general and the probability density function 
is given by
$$
g^{{\rm Gin}}(\z)=
\frac{1}{C_{\rm Gin} (N)}
\prod_{1\le i<j\le N} |z_i-z_j|^2
\exp \left( - \frac{|\z|^2}{2} \right), \quad
\z \in {\bf C}^N,
$$
where $C_{\rm Gin} (N)$ is the normalization constant.
Characterization of $g^{{\rm Gin}}(\z)$
has been intensively studied 
(see for instance \cite{S06,Osa10}).

\vskip 0.3cm
\noindent {\rm 4.} \quad
In the present paper 
noncolliding diffusion processes are discussed.
Noncolliding systems of discrete time Markov processes
have been also studied. 
In particular, the system of noncolliding random walks,
called the {\bf vicious walk model} \cite{Fis84}, 
is an interesting and important model, since
it is related to the representation theory
of symmetry groups
through the Young diagrams, the Young tableaux, and 
the Schur functions \cite{Joh00,Baik00,Joh02,KT03a,
OR03,Katori04a,Hor05,BO06,BO09}. 

\begin{small}
%%%%%%%%%%%%%%%%%%%%%%%%%%%%%%%%%%%%%%%%%%%%%%%%%%%%%%%%%%%%
%%%%%%%%%%BIBLIOGRAPHY%%%%%%%%%%%%%%%%%%%%%%%%%%%%%%%%%%%%%%
%%%%%%%%%%%%%%%%%%%%%%%%%%%%%%%%%%%%%%%%%%%%%%%%%%%%%%%%%%%%
\bibliographystyle{amsplain}

\end{small}

\end{document}